\documentclass[reqno]{amsart}
\usepackage{amsthm,amsmath,amssymb,amscd, mathrsfs
}
\input diagrams
\diagramstyle[PostScript=dvips]
\usepackage[dvips]{graphics}

\newcommand{\ca}{\mathscr{A}} 
\newcommand{\ce}{{\mathscr E}}
\newcommand{\cb}{{\mathscr B}}
\newcommand{\cf}{{\mathscr F}}
\newcommand{\cg}{\mathscr{G}} 
\newcommand{\cl}{{\mathscr L}} 
\newcommand{\co}{{\mathscr O}}
\newcommand{\codim}{\operatorname{codim}} 
 
\newcommand{\crr}{{\mathscr R}}
\newcommand{\crrt}{\widetilde{\crr}} 
\newcommand{\crrfg}{{\mathscr R}_{FG}}
\newcommand{\cc}{{\mathscr C}}
\newcommand{\ch}{{\mathscr H}} 
\newcommand{\ci}{{\mathfrak I}} 
\newcommand{\cs}{{\mathscr S}}
\newcommand{\cso}{\tilde{\mathscr S}} 
\newcommand{\ct}{{\mathscr T}} 
\newcommand{\cv}{{\mathscr V}}
\newcommand{\cx}{\mathscr{X}} 
\newcommand{\cy}{\mathscr{Y}}

\newcommand{\nc}{\mathbb{C}} 
\newcommand{\nq}{\mathbb{Q}}
\newcommand{\nr}{\mathbb{R}} 
\newcommand{\xg}{[X/G]}

\newcommand{\Aut}{\operatorname{Aut}} 

\newcommand{\bk}{\mathbf{k}} 
\newcommand{\bm}{\mathbf{m}}

\newcommand{\bv}{\mathbf{V}}
\newcommand{\cch}{\Ch} 

\newcommand{\cth}{\mathbf{H}} 
 
\newcommand{\cab}{\overline{\mathscr{A}}} 
\newcommand{\Ch}{\mathbf{ch}} 
\newcommand{\CCh}{\mathscr{C}\mathbf{h}} 
\newcommand{\CChf}{\mathsf{C}\mathsf{h}_\mathrm{orb}} 
\newcommand{\CChb}{\overline{\CCh}} 
\newcommand{\ck}{\mathscr{K}}
\newcommand{\ckb}{\overline{\ck}}
\newcommand{\ckt}{\ck^\mathrm{top}} 
\newcommand{\ctkt}{\mathbf{K}^\mathrm{top}} 
 
\newcommand{\cktb}{\ckb_\mathrm{top}}

\newcommand{\chb}{\overline{\ch}}

\newcommand{\ctop}{c_\mathrm{top}}

\newcommand{\e}{\mathbf{e}} 
\newcommand{\ec}{\check{\e}}
\newcommand{\etaa}{{\eta_{\ca}}} 
\newcommand{\etab}{\overline{\eta}}

\newcommand{\etak}{{\eta_{\ck}}}

\newcommand{\etakof}{\eta_{\okf}} 
\newcommand{\etaao}{\eta_{\oa}} 

\newcommand{\Gb}{\overline{G}}
\newcommand{\gb}{\overline{\gamma}}
\newcommand{\grq}[1]{|{#1}|_{str}} 
\newcommand{\grz}[1]{|{#1}|}

\newcommand{\II}{I\!\!I} 

\newcommand{\irightarrow}{\stackrel{\sim}{\longrightarrow}} 

\newcommand{\kt}{{K_{\mathrm{top}}}} 
\newcommand{\ktzero}{K^0_{\mathrm{top}}}
\newcommand{\ktone}{K^1_{\mathrm{top}}}
\newcommand{\ktn}{K^n_{\mathrm{top}}}
\newcommand{\ktm}{K^m_{\mathrm{top}}}
\newcommand{\M}{\overline{\mathscr{M}}}
\newcommand{\MM}{\M^G}

\newcommand{\multipl}{*} 
\newcommand{\nz}{\mathbb{Z}}

\newcommand{\pt}{\widetilde{p}}
 
\newcommand{\res}[2]{\left.{#1}\right|_{#2}}

\newcommand{\vac}{\mathbf{1}}

\newcommand{\vr}{\mathbf{vr}}

\newcommand{\taua}{\tau^{\ca}}

\newcommand{\taub}{\overline{\tau}}
\newcommand{\tauk}{\tau^{\ck}}

\newcommand{\taukof}{\tau^{\okf}} 
\newcommand{\tauao}{\tau^{\oa}} 
\newcommand{\td}{\mathbf{td}}

\newcommand{\ttau}{\boldsymbol{\tau}}

\newcommand{\Xab}{X^{\langle a, b\rangle}}

\newcommand{\stab}{\operatorname{stab}}

\newcommand{\tensor}{\otimes}
\newcommand{\ccx}{{\ci_{\cx}}} 
\newcommand{\cccx}{\ci\!\!\ci_{\cx}} 
\newcommand{\ok}{K_{\mathrm{orb}}} 
\newcommand{\okt}{\ok^{top}} 
\newcommand{\okf}{\mathsf{K}_{\mathrm{orb}}} 
\newcommand{\och}{Ch_{\mathrm{orb}}} 
\newcommand{\oa}{A^{\bullet}_{\mathrm{orb}}} 
\newcommand{\ob}{\widetilde{\crr}} 
\newcommand{\tr}{\operatorname{Tr}} 
\newcommand{\ind}{\operatorname{Ind}} 
\newcommand{\Kcond}{$(\star)$} 
\newcommand{\ee}{\iota} 
\newcommand{\eee}{\iota'} 

\newcommand{\bmp}{{\bm'}} 

\newcommand{\cep}{{\ce'}} 

\newcommand{\Deltap}{{\Delta'}}

\newcommand{\jp}{{j'}} 

\newcommand{\mmp}{{m'}} 

\newcommand{\ta}{{\tilde{a}}}

\newcommand{\tb}{{\tilde{b}}}

\newcommand{\field}{\nq}

\newarrow{Eq}=====

\newtheorem{result}{Main Result}

\newtheorem{thm}{Theorem}[section]
\newtheorem{lm}[thm]{Lemma}
\newtheorem{prop}[thm]{Proposition} 
\newtheorem{crl}[thm]{Corollary}
\newtheorem{conj}[thm]{Conjecture}

\theoremstyle{definition}
\newtheorem{rem}[thm]{Remark} 

\newtheorem{df}[thm]{Definition} 
\newtheorem{ex}[thm]{Example}
\newtheorem{df-pr}[thm]{Definition-Proposition}

\theoremstyle{remark} 
\newtheorem{nota}[thm]{Notation} \renewcommand{\thenota}{\kern-1ex}

\begin{document}
\title[Stringy K-Theory and the Chern Character] {Stringy K-theory and
  the Chern Character}


\author [T. J. Jarvis]{Tyler J. Jarvis} \address {Department of
Mathematics; Brigham Young University; Provo, UT 84602, USA}
\email{jarvis@math.byu.edu} \thanks{Research of the first author was
partially supported by NSF grant DMS-0105788}

\author [R. Kaufmann]{Ralph Kaufmann} \address {Department of
Mathematics; University of Connecticut; 196 Auditorium Road; Storrs,
CT 06269-3009, USA} \email{kaufmann@math.uconn.edu} \thanks{Research
of the second author was partially supported by NSF grant DMS-0070681}

\author [T. Kimura]{Takashi Kimura} \address {Department of
Mathematics and Statistics; 111 Cummington Street, Boston University;
Boston, MA 02215, USA } \email{kimura@math.bu.edu} \thanks{Research of
the third author was partially supported by NSF grant DMS-0204824}

\date{\today}

\begin{abstract}
We  construct two new $G$-equivariant rings: 
$\ck(X,G)$, called the \emph{stringy K-theory} of the
$G$-variety $X$,
 and $\ch(X,G)$, called the  \emph{stringy cohomology} of the
$G$-variety $X$,
for any smooth, projective variety $X$ with an action of a
finite group $G$.  For a  smooth Deligne-Mumford stack $\cx$, we also construct a new ring $\okf(\cx)$ called the \emph{full orbifold 
  K-theory} of $\cx$.  We show that for a global quotient $\cx = [X/G]$, the ring 
of $G$-invariants $\ok(\cx)$ of $\ck(X,G)$ is a subalgebra of $\okf[X/G]$ and  is linearly isomorphic to the ``orbifold
K-theory'' of Adem-Ruan \cite{AR1} (and hence Atiyah-Segal), but carries a
different ``quantum'' product which respects the natural group grading. 

We prove that there is a ring isomorphism $\CCh:\ck(X,G)\to \ch(X,G)$,
  which we call the \emph{stringy Chern character}.  We also show that there is a ring homomorphism
  $\CChf:\okf(\cx) \rTo H^\bullet_{\mathrm{orb}}(\cx)$, which we call the
  \emph{orbifold Chern character}, which induces an isomorphism $\och:\ok(\cx)\rTo H^\bullet_{\mathrm{orb}}(\cx)$ when restricted to the sub-algebra $\ok(\cx)$.  
Here $H_{\mathrm{orb}}^\bullet(\cx)$ is the Chen-Ruan
  orbifold cohomology.  We further show that $\CCh$ and $\CChf$ preserve
many properties of these algebras
and satisfy the Grothendieck-Riemann-Roch Theorem with respect to \'etale maps.  All of these results hold both in the algebro-geometric category and in the topological category for equivariant almost complex
manifolds.

We further prove that $\ch(X,G)$ is isomorphic to Fantechi and G\"ottsche's construction \cite{FG,JKK}.  Since our constructions do not use complex curves, stable maps, admissible covers, or moduli spaces, our results greatly simplify the definitions of the Fantechi-G\"ottsche ring, Chen-Ruan orbifold cohomology, and the Abramovich-Graber-Vistoli orbifold Chow ring.  

We conclude by showing that a K-theoretic version of Ruan's
Hyper-K\"ahler Resolution Conjecture holds for the symmetric product of a complex projective surface with trivial first Chern class.
\end{abstract}

\maketitle \setcounter{tocdepth}{1}

\tableofcontents

\section{Introduction}
The first main result of this paper is the construction of two new
$G$-Frobenius algebras $\ch(X,G)$ and $\ck(X,G)$, called the \emph{stringy
cohomology of $X$} and the \emph{stringy K-theory of $X$}, respectively,
where $X$ is a 
manifold with an action of a finite group
$G$. The rings of $G$-invariants of these algebras
bear some resemblance to equivariant cohomology and equivariant
K-theory, but they carry different information and generally produce a more
refined invariant than their equivariant counterparts. 

The most important part of these constructions is the multiplication, which
is defined purely in terms of the the $G$-equivariant tangent bundle $TX$
restricted to various fixed point loci of $X$.

While our stringy K-theory $\ck(X,G)$ is an entirely new construction, we
prove that our stringy cohomology $\ch(X,G)$ is equivalent to 
Fantechi and G\"ottsche's construction of stringy cohomology \cite{FG,JKK}.  Since our 
definition avoids any mention of complex curves, admissible covers, or moduli spaces, it greatly simplifies the computations of stringy cohomology and allows us to give elementary proofs of associativity and the trace axiom.  

Because our
constructions are completely functorial, an analogous construction yields the \emph{stringy Chow ring of
$X$}, which we denote by $\ca(X,G)$.  The algebra $\ca(X,G)$ is a \emph{pre-$G$-Frobenius algebra}, a generalization of a
$G$-Frobenius algebra which allows the ring to be of infinite rank and the
metric to be degenerate.

The second main result of this paper is the introduction of a new  \emph{stringy Chern
character} $\CCh:\ck(X,G)\to \ch(X,G)$.   We prove that $\CCh$ is a
\emph{ring} isomorphism  which preserves all of the properties of a
pre-$G$-Frobenius algebra except those involving the metric.

The third main result of this paper is the introduction of two new \emph{orbifold
K-theories}.  The first we call \emph{full orbifold K-theory} and is defined for a general almost complex orbifold (or a smooth Deligne-Mumford stack).  We denote it by  $\okf(\cx)$.  The second algebra is defined when $\cx = [X/G]$ is  a global quotient by a finite group as the algebra of invariants $\ck(X,G)^G$ of the stringy K-theory of $X$.  We  denote this algebra by $\ok(\cx)$ and call it the \emph{small orbifold K-theory of $\cx$}.  It is  
linearly isomorphic to the construction of Adem and Ruan \cite{AR1}, but our
construction possesses a different, ``quantum,'' product.  We show there is a natural homomorphism of algebras $\okf(\cx) \rTo^{\pi^*} \ok(\cx)$, and an \emph{orbifold Chern character} $\CChf:\okf(\cx) \to H^\bullet_{\mathrm{orb}}(\cx)$ which, like the stringy Chern character, is 
a ring homomorphism which preserves all of the properties of a Frobenius
algebra that do not involve the metric.  In the special case that the
orbifold is a global quotient $\cx=[X/G]$, the orbifold Chern character induces an \emph{isomorphism} $\och:\ok(\cx) \rTo H^\bullet_{\mathrm{orb}}(\cx)$ which agrees with that induced by the stringy Chern character on the rings of
$G$-invariants.

Our results are initially formulated and proved in
the algebro-geometric category, with Chow rings and algebraic
K-theory, but they also hold in the topological category, with cohomology
and topological K-theory (see Section~\ref{sec:top}) for almost
complex manifolds with a $G$-equivariant almost complex structure. In fact, these algebraic structures depend only upon the homotopy
class of the $G$-equivariant almost complex structure.  Our results
can also be generalized to equivariant stable complex manifolds (see
Remark \ref{rem:StableACM}).

\subsection{Notation and conventions}

Unless otherwise specified, we assume throughout the paper that all
cohomology rings have coefficients in the rational numbers $\nq$.
Also, unless otherwise specified, all groups are finite and all group
actions are left actions.

The stack (or orbifold) quotient of a variety (or manifold) $X$ by $G$ will be denoted $\xg $ and the coarse moduli space (i.e., underlying space) of this quotient will be denoted $X/G$.

The conjugacy class of any element $g$ in a group $G$ will be denoted
$[\![g]\!]$ and the commutator $aba^{-1}b^{-1}$ of two elements $a,b\in
G$ is denoted $[a,b]$. 

\subsection{Background and motivation}

We now describe part of our motivation for studying stringy K-theory.
For convenience, we assume throughout this subsection that the
coefficient ring is $\nc$ rather than $\nq$.

Let $Y$ be a projective, complex surface such that $c_1(Y)=0$.  For
all $n$, consider the product $Y^n$ with the symmetric group $S_n$ acting by permuting its factors. The quotient orbifold $[Y^n/S_n]$ is called the
symmetric product of $Y$. Let $Y^{[n]}$ denote the Hilbert scheme of
$n$ points on $Y$. The morphism $Y^{[n]}\to Y^n/S_n$ is a crepant
resolution of singularities and is, furthermore, a hyper-K\"ahler
resolution \cite{Ru}. Fantechi and G\"ottsche \cite{FG} proved that
there is a ring isomorphism $\psi':H^\bullet_\mathrm{orb}([Y^n/S_n])\to
H^\bullet(Y^{[n]})$, where $H^\bullet(Y^{[n]})$ is the ordinary
cohomology ring (see also \cite{Ka4,Ur}).\footnote{In fact, they proved
  that the isomorphism holds over $\nq$, provided that the
  multiplication on $H^\bullet_\mathrm{orb}([Y^n/S_n])$ is twisted by signs.
  This sign change can be regarded as a kind of discrete torsion (see
  Section~\ref{eq:SubsectionSP} for more details).}

The previous example is a verification, in a special case, of the
following conjecture of Ruan \cite{Ru}, which was inspired by the work
of string theorists studying topological string theory on orbifolds.
\begin{conj}[Cohomological Hyper-K\"ahler Resolution
  Conjecture]\label{conj:Ruan} Suppose that $\tilde{V}\to V$ is a
  hyper-K\"ahler resolution of the coarse moduli space $V$ of an orbifold $\cv$.
 The ordinary cohomology ring
  $H^\bullet(\tilde{V})$ of $\tilde{V}$ is isomorphic (up to discrete torsion) to the Chen-Ruan
  orbifold cohomology ring $H^\bullet_\mathrm{orb}(\cv)$ of $\cv$.
\end{conj}
Let us return again to the example of the symmetric product. The
algebra isomorphism $\psi':H^\bullet_\mathrm{orb}([Y^n/S_n])\to
H^\bullet(Y^{[n]})$ suggests that there should exist a K-theoretic
analogue $\ok([Y^n/S_n])$ of $H^\bullet_\mathrm{orb}([Y^n/S_n])$, a stringy
Chern character isomorphism $\och:\ok([Y^n/S_n])\to
H^\bullet_\mathrm{orb}([Y^n/S_n])$, and an algebra isomorphism
$\psi:\ok([Y^n/S_n]) \to K(Y^{[n]})$, such that the following diagram
commutes:
\begin{equation}
\begin{diagram}
  \ok([Y^n/S_n]) & \rTo^{\och} & H^\bullet_\mathrm{orb}([Y^n/S_n]) \\
  \dTo^{\psi} & & \dTo^{\psi'} \\ K(Y^{[n]}) & \rTo^{\cch} &
  H^\bullet(Y^{[n]}).\\
\end{diagram}
\end{equation}

In this paper we construct an \emph{orbifold K-theory} analogous to the Chen-Ruan orbifold cohomology, and we construct an orbifold Chern character which is a ring isomorphism (see Theorem~\ref{thm:orbFACh}).  This leads us to pose the following K-theoretic analogue of the Ruan conjecture.
\begin{conj}[K-theoretic Hyper-K\"ahler Resolution
  Conjecture]\label{conj:KHRC} Suppose that $\tilde{V}\to V$ is a
  hyper-K\"ahler resolution of the coarse moduli space $V$ of an orbifold $\cv$.
 The ordinary K-theory 
  $K(\tilde{V})$ of the resolution $\tilde{V}$ is isomorphic (up to discrete torsion) to the (small) orbifold K-theory $\ok(\cv)$ of $\cv$.
\end{conj}

The method that Fantechi and G\"ottsche use to prove their result involves
the construction of a new ring $\ch(X,G)$, which we call \emph{stringy
cohomology}, associated to any smooth, projective manifold $X$ with an action
by a finite group $G$.  They show that for a global quotient orbifold
$\cx:=[X/G]$, the Chen-Ruan orbifold cohomology $H^\bullet_\mathrm{orb}(\cx)$ is 
isomorphic to the ring of 
invariants $\ch(X,G)^G$ of the stringy cohomology. 

Their construction suggests that a similar construction in K-theory should be
possible and that the two constructions might be related by a stringy Chern 
character. 

\subsection{Summary and discussion of main results}\label{sec:summary}

We will now briefly describe the main results and constructions of the paper.

Let $X$ be a smooth, projective variety with an action of a finite group
$G$. For each $m\in G$ we denote the fixed locus of $m$ in $X$ by $X^m$, and we let $$I_G(X):=\coprod_{m\in G} X^m\subset X\times G$$ denote the
\emph{inertia variety} of $X$.
The inertia variety $I_G(X)$ should not be confused with the inertia
\emph{orbifold}, 
or inertia \emph{stack}, $\coprod_{[\![ g ]\!]}[ X^{g}/Z_G(g)]$, where the sum runs over
conjugacy classes $[\![ g ]\!]$ in $G$.
Note that the $G$-variety $I_G(X)$ contains $X = X^1$ as a connected component. 

As a $G$-graded $G$-module, the
\emph{stringy Chow ring $\ca(X,G)$ of $X$} is the Chow ring of $I_G(X)$,
i.e., $$\ca(X,G) = \bigoplus_{g\in G} \ca_g(X) = \bigoplus_{g\in G} A^\bullet(X^g).$$ 
The inertia variety has a canonical 
$G$-equivariant involution $\sigma:I_G(X)\to I_G(X)$ 
which maps $X^m$ to $X^{m^{-1}}$ via 
\begin{equation}\label{eq:invol}
\sigma: (x,m)\mapsto (x,m^{-1})
\end{equation}
for all $m$ in $G$. We define a
pairing $\etaa$ on $\ca(X,G)$ by
\[
\etaa(v_1,v_2) := \int_{[I_G(X)]} v_1\cup\sigma^*v_2
\]
for all $v_1,v_2$ in $\ca(X,G)$.

In a similar fashion, we define the \emph{stringy K-theory $\ck(X,G)$ of $X$},
as a $G$-graded $G$-module, to be the K-theory of the inertia variety, i.e.,  
$$\ck(X,G) = \bigoplus_{g\in G} \ck_g(X) = \bigoplus_{g\in G} K(X^g).$$ 
We define a 
pairing $\etak$ on $\ck(X,G)$
by
\[
\etak(\cf_1,\cf_2) := \chi (I_G(X), \cf_1\otimes\sigma^*\cf_2)
\]
for all $\cf_1,\cf_2$ in $\ck(X,G)$,
where $\chi(I_G(X), \cf) $ denotes the Euler characteristic of $\cf \in K(I_G(X))$.

The definition of the multiplicative structure on $\ca(X,G)$ and $\ck(X,G)$ requires the following new constructions.
\begin{df}\label{df:kage}
Define $\cs$ in $K(I_G(X))$ (the rational K-theory) to be such that for any
$m$ in $G$, its restriction $\cs_m$ in $K(X^m)$ is given by
\begin{equation}\label{eq:S}
\cs_m:= \res{\cs}{X^m}:= \bigoplus_{k=0}^{r-1}\frac{k}{r}W_{m,k},
\end{equation} 
where $r$ is the order of $m$, and $W_{m,k}$ is the eigenbundle of
$W_m := \res{TX}{X^{m}}$ such that $m$ acts with eigenvalue 
$\exp(2\pi k i/r)$.  
\end{df}
The virtual rank $a(m)$ of $\cs_m$ is called the \emph{age of $m$} and is a
locally constant $\nq$-valued function on $X^m$.
\begin{rem}\label{rem:TangentNormal}
It is worth pointing out that $\cs$
would remain the the same if, in the definition of $\cs_m$, the restriction
of $TX$ to $X^m$ were replaced by the normal bundle of $X^m$ in $X$.  For
this reason, the construction of $\cs$ and the construction of stringy
cohomology and stringy K-theory still works over stable complex manifolds.
\end{rem}

For any triple $\bm := (m_1,m_2,m_3)$ in $G^3$ such that $m_1 m_2 m_3 = 1$, we let $X^\bm := X^{m_1}\cap X^{m_2} \cap X^{m_3}$,  where $X^{m_i}$ is
regarded as a subvariety of $X$.  
\begin{df}\label{df:R}
Define the element $\crr(\bm)$ in $K(X^\bm)$
by
\begin{equation} \label{eq:R}
\crr (\bm) :=
  TX^\bm\ominus \res{TX}{X^\bm}\oplus\bigoplus_{i=1}^3
  \res{\cs_{m_i}}{X^{\bm}}.
\end{equation}
\end{df}

It is central to our theory, but not at all obvious, that $\crr(\bm)$ is actually 
represented by  a vector bundle on $X^\bm$.  
In general, the only way we know how to establish this 
key fact is through our proof in Section~\ref{sec:obs}, which uses the
Eichler trace formula (a special case of the holomorphic Lefschetz Theorem) 
to show that $\crr(\bm)$ is equal to the
obstruction bundle $R^1\pi^G_*f^*(TX)$ arising in the Fantechi-G\"ottsche
construction of stringy cohomology.  However, once one knows that $\crr(\bm)$
is always represented by a vector bundle, 
all of the properties of a pre-G-Frobenius algebra 
can be established (see Definition~\ref{df:preGFA} for details).
We first use $\crr(\bm)$ to define 
the multiplication in $\ca(X,G)$ and
$\ck(X,G)$ as follows. 
For all $i=1,2,3$, let
$$\e_{m_i} : X^\bm\to X^{m_i}$$ be the canonical inclusion
morphisms, and define
$$\ec_{m_i}:=\sigma\circ\e_{m_i}:X^\bm\to X^{m_i^{-1}},$$ 
where $\sigma$ is the canonical involution (see Equation~(\ref{eq:invol})).
\begin{df}\label{df:chow-mult}
 Given $m_1, m_2 \in G$, let $m_3:=(m_1m_2)^{-1}$.  For any
$v_{m_1}\in \ca_{m_1}(X)$ and $v_{m_2}\in \ca_{m_2}(X)$, we define the
\emph{stringy product 
(or multiplication)} of $v_{m_1}$ and $v_{m_2}$  in $\ca(X,G)$ to be
\begin{equation}\label{eq:DefChowMult}
v_{m_1} \multipl v_{m_2} := \ec_{m_3 *}\left(\e_{m_1}^*v_{m_1} \cup \e_{m_2}^*
v_{m_2} \cup \ctop\left(\crr(\bm)\right)\right).
\end{equation}
The product is then extended linearly to all of $\ca(X,G)$.  
\end{df}

We define the \emph{stringy product 
on $\ck(X,G)$} analogously.   
\begin{df} \label{df:K-mult}
Given $m_1, m_2 \in G$, let $m_3:=(m_1m_2)^{-1}$.  For any $\cf_{m_1}\in \ck_{m_1}(X)$ and $\cf_{m_2}\in \ck_{m_2}(X)$, we define the \emph{stringy product} of $\cf_{m_1}$ and $\cf_{m_2}$  in $\ck(X,G)$ to be
\begin{equation}\label{eq:DefKMult}
\cf_{m_1} \multipl \cf_{m_2} := \ec_{m_3 *}\left(\e_{m_1}^*\cf_{m_1} \otimes \e_{m_2}^*
\cf_{m_2} \otimes \lambda_{-1}\left(\crr(\bm)^*\right)\right),
\end{equation}
and again the product is extended linearly to all of $\ck(X,G)$.  
\end{df}

The stringy Chow ring and stringy K-theory are almost $G$-Frobenius algebras,
but they are generally infinite dimensional and have degenerate pairings.  An
algebra which satisfies essentially all of the axioms of a $G$-Frobenius
algebra except those involving finite dimensionality and a non-degenerate
pairing is called a \emph{pre-$G$-Frobenius algebra} (see
Definition~\ref{df:preGFA} for details).   

Our first main result is the following. 
\begin{result} (See Theorems~\ref{thm:ChowGFA} and \ref{thm:KGFA} for complete details.)
For any smooth, projective variety $X$ with an action of a finite
group $G$, the 
ring $\ca(X,G)$ is a $\nq$-graded, pre-$G$-Frobenius algebra
which contains the ordinary Chow ring $A^\bullet(X) = \ca_1(X)$ of $X$ as a sub-algebra.
In particular, $\ca(X,G)$ is a $G$-equivariant, associative ring (generally non-commutative) with a $\nq$-grading that respects the multiplication and the metric. 

Similarly, the 
ring $\ck(X,G)$
       is a pre-$G$-Frobenius algebra
which contains the ordinary K-theory $K(X) = \ck_1(X)$ of $X$ as a sub-algebra.  
\end{result}

Unfortunately, the 
ordinary Chern character $\Ch:\ck(X,G)\to \ca(X,G)$ does not respect the
stringy multiplications.
We repair this problem by defining
the \emph{stringy Chern character} $\CCh: \ck(X,G) \to \ca(X,G)$ to be
a deformation of the ordinary Chern character. 
That is, for every element $m\in G$ and every $\cf_m \in \ck_m(X)$ we define
\begin{equation}\label{eq:DefCCh}
\CCh(\cf_m) := \Ch(\cf_m)\cup \td^{-1}(\cs_m)
= \Ch(\cf_m)\cup (\vac - c_1(\cs_m)/2 +  \cdots),
\end{equation}
where $\cs_m$ 
is defined in Equation (\ref{eq:S}).
This yields our second main result.
\begin{result} (See Theorem~\ref{thm:ChernHomo} and Theorem~\ref{thm:RRfunct}
for complete details.) 
  The stringy Chern character $\CCh:\ck(X,G)\to \ca(X,G)$ 
is a $G$-equivariant algebra isomorphism.
Moreover, $\CCh$ is natural and satisfies a form of the
Grothendieck-Riemann-Roch Theorem
with respect to $G$-equivariant \'etale maps.
\end{result}

It is natural to ask whether the rings of $G$-invariants of the stringy Chow ring and stringy K-theory are presentation independent, and if so, whether these rings can be constructed for orbifolds which are not global quotients of a variety by a finite group. The answer is yes in both cases.  
It is already known that 
the ring of $G$-invariants of the stringy Chow ring $\ca(X,G)^G$
is isomorphic to the Abramovich-Graber-Vistoli orbifold Chow
ring $\oa([X/G])$ of the quotient orbifold $[X/G]$.

The third main result of this paper 
has three parts:
first, the construction of a new
\emph{full orbifold K-theory} $\okf(\cx)$ for smooth Deligne-Mumford stacks, and a second \emph{small orbifold K-theory} $\ok([X/G])$ for global quotients by finite groups; second, the construction of an orbifold Chern character $\och:\ok(\cx) \rTo \oa(\cx)$ which is a
ring homomorphism; and third, a demonstration of the relations between the two theories.

\begin{result}(See Theorems~\ref{thm:ccx-preFA} and \ref{thm:orbFACh} for complete details.)
For a smooth Deligne-Mumford stack $\cx$ satisfying the resolution property,
the full orbifold K-theory $\okf(\cx)$ is a pre-Frobenius algebra.  Moreover, there is a \emph{full orbifold Chern character}
$\CChf:\okf(\cx) \to \oa(\cx)$ which, like the stringy Chern character, is 
a ring homomorphism which preserves all 
of the properties of a pre-Frobenius algebra that do not involve the metric.  

For a global quotient $\cx = [X/G]$ by a finite group, $G$ the small orbifold K-theory $\ok(\cx)$ is also a pre-Frobenius algebra, independent of the choice of resolution. 
There is an \emph{orbifold Chern character} $\och:\ok(\cx) \rTo \oa(\cx)$ which is an algebra \emph{isomorphism}, and 
there is a natural algebra homomorphism $\pi^*:\okf(\cx) \rTo \ok(\cx)$
making the following diagram commute:
$$\begin{diagram}
\okf(\cx) & \rTo^{\pi^*}      & \ck(X,G)^G & = & \ok(\cx)\\
\dTo^{\CChf} &		     & & &\dTo_{\och}^{\cong} &\\
 H_\mathrm{orb}^{\bullet}(\cx) & \rEq & &&\ch^\bullet(X,G)^G.
\end{diagram}
$$
\end{result}

As we mentioned above, all these results are proved initially in the
algebro-geometric category, but we prove in Section~\ref{sec:top} that
their analogues in the topological category also hold. 
That is, we define stringy cohomology, stringy topological K-theory, 
orbifold cohomology, orbifold topological K-theory, and their corresponding
Chern characters. We prove theorems analogous to the above for these 
topological constructions. Furthermore, we prove that stringy cohomology of a complex (or almost complex) orbifold $\cx$ 
is equal to the $G$-Frobenius algebra $\ch(X,G)$ described in \cite{FG,JKK}. Similarly, we prove that
the orbifold cohomology $H^\bullet_{\mathrm{orb}}(\cx)$ of a complex (or almost complex) orbifold
$\cx$ is equal to its Chen-Ruan orbifold cohomology \cite{CR1,AGV}.

We conclude with an application of 
these results to the case of the symmetric product
of a smooth, projective surface $Y$ with trivial canonical bundle and
verify that our K-theoretic Hyper-K\"ahler Resolution Conjecture
(Conjecture~\ref{conj:KHRC}) holds in this case; that is,
$\ok([Y^n/S_n])$ is isomorphic to $K(Y^{[n]})$.

\subsection{Directions for further research}

These results suggest many different directions for further research.
The first is to generalize to the case where $G$ is a Lie group and to
higher-degree Gromov-Witten invariants. This will be explored
elsewhere.  It would also be interesting to study stringy
generalizations of the usual algebraic structures of K-theory, e.g.,
the Adam's operations and $\lambda$-rings. Another interesting
direction would be to study stringy generalizations of other
K-theories, including algebraic K-theory and higher K-theory.  It
would also be very interesting to find an analogous construction in
orbifold conformal field theory, e.g., twisted vertex algebras and the
chiral de~Rham complex \cite{FrSz}. Finally, it would be interesting
to see if our results can shed light upon the relationship between
Hochschild cohomology and orbifold cohomology \cite{DoEt} in the
context of deformation quantization.

\subsection{Acknowledgments}

We would like to thank D.~Fried, J.~Morava, S.~Rosenberg, and Y.~Ruan
for helpful discussions.  We would also like to thank J.~Stasheff for
his useful remarks about the exposition. The second and third author
would like to thank the Institut des Hautes \'Etudes Scientifiques,
where much of the work was done, for its financial support and
hospitality, and the second author would also like to thank the
Max-Planck Institut f\"ur Mathematik in Bonn for its financial support
and hospitality.  Finally, we thank the referees for their helpful suggestions.

\section{The ordinary Chow ring and K-theory of a variety}

In this section, we briefly review some basic facts about the classical Chow
ring, K-theory, and certain characteristic classes that we will need.
Throughout this section, all varieties we consider will be smooth,
projective varieties over $\nc$.

Recall that a \emph{Frobenius algebra} is a finite dimensional, unital, commutative,
associative algebra with an invariant (non-degenerate) metric. To each smooth,
projective variety $X$, one can associate two algebras, which are almost
Frobenius algebras,  namely, the \emph{Chow ring $A^\bullet(X)$ of $X$}, and
\emph{the K-theory $K(X)$ of $X$}. These fail to be Frobenius algebras in
that they may be infinite dimensional and their symmetric
pairing may be degenerate.
Both $A^\bullet(X)$ and $K(X)$ also possess an additional structure, which we call a \emph{trace element}, which is closely related to the Euler characteristic.
We call such algebras (with a trace element) \emph{pre-Frobenius algebras}
(see Definition~\ref{df:preFA}).

Furthermore, there is an isomorphism of unital, commutative, associative
algebras $\Ch:K(X)\to A^\bullet(X)$ called the \emph{Chern character}.  The
Chern character does not  preserve the metric, but it does preserve the 
trace elements.  We call such an isomorphism \emph{allometric}.  We will
now briefly review these constructions in order to fix notation and
conventions, referring the interested reader to \cite{Fu,FuLa} for more details.

\subsection{The Chow ring}

The Chow ring of a smooth, projective variety $X$ is additively a
$\nz$-graded Abelian group $A^\bullet(X,\nz) = \bigoplus_{p=0}^D
A^p(X,\nz)$, where $D$ is the dimension of $X$, and $A^p(X)$ is the
group of finite formal sums of $(D-p)$-dimensional subvarieties of
$X$, modulo rational equivalence.

In this paper we will always work with rational coefficients, and we
write
$$A^\bullet(X):=A^\bullet(X,\nz)\tensor_\nz \nq.$$
The vector space
$A^\bullet(X)$ is endowed with a commutative, associative
multiplication which preserves the $\nz$-grading, arising from the
intersection product, and possesses an identity element $\vac := [X]$
in $A^0(X)$. The intersection product $A^p(X)\otimes A^q(X)\to
A^{p+q}(X)$ is denoted by $v\otimes w\mapsto v\cup w$ for all $p,q$.

Given a proper morphism $f:X\to Y$ between two varieties, there is an
induced pushforward morphism $f_*:A^\bullet(X)\to A^\bullet(Y)$.  In
particular, if $Y$ is a point and $f:X\to Y$ is the obvious map, then
one can define integration via the formula
\[
\int_{[X]} v := f_*(v)
\]
for all $v$ in $A^\bullet(X)$. 
Whenever $X$ is equidimensional of dimension $D$, 
the integral vanishes unless $v$
belongs to $A^D(X)$. Define a symmetric, bilinear form
$\eta_A:A^\bullet(X)\otimes A^\bullet(X)\to \nq$ via $\eta_A(v,w) :=
\int_{[X]} v\cup w$. 
Finally, we define 
a special element $\tau^A  \in (A^{\bullet}(X))^*$ in the dual,
the \emph{trace element} 
$$\tau^A(v) := \int_X (v\cup \ctop(TX)),$$ 
where $\ctop(TX)$ is the top Chern class of the tangent bundle
$TX$. The integer $\tau^A(\vac)$ is the usual Euler characteristic of $X$.

Although the Chow ring resembles the cohomology ring in many ways, it is
important to note that $A^\bullet(X)$ 
is generally infinite dimensional, 
and that the pairing $\eta_A$ is often degenerate.
This motivates the following definition.

\begin{df}\label{df:preFA}\footnote{We thank the referee 
for help clarifying these details.}
Consider a  tuple $(R,\multipl,\eta,\vac,\tau)$ consisting of 
a commutative, associative
algebra
$(R,\multipl)$ (possibly infinite dimensional) 
with unity $\vac \in R$, 
a symmetric bilinear pairing $\eta$ (possibly degenerate),
and $\tau$ in $R^*$, called the \emph{trace element}.
We say that $(R,\multipl,\eta,\vac,\tau)$ is a \emph{pre-Frobenius algebra} if
the pairing
$\eta$ is multiplicatively invariant: 
$$ \eta(r\multipl s,t) = \eta(r,s\multipl t)$$ for all $r,s,t \in R$.  

Every Frobenius algebra $(R,\multipl,\eta,\vac)$ has a \emph{canonical trace
element} $\tau(v) := \tr_R(L_v)$ where $L_v$ is left multiplication by $v$ in
$R$. We call $(R,\multipl,\eta,\vac,\tau)$ the \emph{canonical pre-Frobenius
algebra structure} associated to the Frobenius algebra $(R,\multipl,\eta,\vac).$ 
\end{df}

\begin{prop}(see \cite[\S1]{Kl} and \cite[\S19.1]{Fu})
Let $A^\bullet(X)$ be the Chow ring of an irreducible, smooth, projective
variety $X$. 
\begin{enumerate}
\item The triple $(A^\bullet(X),\cup,\vac,\eta_A,\tau^A)$ is a
      pre-Frobenius algebra
      graded by $\nz$.  
\item If $f:X\to Y$ is any morphism, then the associated pullback
  morphism $f^*:A^\bullet(Y)\to A^\bullet(X)$ is a morphism of
commutative, associative algebras graded by $\nz$. 
\item (Projection formula) For any proper morphism $f:X\to Y$, if
  $\alpha \in A^\bullet(X)$ and $\beta \in A^\bullet(Y),$ we have
$$ f_*(\alpha \cup f^*(\beta)) = f_*(\alpha) \cup \beta.$$
\end{enumerate}
\end{prop}

\subsection{K-theory}

$K(X;\nz)$ is additively equal to the free Abelian group generated by
isomorphism classes of (complex algebraic) vector bundles on $X$,
modulo the subgroup generated by
\begin{equation}\label{eq:KRelation}
[E]\ominus [E'] \ominus [E'']
\end{equation}
for each exact sequence of vector bundles
\begin{equation}\label{eq:ModExact}
0\to E'\to E\to E''\to 0.
\end{equation}
Here $\ominus$ denotes subtraction and $\oplus$ denotes addition in
the free Abelian group. We define $$K(X) :=
K(X;\nz)\otimes_{\nz}\nq.$$
The multiplication operation, also denoted
by $\otimes$, taking $K(X)\otimes K(X)\to K(X)$ is the usual tensor
product $[E]\otimes [E']\mapsto [E\otimes E']$ for all vector bundles
$E$ and $E'$.  We denote the multiplicative identity by $\vac
:=[\co_X]$.

Given a proper morphism $f:X\to Y$ between two smooth varieties, there
is an induced pushforward morphism $f_*:K(X)\to K(Y)$ given by
$f_*([E])=\sum^D_{i=0}(-1)^i R^if_*E,$ 
where $D$ is the relative dimension of $f$. 
In particular, if $Y$ is a
point and $f:X\to Y$ is the obvious map, then the Euler characteristic
of $v \in K(X)$ is the pushforward
\[
\chi(X,v)=f_*(v). 
\]
If $X$ is irreducible, we define a symmetric 
bilinear form
$\eta_K:K(X)\otimes K(X)\to \nq$ via $\eta_K(v,w) := \chi(X,v \otimes w).$ 

While $K(X)$ does not have a $\nz$-grading like $A^\bullet(X)$, it
does have a \emph{virtual rank (or augmentation)}. That is, for each
connected component $U$ of $X$, there is a surjective ring
homomorphism $\vr:K(U)\to \nq$ which assigns to each vector bundle $E$
on $U$ its rank.  In addition, $K(X)$ has an 
involution which takes a vector bundle $[E]$ to its dual $[E^*]$.

Another important property of $K$-theory is that it is a so-called
\emph{$\lambda$-ring}. That is, for every non-negative integer $i$,
there is a map $\lambda^i:K(Y)\to K(Y)$ defined by $\lambda^i([E]) :=
[\bigwedge^i E]$, where $\bigwedge^i E$ is the $i$-th exterior power of
the vector bundle $E$. In particular, $\lambda^0([E]) = \vac$, and
$\lambda^i([E]) = 0$ if $i$ is greater than the rank of $E$.

These maps satisfy the relations
\[
\lambda^k(\cf\oplus\cf') = \bigoplus_{i=0}^k
\lambda^i(\cf)\lambda^{k-i}(\cf')
\]
for all $k=0,1,2,\dots$ and all $\cf$, and $\cf'$ in $K(Y)$.  These
relations can be neatly stated in terms of the universal formal power
series in $t$
\begin{equation}
\lambda_t(\cf) := \bigoplus_{i=0}^\infty \lambda^i(\cf) t^i
\end{equation}
by demanding that $\lambda_t$ satisfy the multiplicativity relation
\begin{equation}\label{eq:LambdaMult}
\lambda_t(\cf\oplus\cf') = \lambda_t(\cf)\lambda_t(\cf').
\end{equation}

If $E$ is a rank-$r$ vector bundle over $X$, then one can define
\[
\lambda_{-1}([E]) := \bigoplus_{i=0}^r (-1)^i \lambda^i([E])
\]
in $K(X)$, which will play an important role in this paper.
In particular, $\lambda_{-1}([E^*])$ is the \emph{K-theoretic Euler class of $E$}.

Like the Chow ring, the ring $K(X)$ is not quite a Frobenius algebra, because it is generally infinite dimensional, and its pairing may be degenerate; however, if we define the the \emph{trace element} as  
\begin{equation}
\tau^K(v) := \chi(X,\lambda_{-1}(T^*X)\otimes v)
\end{equation} 
for all $v$ in $K(X)$, then we have the following proposition.
\begin{prop}\label{prop:BasicK}
Let $K(X)$ be the K-theory of an irreducible, smooth, projective variety $X$.
\begin{enumerate}
\item The tuple $(K(X),\otimes,\vac,\eta_K,\tau^K)$ is a pre-Frobenius
  algebra.
\item If $f:X\to Y$ is any morphism, then the associated pullback morphism
      $f^*:K(Y)\to K(X)$ is a 
morphism of commutative, associative algebras. 
\item (Projection formula) For any proper morphism $f:X\to Y$, if
  $\alpha \in K(X)$ and $\beta \in K(Y)$ we have $$f_*(\alpha \cup
  f^*(\beta)) = f_*(\alpha) \cup \beta.$$
\end{enumerate}
\end{prop}

\subsection{Chern classes, Todd classes, and the Chern character}

The \emph{Chern polynomial of $\cf$} in $K(X)$ is defined to be the
universal formal power series in $t$
\[
c_t(\cf) := \sum_{i=0}^\infty c_i(\cf) t^i,
\]
where $c_i(\cf)$, the \emph{$i$-th Chern class of $\cf$}, belongs to
$A^i(X)$  for all $i$, 
and $c_t$ and the $c_i$ satisfy the following axioms:
\begin{enumerate}
\item If $\cf = [\co(D)]$ is a line bundle defined by a divisor $D$,
  then $$c_t(\cf) = \vac + Dt.$$
\item The Chern classes commute with pullback, i.e., if $f:X\to Y$ is
  any morphism, then $c_i(f^*\cf) = f^* c_i(\cf)$ for all $\cf$ in
  $K(X)$ and all $i$.
\item If $$0 \to \cf' \to \cf \to \cf'' \to 0$$
  is an exact sequence, then $$ c_t(\cf ) =
  c_t(\cf')c_t(\cf'').\label{eq.chern-mult}
$$
\end{enumerate}
In particular, $c_0(\cf) = \vac$ for all $\cf$.

A fundamental tool is the \emph{splitting principle}, which says that
for any vector bundle $E$ on $X$ of rank $r$, 
there is a morphism $f: Y \to X$,
such that $f^*: A^\bullet(X) \to A^\bullet(Y)$ is injective, and
$f^*([E])$ splits (in K-theory) as a sum of line bundles:
\begin{equation}
f^*([E]) = [\cl_1] \oplus \cdots \oplus [\cl_r].
\end{equation}
We define the \emph{Chern roots} of $[E]$ to be
$a_i:= c_1(\cl_i)$, and thus by Property~(\ref{eq.chern-mult}) of the
Chern polynomial, we have
\begin{equation}
c_t([E])  =  \prod_{i=1}^r (\vac + a_i t).
\end{equation}
Of course, the Chern roots depend on the choice of $f$, but any
relations derived in this way among the Chern classes of $[E]$ will
hold in $A^{\bullet}(X)$ regardless of the choice of $f$.

From the Chern classes, one can construct the \emph{Chern character}
$$\Ch:K(X)\to A^\bullet(X)$$
by associating to a rank-$r$ vector
bundle $E$ over $X$ the element
\begin{equation}\label{eq:ChernCharDef}
\Ch([E]) := \sum_{i=1}^r \exp(a_i) = 
r + c_1([E]) + \frac12(c_1^2([E]) -2c_2([E])) + \cdots,
\end{equation}
where $a_1,\ldots,a_r$ are the Chern roots of $[E]$.

For each connected component $U$ of $X$, the \emph{virtual rank} is the
algebra homomorphism $\vr:K(U)\to\nq$, which is the composition of
$\Ch:K(U)\to A^\bullet(U)$ with the canonical projection
$A^\bullet(U)\to A^0(U)\cong\nq$.

\begin{rem}\label{rem:ChNotMetricPres}
  In general, the Chern character does not commute with pushforward.
  That is the content of the Grothendieck-Riemann-Roch theorem, which
  we will review shortly. Since the pairings
of both $K(X)$ and
  $A^{\bullet}(X)$ are defined by pushforward, this means the Chern
  character does not respect the 
pairings.
\end{rem}

To state the Grothendieck-Riemann-Roch theorem, we need the
\emph{Todd class} $\td:K(X)\to A^\bullet(X)$, which is defined by
imposing the multiplicativity condition
\[
\td(\cf\oplus\cf') = \td(\cf)\td(\cf')
\]
for all $\cf$, $\cf'$ in $K(X)$, and by also demanding that if $E$ is
a rank $r$ vector bundle on $X$, then
\[
\td([E]) := \prod_{i=1}^r \phi(a_i),
\]
where $a_i=1,\ldots,r$ are the Chern roots of $[E]$ and
\[
\phi(t) := \frac{t e^{t}}{e^{t} -1}
\]
is regarded as a element in $\nq[\![t]\!]$. Therefore, $\td(\cf) =
\vac+x$, where $x 
= c_1(\cf) + (c_1^2(\cf) + c_2(\cf))/12 + \cdots$ 
belongs to $\bigoplus_{i=1}^D A^i(X)$. 

\begin{thm}[Grothendieck-Riemann-Roch]\label{thm:GRR}
  For any proper morphism $f:X\to Y$ of non-singular varieties and any
  $\cf \in K(X)$, we have
\begin{equation}
\Ch(f_*(\cf)) \cup \td (TY) = f_*(\Ch(\cf) \cup \td(TX)),
\end{equation}
where $TX$ and $TY$ are the tangent bundles of $X$ and $Y$,
respectively.
\end{thm}

The following useful proposition intertwines many of the structures discussed in this section.
\begin{prop}\label{prop:ToddCh} \cite[Prop I.5.3]{FuLa}
If $E$ is a vector bundle of rank $r$ over $X$, then
the following identity holds in $A^\bullet(X)$:
\begin{equation}\label{eq:UsefulMix}
\td([E]) \Ch(\lambda_{-1}([E^*])) = \ctop([E]),
\end{equation}
where $\ctop([E])$ is the top Chern class $c_r([E])$.
\end{prop}

\begin{nota}
When $E$ is a vector bundle over $X$, we will often write $c_t(E)$
instead of $c_t([E])$, and similarly for $\lambda_t$, $\td$ and $\Ch$.
\end{nota}

We are now ready to state the key property of the Chern
character. 
As mentioned in 
Remark~\ref{rem:ChNotMetricPres}, 
the Grothendieck-Riemann-Roch Theorem 
implies that the Chern character cannot preserve the pairings, but it does preserve all of the other structures.

\begin{df}\label{df:allometric-FA}
An \emph{allometric isomorphism} $\phi:(R,\multipl,\eta,\vac,\tau) \rTo
(R',\star,\eta',\vac',\tau')$ of pre-Frobenius algebras is an isomorphism of
unital, associative algebras that does not necessarily preserve the metric
but does preserve the trace elements:
$$\phi^*\tau' = \tau.$$
\end{df}

We have the following theorem.
\begin{thm}\label{thm.ch-isom}
  The Chern character $\Ch:K(X)\to A^\bullet(X)$ is an
allometric isomorphism of  pre-Frobenius algebras. Furthermore, if $f:X\to Y$
is any morphism, then the following diagram commutes: 
\begin{equation}
\begin{diagram}
K(Y) & \rTo^{f^*} & K(X) \\ \dTo^{\Ch} & & \dTo^{\Ch} \\
A^{\bullet}(Y) & \rTo^{f^*} & A^{\bullet}(X).\\
\end{diagram}
\end{equation}
\end{thm}
\begin{proof}
The only nonstandard part of this statement is that $\Ch$ preserves the trace
elements.  This can be seen as follows. For all $\cf$ in $K(X)$, 
\begin{eqnarray*}
\tau^K(\cf) &=&  \chi(X,\lambda_{-1}(T^* X)\otimes \cf) \\
&=&\int_X \td(TX)\cup\Ch(\lambda_{-1}(T^* X)\otimes\cf) \\
&=& \int_X \td(TX)\cup\Ch(\lambda_{-1}(T^* X))\cup \Ch(\cf)
\\ 
&=& \int_X \td(TX)\cup\Ch(\lambda_{-1}(T^* X))\cup \Ch(\cf)\\ 
&=& \int_X \ctop(TX)\cup \Ch(\cf)\\ 
&=& \tau^A(\Ch(\cf)),
\end{eqnarray*}
where the second equality holds by the 
Hirzebruch-Riemann-Roch Theorem (a special case of the
Grothendieck-Riemann-Roch Theorem),
the third because $\Ch$ preserves multiplication, and the fifth by
Equation (\ref{eq:UsefulMix}).
\end{proof}

\begin{rem}
Since $K(X)$ is a $\nq$-vector space, we will need to make sense of
expressions such as $\td(\frac{1}{n}[E])$, where $n$ is a positive
integer and $E$ is a rank $r$ vector bundle over $X$. Observe that
\[
\td([E]) = \td\left(\bigoplus_{i=1}^n \frac{1}{n} [E]\right) = 
\left(\td(\frac{1}{n} [E])\right)^n.
\]
Consider the formal power series $\Phi(t_1,\ldots,t_r)$ in
$\nq[\![t_1,\ldots,t_r]\!]$ defined by 
\[
\Phi(t_1,\ldots,t_r) := \prod_{i=1}^r \phi(t_i).
\]
In particular, $\td([E]) = \Phi(a_1,\ldots,a_r)$. Since
$\Phi(t_1,\ldots,t_r)$ is equal to $\vac$ plus higher order terms, we
can define $\Phi^{\frac{1}{n}}(t_1,\ldots,t_r)$ to be the unique
formal power series in $\nq[\![t_1,\ldots,t_r]\!]$ equal to $\vac$ plus
higher order terms such that
\[
(\Phi^{\frac{1}{r}}(t_1,\ldots,t_r))^r = \Phi(t_1,\ldots,t_r).
\]
We define
\[
\td^{\frac{1}{r}}([E]) := \Phi^{\frac{1}{r}}(a_1,\ldots,a_r).
\]
\end{rem}

\section{$G$-graded $G$-modules and $G$-(equivariant) Frobenius algebras}

In this section we introduce some algebraic structures which we will need
throughout the rest of the paper.

\begin{df} \label{df:module} A
   $G$-graded vector space $\ch := \bigoplus_{m\in
    G} \ch_{m}$ endowed with the structure of a 
  $G$-module by isomorphisms $\rho(\gamma): \ch\irightarrow \ch$ for
  all $\gamma$ in $G$ is said to be a \emph{$G$-graded $G$-module} if
  $\rho(\gamma)$ takes $\ch_{m}$ to 
$\ch_{\gamma m\gamma^{-1}}$
 for
  all $m$ in $G$.
 \end{df}

$G$-graded $G$-modules form a category whose objects are $G$-graded
$G$-modules and whose morphisms are homomorphisms of $G$-modules which
respect the $G$-grading. Furthermore, the dual of a $G$-graded
$G$-module inherits the structure of a $G$-graded $G$-module.

Let us adopt the notation that $v_m$ is a vector in $\ch_m$ for any
$m\in G$.

\begin{df}\label{df:preGFA}
  A tuple $((\ch,\rho),\multipl,\vac,\eta,\tau)$ is said to be a
  \emph{pre-$G$-(equivariant) Frobenius algebra} provided that the following  properties hold: 
\begin{enumerate}
\item ($G$-graded $G$-module) $(\ch,\rho)$ is a (possibly infinite-dimensional) $G$-graded $G$-module.
\item (Self-invariance) For all $\gamma$ in $G$,
  $\rho(\gamma):\ch_{\gamma}\to\ch_{\gamma}$ is the identity map.
\item 
($G$-graded 
Pairing) $\eta$ is a symmetric, (possibly degenerate) bilinear form on
  $\ch$ such that $\eta(v_{m_1},v_{m_2})$ is nonzero only if $m_1 m_2
  = 1$.
\item ($G$-graded Multiplication) The binary product $(v_{1},
  v_{2})\mapsto v_{1}\multipl v_{2}$, called the \emph{multiplication} on
  $\ch$, preserves the $G$-grading (i.e., the multiplication is a map
  $\ch_{m_1}\otimes \ch_{m_2}\to\ch_{m_1 m_2}$) and is distributive
  over addition.
\item (Associativity) The multiplication is associative; i.e.,
\[
(v_{1}\multipl v_{2})\multipl v_{3} = v_{1}\multipl (v_{2}\multipl v_{3})
\]
for all $v_{1}$, $v_{2}$, and $v_{3}$ in $\ch$.
\item (Braided Commutativity) The multiplication is invariant with
  respect to the braiding:
\[
v_{m_1}\multipl v_{m_2} = (\rho(m_1) v_{m_2})\multipl v_{m_1}
\]
for all $m_i \in G$ and all $v_{m_i}\in \ch_{m_i}$ with $i=1,2$.
\item ($G$-equivariance of the Multiplication)
\[ (\rho(\gamma) v_{1})\multipl (\rho(\gamma) v_{2}) =
\rho(\gamma)(v_{1}\multipl v_{2})
\] for all $\gamma$ in $G$
and all $v_1, v_2 \in \ch$.
\item ($G$-invariance of the 
Pairing)
\[  \eta(\rho(\gamma) v_{1},\rho(\gamma) v_{2}) = \eta(v_{1}, v_{2})
\]
for all $\gamma$ in $G$ and all $v_1, v_2 \in \ch$.
\item (Multiplicative Invariance of the 
Pairing)
\[
\eta(v_{1}\multipl v_{2},v_{3}) = \eta(v_{1}, v_{2}\multipl v_{3})
\]
for all $v_1,v_2,v_3 \in \ch$.
\item \label{identity} ($G$-invariant Identity) The element $\vac$ in
$\ch_1$ is the identity element under the multiplication, and it
satisfies
\[
\rho(\gamma)\vac = \vac
\]
 for all $\gamma$ in $G$.
\item\label{equivtrace} ($G$-equivariant Trace Element) The \emph{trace
      element} $\tau$ is a collection $\{\,\tau_{a,b}\,\}_{a,b\in G}$ 
      of components $\tau_{a,b} \in \ch^*$, 
      such that $\tau_{a,b}(v_m)$ is nonzero
only if $m = [a,b]$, and is $G$-equivariant, i.e.,
\[
\tau_{\gamma a \gamma^{-1},\gamma b\gamma^{-1}} \circ\rho(\gamma) = \tau_{a,b}
\]
for all $a,b,\gamma$ in $G$.
\item (Trace Axiom)\label{NewTraceAxiom} For all $a,b$ in $G$, the trace
      element $\tau$ satisfies 
\[
\tau_{a,b} = \tau_{a b a^{-1},a^{-1}}.
\]
\end{enumerate}
We define the \emph{characteristic element} $\ttau$ in $\ch^*$
 to be
\[
\ttau := \frac{1}{|G|}\sum_{a,b\in G} \tau_{a,b},
\]
and we call the element $\ttau(\vac)\in \nq$ the  \emph{characteristic of the pre-$G$-Frobenius
algebra}.

A \emph{pre-Frobenius algebra} is a pre-$G$-Frobenius algebra with a
trivial group $G$. In this case, the trace element and characteristic element are
equal. 
\end{df}

\begin{rem}\label{rem:GInvTau}
The $G$-equivariance of $\tau$ insures that the characteristic element $\ttau$ is
$G$-invariant, i.e., 
\begin{equation}
\ttau\circ\rho(\gamma) = \ttau.
\end{equation}
\end{rem}

\begin{rem}
Any pre-$G$-Frobenius algebra $\ch$, has a pre-Frobenius subalgebra $\ch_1$
with a $G$-action which preserves the multiplication, unity, 
pairing,
and
trace element.
\end{rem}

\begin{rem}
  One can readily generalize the above definition to a 
pre-$G$-Frobenius
  superalgebra by introducing an additional $\nz/2\nz$-grading and
  inserting signs in the usual manner.
\end{rem}

\begin{df}\label{df:GFA}
For a tuple $((\ch,\rho),\multipl,\vac,\eta)$
satisfying all of the properties of a pre-$G$-Frobenius algebra which do not
involve the trace element and where $\ch$ is finite dimensional, we define the \emph{canonical
trace} to be 
\begin{equation}\label{eq:FAtrace}
\tau_{a,b}(v) := \tr_{\ch_{a}}(L_v\circ\rho(b))
\end{equation}
for all $a,b$ in $G$ and $v$ in $\ch_{[a,b]}$, where $L_v$ denotes left
multiplication by $v$.

We define a \emph{$G$-Frobenius algebra \cite{Ka,Ka2,Tu}} to be a tuple $((\ch,\rho),\multipl,\vac,\eta)$, such that $\ch$ is finite dimensional, the metric $\eta$ is nondegenerate, and such that the tuple, together with the canonical trace, forms a pre-$G$-Frobenius algebra.
\end{df}

\begin{rem}
The trace axiom (axiom~(\ref{NewTraceAxiom}))
for a $G$-Frobenius algebra $((\ch,\rho),\multipl,\vac,\eta)$
with the canonical trace is easily seen to be equivalent to the more familiar form
\begin{equation}\label{eq:OldTraceAxiom}
\tr_{\ch_{a}}(L_v\circ \rho(b)) = \tr_{\ch_{b}}(\rho(a^{-1})\circ L_v)
\end{equation}
for all $a,b$ in $G$ and $v$ in
  $\ch_{[a,b]}$, where $L_v$ denotes left multiplication by $v$.
\end{rem}

A $G$-Frobenius algebra with trivial group $G$ is nothing more than a Frobenius algebra.
Moreover, in this case the canonical trace of the trivial $G$-Frobenius algebra reduces to the canonical trace of the Frobenius algebra.

Later in the paper we will construct a \emph{stringy Chern character} which
maps the pre-$G$-Frobenius algebra of stringy K-theory to the
pre-$G$-Frobenius algebra of the stringy Chow ring.  We will see that, as in the case of the ordinary Chern character, the stringy
Chern character preserves all of the structure of a pre-$G$-Frobenius algebra
except the 
pairing.
This inspires the following definition:
\begin{df}
\label{df:allometric-GFA}
An \emph{allometric isomorphism} $\phi:((\ch,\rho),\multipl,\eta,\vac,\tau) \rTo
((\ch',\rho'),\star,\eta',\vac',\tau')$ of pre-$G$-Frobenius algebras is a
$G$-equivariant isomorphism of unital algebras that does not necessarily
preserve the 
pairing
but does preserve 
the trace element: 
$$\phi(\rho(m)v) = \rho(m)\phi(v)$$
for all $m\in G$ and all $v\in \ch$, 
and $$\phi^*\tau' = \tau.$$
\end{df}

\begin{df} Let $(\ch,\rho)$ be a $G$-graded $G$-module. Let
  $\pi_G : \ch\to\ch$ be the averaging map
\[
\pi_G(v) := \frac{1}{|G|}\sum_{\gamma\in G} \rho(\gamma) v
\]
for all $v$ in $\ch$. Let $\chb$ be the image of $\pi_G$. The vector
space $\chb$ is called the \emph{space of $G$-coinvariants of $\ch$,}
and it inherits a grading by the set $\Gb$ of conjugacy classes of
$G$:
\[
\chb = \bigoplus_{\gb\in\Gb} \chb_{\gb}.
\]
Since the group $G$ is finite, the space $\chb$ is equal to the space $\ch^G$ of $G$-invariants of $\ch$.

For any bilinear form
$\eta$ on $\ch$, we define $\etab$ to be the
restriction of the 
bilinear form
 $\frac{1}{|G|} \eta$ to $\chb$. 
Finally, define the trace element $\taub$ on $\chb$ to be the restriction
of the characteristic element $\ttau$ to $\chb$.
\end{df}

We have the following proposition.
\begin{prop}\label{prop:InducedPreGFA}
If the tuple $((\ch,\rho),\multipl,\vac,\eta,\tau)$ is a pre-$G$-Frobenius
algebra, then its $G$-coinvariants $(\chb,\multipl,\vac,\etab,\taub)$ form a
pre-Frobenius algebra, where $\multipl$ is induced from $\ch$. Moreover, if the
tuple $((\ch,\rho),\multipl,\vac,\eta)$ is a $G$-Frobenius algebra
with canonical trace element $\tau$, as defined in Equation
(\ref{eq:FAtrace}), then its ring of $G$-coinvariants $(\chb,\multipl,\vac,\etab)$ is a Frobenius algebra 
whose induced trace element $\taub$ is equal to its canonical trace element,
i.e.,
\begin{equation}\label{eq:Canonicaltaub}
\taub(\overline{v}) = \tr_{\chb}(L_{\overline{v}})
\end{equation}
for all $\overline{v}$ in $\chb$. 
\end{prop}
\begin{proof}
All that must be shown is Equation~(\ref{eq:Canonicaltaub}).
For all 
$m$ in $G$, and for all
$v_m$ in $\ch_m$, we have
\begin{eqnarray*}
\tr_{\chb}(L_{\pi_G(v_m)}) &= & \tr_\ch(L_{\pi_G(v_m)}\circ\pi_G) \\
&=&  \frac{1}{|G|^2}\sum_{\beta,\gamma\in G}
\tr_\ch(L_{\rho(\gamma)v_{m}}\circ\rho(\beta)) \\
&=& \frac{1}{|G|^2}\sum_{\beta,\gamma\in G}
\tr_\ch(L_{\rho(\gamma)v_{m}}\circ\rho(\gamma)\circ\rho(\gamma^{-1})\circ
\rho(\beta))  
\\ 
&=& \frac{1}{|G|^2}\sum_{\beta,\gamma\in G}
\tr_\ch(\rho(\gamma)\circ L_{v_{m}}\circ\rho(\gamma^{-1})\circ\rho(\beta))
\\ 
&=& \frac{1}{|G|^2}\sum_{\beta,\gamma\in G}
\tr_\ch(L_{v_{m}}\circ\rho(\gamma^{-1})\circ\rho(\beta)\circ \rho(\gamma))
\\ 
&=& \frac{1}{|G|^2}\sum_{\beta,\gamma\in G}
\tr_\ch(L_{v_{m}}\circ \rho(\gamma^{-1} \beta\gamma))\\
&=& \frac{1}{|G|^2}\sum_{b,\gamma\in G}
\tr_\ch(L_{v_{m}}\circ \rho(b))\\
&=& \frac{1}{|G|}\sum_{b\in G}
\tr_\ch(L_{v_{m}}\circ \rho(b)),
\end{eqnarray*}
where the second equality follows from
the definition of $\pi_G$,
the fourth from
the $G$-equivariance of the multiplication, and the fifth 
from
the cyclicity of the trace. Now, for all $a$ and $b$ 
in $G$, let
$\phi_a:\ch_a\rTo\ch_{m b a b^{-1}}$ be the restriction of
$L_{v_m}\circ\rho(b)$
 to $\ch_a$. The map $\phi_a$ preserves $\ch_a$ 
if and only if $m b a b^{-1} = a$ or, equivalently, if $m = [a,b]$. 
Furthermore, $\phi_a$ only contributes to the trace 
$\tr_\ch(L_{v_{m}}\circ \rho(b))$  when $m=[a,b]$. Therefore,
\[
\frac{1}{|G|}\sum_{b\in G} \tr_\ch(L_{v_{m}}\circ \rho(b)) =
\frac{1}{|G|}\sum_{a,b} \tr_{\ch_a}( L_{v_m}\circ \rho(b)) =
\frac{1}{|G|}\sum_{a,b} \tau_{a,b}(v_m),
\]
where the last two sums are over all $a, b\in G$ such that $[a,b] = m$. 
\end{proof}

\section{The stringy Chow ring 
and stringy K-theory 
of a variety with $G$-action}
\label{sec:Chow}\label{sec:K}

In this section, we discuss the main properties of the 
stringy Chow ring $\ca(X,G)$ and the stringy K-theory $\ck
(X,G)$ of a smooth, projective variety with an action of a finite group $G$.

As discussed in the introduction (Section~\ref{sec:summary}), the vector 
spaces underlying the 
stringy Chow ring and stringy K-theory of $X$ are just the usual Chow ring and K-theory, respectively, of the \emph{inertia variety}
$$I_G(X) := \coprod_{m\in G} X^m \subseteq X\times G,$$
where $X^m :=
\{ (x,m) | \rho(m) x = x \}$ with its induced $G$-action. 
Again, the reader should beware that the inertia 
\emph{variety} is not the same as the inertia
\emph{orbifold} $$[I_G(X)/G] = \coprod_{[\![ g ]\!]}
[X^{g}/Z_G(g)]$$
of \cite{CR1, AGV}, which is the stack quotient of
the inertia variety $I_G(X)$ by the action of $G$.

Recall (see Definition~\ref{df:kage}) that one of the key elements in the
construction of both the stringy multiplication and the stringy Chern
character is the element $\cs\in I_G(X)$, defined as  
\begin{equation}
\cs_m:= \res{\cs}{X^m}:= \bigoplus_{k=0}^{r-1}\frac{k}{r}W_{m,k}.
\end{equation} 

The $G$-equivariant involution $\sigma:X^m\to X^{m^{-1}}$ yields a
$G$-equivariant isomorphism $\sigma^*:W_{m^{-1}}\to W_m$ for all $m$
in $G$.  If $m$ acts by multiplication by $\zeta^k$, then $m^{-1}$
acts by $\zeta^{r-k}$, so we have
\begin{equation}
\sigma^*W_{m^{-1},0} = W_{m,0}
\end{equation}
and
\begin{equation}\label{eq:NSigma}
\sigma^*W_{m^{-1},k} = W_{m,r-k}
\end{equation}
for all $k \in \{1,\ldots,r-1\}$. Consequently, the induced map
$\sigma^*:K(X^{m^{-1}})\to K(X^m)$ satisfies
\begin{equation}\label{eq:SSigma}
\cs_m\oplus \sigma^*\cs_{m^{-1}} = N_m,
\end{equation}
since the normal bundle, $N_m$, of $X^m$ in
$X$ satisfies the equation $N_m = W_m\ominus W_{m,0}$. 

The virtual rank of $\cs_m$ on a connected component $U$ of $X^m$ is
the age $a(m,U)$ (see Definition~\ref{df:kage}). Taking the virtual rank
of both sides of Equation~(\ref{eq:SSigma}) yields the well-known equation 
\begin{equation}\label{eq:AgePlusAgeInv}
a(m,U) + a(m^{-1},U) = \codim(U\subseteq X).
\end{equation}
This supports the interpretation of $\cs_m$ as a K-theoretic version of the age.

Recall that one may use the age to define a rational grading on $\ca(X,G)$. 
\begin{df}
For all $m$ in $G$, all connected components $U$ of $X^m$, and all
elements $v_m$ in $A^p(U)\subseteq \ca_m(X)$, for $p$ the usual
integral degree in the Chow ring, we define a $\nq$-grading which we call the 
\emph{stringy grading}
$\grq{v_m}$ on $\ca_m(X)$ by
\begin{equation}\label{eq:DefGrading}
\grq{v_m} := a(m,U)+p.
\end{equation}
\end{df}

\begin{rem}\label{rem:integrality}
  Sometimes the $\nq$-grading just happens to be integral. For example, 
  if $X$ is
  $n$-dimensional and its canonical bundle $K_X$ has a
  nowhere-vanishing section $\Omega$, then for all $m$ in $G$, we have
\[ \rho(m)^*\Omega = \exp(2\pi i a(m)) \Omega. \]
Thus, if $G$ preserves $\Omega$, then $a(m)$ must be an integer.

A special case is when $X$ is $2n$-dimensional, possessing a (complex
algebraic) symplectic form $\omega$ in $\bigwedge^2 T^* X$. This can
arise if $X$ happens to be a hyper-K\"ahler manifold. If, in
addition,  $G$ preserves $\omega$, then $G$ preserves the nowhere
vanishing section $\omega^n$ of $K_X$. In this case, for all $m$ in
$G$ and for every connected component $U$ of $X^m$, the associated age
\cite{Kal} is the integer
\[
a(m,U) = \frac{1}{2}\codim(U\subseteq X).
\]
\end{rem}

\begin{rem}
Unlike the stringy Chow ring, the ring $\ck(X,G)$
lacks a $\nq$-grading. This should not be surprising, as even ordinary
K-theory lacks a grading by ``dimension.'' In particular, the virtual rank does not enjoy the same good properties in K-theory that grading by codimension has in the Chow ring.   
\end{rem}

The multiplication in the string Chow ring and stringy K-theory were already defined in Definitions~\ref{df:chow-mult} and \ref{df:K-mult}, but to see that these form pre-$G$-Frobenius algebras, we also need to define their trace element.
\begin{df} The \emph{trace element $\taua$ of $\ca(X,G)$} is a collection of
\emph{components} $\{ \tau_{a,b} \}_{a,b\in G}$, where
$\tau_{a,b}\in\ca(X,G)^*$ is defined to be
\begin{equation}\label{eq:taua}
\tau_{a,b}(v_m) := 
\begin{cases} 
\int_{\Xab} \res{v_m}{\Xab}\cup\ctop(T\Xab\oplus
                         \res{\cs_{[a,b]}}{\Xab}) & \text{if $m = [a,b]$}\\
0 &\text{if $m\not=[a,b]$}
\end{cases}
\end{equation}
for all $a$, $b$, m in $G$ and $v_m$ in $\ca_m(X)$. 
\end{df}
We define the trace element of stringy K-theory similarly.
\begin{df} The \emph{trace element $\tauk$ of $\ck(X,G)$} is a collection of \emph{components} $\{ \tau_{a,b} \}_{a,b\in G}$, where
$\tau_{a,b}\in\ck(X,G)^*$ is defined to be
\begin{equation}\label{eq:tauk}
\tau_{a,b}(\cf_m) := 
\begin{cases}
\chi(\Xab, \res{\cf_m}{\Xab}\cup\lambda_{-1}(T\Xab\oplus
                         \res{\cs_{[a,b]}}{\Xab})^*) & \text{if $m=[a,b]$}\\
0 &\text{if $m\not=[a,b]$}
\end{cases}
\end{equation}
for all $a$, $b$, $m$ in $G$, and $\cf_m$ in $\ck_m(X)$.  Here  $\chi$ denotes the Euler characteristic. 
\end{df}

We will show in Section~\ref{subsec:TraceAxiom} that the element $T\Xab\oplus
\res{\cs_{[a,b]}}{\Xab}$ in $K(\Xab)$ can be represented by a vector bundle. 
 Similarly, we will prove in Section~\ref{sec:FG=us} that $\crr(\bm)$ is represented by a vector bundle.  In the meantime, we will use that fact to prove the following results. 
\begin{thm}\label{thm:ChowGFA}
Let $X$ be a smooth, projective variety with an action of a finite
group $G$.
\begin{enumerate}
\item \label{thm:ChowGFAOne} The tuple
      $((\ca(X,G),\rho),\multipl,\vac,\etaa,\taua)$
       is a pre-$G$-Frobenius algebra.
\item \label{thm:ChowGFATwo} $\grq{\vac} = 0$.
\item \label{thm:ChowGFAThree} The multiplication respects the
  $\nq$-grading, i.e., for all homogeneous elements $v_{m_i}$ in
  $\ca_{m_i}(X)$, for $i=1,2$, we have
\[
\grq{v_{m_1}\multipl v_{m_2}} = \grq{v_{m_1}} + \grq{v_{m_2}}.
\]
\item \label{thm:ChowGFAFour}The 
pairing
has a definite $\nq$-grading,
  i.e., for all homogeneous elements $v_{m_1}$ in $\ca_{m_1}(X)$ and $v_{m_2}$ in $\ca_{m_2}(X)$ we have
$\etaa(v_{m_1},v_{m_2}) = 0$
unless $m_1 m_2 = 1$ and 
\begin{equation}\label{eq:etaatwo}
\grq{v_{m_1}} + \grq{v_{m_2}} = \dim X.
\end{equation}
\item \label{thm:ChowGFAFive}
The components $\{\tau_{a,b}\}$ of $\taua$,  satisfy
$\tau_{a,b}(v_m) = 0$
unless 
$\grq{v_m} = 0$
and $m = [a,b]$.
\end{enumerate}
\end{thm}

\begin{thm}\label{thm:KGFA}
  Let $X$ be a smooth, projective variety with an action of a finite
  group $G$.  The tuple $((\ck(X,G),\rho),\multipl,\vac,\etak,\tauk)$ is a
  pre-$G$-Frobenius algebra, where the trace element $\tauk$ is defined by
Equation (\ref{eq:tauk}).
\end{thm}

For both Theorems (\ref{thm:ChowGFA} and \ref{thm:KGFA}),
the only nontrivial parts are the
trace axiom and the associativity of multiplication.  These are proved in 
Lemmas~\ref{thm:traceaxiom} and ~\ref{thm:SecondAssoc}, respectively. 

\begin{ex}
Consider the case where  $m_i = 1$ for some $i=1,2,3$. 
In this case the bundle $\crr$ on $X^\bm$
is trivial. 

If $m_1 = 1$
and $m_2 m_3 = 1$, then the stringy multiplication is given by the
restriction to $X^{m_3}$ of the ordinary multiplication in ordinary
K-theory, i.e.,  
\begin{equation}\label{eq:m1Equals1}
\cf_{m_1=1} \multipl \cf_{m_2} = \res{\cf_{m_1}}{X^{m_3}} \otimes
\sigma^*\cf_{m_2}.  
\end{equation}
A similar result holds if $m_2 = 1$ and $m_1 m_3 = 1$. In particular, this
means that stringy multiplication on the untwisted sector $\ck_1(X)$
coincides with the ordinary multiplication on $\ck_1(X)$.

More interesting is the case where $m_3 = 1$ and $m_1 m_2 = 1$.  In this case, we have
\begin{equation}\label{eq:m3Equals1}
\cf_{m_1}\multipl\cf_{m_2} = \ec_{m_3 *}( \e_{m_1}^*\cf_{m_1}\otimes
\e_{m_2}^*\cf_{m_2}).
\end{equation}
Here, even though the bundle $\crr(\bm)$ is trivial, the stringy multiplication is nontrivial, since the map
$\ec_{m_3}$  will generally be between varieties of 
different dimensions.
\end{ex}

\begin{rem}
If $a=b=1$ in $\ca(X,G)$, then for all $v_1$ in $\ca_1(X)$, we have 
\begin{equation}
\tau_{1,1}(v_1) = \int_{X} v_1\cup\ctop(TX).
\end{equation}
Therefore, the component $\tau_{1,1}$ of the trace element on the untwisted sector
$\ca_1(X)$ of the stringy Chow ring agrees with the trace element of the
ordinary Chow ring $A^\bullet(X)$. 
\end{rem}

\begin{rem}
The characteristic of $((\ca(X,G),\rho),\multipl,\vac,\etaa,\taua)$ is
\begin{equation}\label{eq:Characteristic}
\ttau(\vac) = \frac{1}{|G|}\sum_{ a b = b a} \chi(X^{\left< a, b\right>}),
\end{equation}
where the sum is over all commuting pairs $a$, $b$ in $G$, and 
\[
\chi(X^{\left< a, b\right>}) = \int_{[X^{\left< a, b\right>}]}
\ctop(TX^{\left< a, b\right>})
\]
is the usual Euler characteristic. This expression~(\ref{eq:Characteristic}) coincides with the ``stringy
Euler characteristic'' introduced by physicists \cite{DHVW} (see also
\cite{AS}).
\end{rem}

\section{Associativity and the trace axiom}\label{sec:AssocTrace}

In this section, we use the fact that the element $\crr$ defined in Equation~(\ref{eq:R}), is a vector bundle
(Proved in Corollary~\ref{cor:RVB})
to give an elementary proof of associativity and the trace axiom for both the
stringy Chow ring and stringy K-theory. 

\subsection{Associativity}
Let us recall some excess intersection theory. Consider smooth,
projective varieties $V$, $Y_1$, $Y_2$, and $Z$ which form the
following Cartesian square
\begin{equation}
\begin{diagram}
V  & \rTo^{i_1} & Y_1 \\
\dTo^{j_2} & & \dTo^{j_1}
\\
Y_2  & \rTo^{i_2} & Z
\end{diagram},
\end{equation}
where $i_1$, $i_2$ are regular embeddings and $j_1,j_2$ are morphisms
of schemes.  

Let $E(V,Y_1,Y_2)\to V$ be the \emph{excess normal (vector) bundle},
which is the cokernel of the 
map $N_{V/Y_1}\to \res{N_{Y_2/Z}}{V}$, where $N_{V/Y_1}$ and $N_{Y_2/Z}$ are the
normal bundles of $V$ in $Y_1$ and of $Y_2$ in $Z$, respectively. In
$K(V)$ one thus obtains the equality
\begin{equation}
[E(Z,Y_1,Y_2)] = \res{TZ}{V}\ominus\res{TY_1}{V}\ominus\res{TY_2}{V}\oplus TV.
\end{equation}

Under these hypotheses, the following theorem holds (see Theorems 1.3
and 1.4 in \cite[Chapter IV.1]{FuLa}).
\begin{thm} \label{thm:excessintersection}
For all $\cf$ in $K(Y_2)$ and $v$ in $A^\bullet(Y_2)$, 
\begin{equation}
j_1^* i_{2 *} \cf = i_{1 *}(\lambda_{-1}(E(Z,Y_1,Y_2)^*) \tensor j_2^*\cf)
\end{equation}
and
\begin{equation}
j_1^* i_{2 *} v = i_{1 *}(\ctop(E(Z,Y_1,Y_2)) \cup j_2^*v).
\end{equation}
\end{thm}

The previous theorem gives rise to the following fact about $\crr$.

\begin{lm}\label{thm:crr}
Let $\bm := (m_1,\ldots,m_4)$ in $G^4$ such that $m_1 m_2 m_3 m_4 =
1$. Let $X^\bm$ consist of those points in $X$ which are fixed by
$m_i$ for all $i \in \{1,\ldots,4\}$. The following equation holds in
$K(X^\bm)$:
\begin{multline}\label{eq:MagicRelation}
\res{\crr(m_1,m_2,(m_1 m_2)^{-1})}{X^\bm}\oplus
\res{\crr(m_1 m_2,m_3,m_4)}{X^\bm} \oplus 
E_{m_1,m_2}  \\
=\res{\crr(m_1,m_2 m_3,m_4)}{X^\bm}\oplus
\res{\crr(m_2,m_3,(m_2 m_3)^{-1})}{X^\bm}\oplus
E_{m_2,m_3},
\end{multline}
where
\begin{equation}\label{eq:Em1m2}
E_{m_1,m_2} := E(X^{m_1 m_2},X^{\langle m_1,m_2 \rangle},X^{\langle m_1 m_2,m_3 \rangle})
\end{equation}
and
\begin{equation}\label{eq:Em2m3}
E_{m_2,m_3} := E(X^{m_2 m_3},X^{\langle m_1,m_2 m_3 \rangle},X^{\langle m_2, m_3 \rangle}).
\end{equation}
Furthermore, both sides of Equation~(\ref{eq:MagicRelation}) are equal in
$K(X^\bm)$ to
\begin{equation}\label{eq:crtFour}
TX^\bm \ominus \res{TX}{X^\bm} \oplus \bigoplus_{i=1}^4 \res{\cs_{m_i}}{X^\bm}.
\end{equation}
\end{lm}
\begin{proof}
Plug in the definitions of the excess normal bundles and the formula for the
obstruction bundle $\crr$ from Equation~(\ref{eq:Magic}), then apply 
Equation~(\ref{eq:SSigma}) and simplify the result. One discovers that both
the  right hand and left hand sides of Equation~(\ref{eq:MagicRelation}) are
equal in $K(X^\bm)$ to Equation~(\ref{eq:crtFour}).
\end{proof}

\begin{rem}\label{rem:FourPoint}
  For the reader familiar with the $G$-stable maps of \cite{JKK}, we
  note that the element $TX^\bm \ominus \res{TX}{X^\bm} \oplus
  \bigoplus_{i=1}^4 \res{\cs_{m_i}}{X^\bm}$ in
  Equation~(\ref{eq:crtFour}) may be interpreted as the fiber of the
  obstruction bundle over $\{ q \}\times X^\bm$ in
  $\xi_{0,4}(\bm)\times X^\bm =\xi_{0,4} (X,0,\bm)$, where $q$ is any
  point in $\xi_{0,4}(\bm)$.  This can be seen by an argument similar
  to that in the proof of \cite[Prop 6.21]{JKK}.
\end{rem}

\begin{lm}\label{thm:SecondAssoc}
  Let $X$ be a smooth, projective variety with an action of a finite
  group $G$. The multiplications in stringy K-theory
  $((\ck(X,G),\rho),\multipl,\vac,\etak)$ and in the stringy Chow ring
  $((\ca(X,G),\rho),\multipl,\vac,\etaa)$ are both associative.
\end{lm}
\begin{proof}
  Consider $\bm = (m_1,m_2,m_3,m_4)$ in $G^4$ such that $m_1 m_2 m_3
  m_4 = 1$. If $E_{m_1,m_2}$ and $E_{m_2,m_3}$ are defined as in
  Equations (\ref{eq:Em1m2}) and~(\ref{eq:Em2m3}), then the following
  equalities hold:
\begin{multline}\label{eq:ChowAssoct}
\res{\ctop(\crr(m_1,m_2,(m_1 m_2)^{-1}))}{X^\bm} 
\cup \res{\ctop(\crr(m_1 m_2,m_3,m_4))}{X^\bm}  
\cup \ctop(E_{m_1,m_2}) = \\
\res{\ctop(\crr(m_1,m_2 m_3,m_4))}{X^\bm} 
\cup \res{\ctop(\crr(m_2,m_3,(m_2 m_3)^{-1}))}{X^\bm} 
\cup \ctop(E_{m_2,m_3})
\end{multline}
and
\begin{multline}\label{eq:KAssoct}
\res{\lambda_{-1}(\crr(m_1,m_2,(m_1 m_2)^{-1})^*)}{X^\bm} \tensor
\res{\lambda_{-1}(\crr(m_1 m_2,m_3,m_4)^*)}{X^\bm}  \tensor 
\lambda_{-1}(E_{m_1,m_2}^*) = \\
\res{\lambda_{-1}(\crr(m_1,m_2 m_3,m_4)^*)}{X^\bm}  \tensor
\res{\lambda_{-1}(\crr(m_2,m_3,(m_2 m_3)^{-1})^*)}{X^\bm}  \tensor
\lambda_{-1}(E_{m_2,m_3}^*).
\end{multline}
Equation~(\ref{eq:ChowAssoct}) follows by taking the top Chern class
of both sides of Equation~(\ref{eq:MagicRelation}) and then using
multiplicativity of $\ctop$.  Equation~(\ref{eq:KAssoct}) follows by
taking the dual of Equation~(\ref{eq:MagicRelation}), applying
$\lambda_{-1}$, and then using multiplicativity of $\lambda_{-1}$.

Associativity will follow from Equations~(\ref{eq:ChowAssoct}) and
(\ref{eq:KAssoct}) and the definitions of the multiplications
as follows.

Let $m_+ = (m_1m_2)^{-1}$ and $m_- = (m_1m_2)$. Consider the following
diagram:
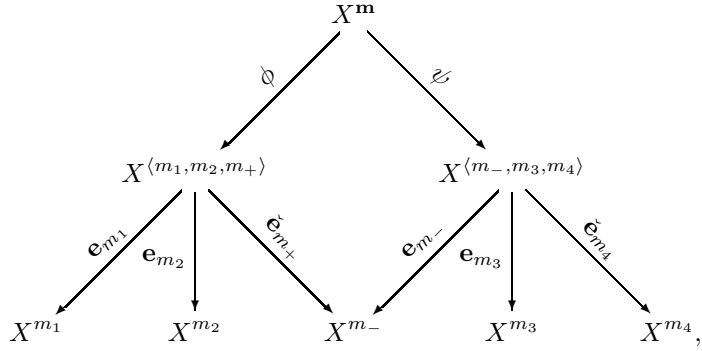
\begin{diagram}
        & & & & X^{\bm} & & & & \\ & & &\ldTo^{\phi} & & \rdTo^{\psi}
        & & & \\ & & X^{\langle m_1, m_2, m_+ \rangle}& & & &
        X^{\langle m_-, m_3,m_4\rangle} \\ & \ldTo^{\e_{m_1}} &
        \dTo^{\e_{m_2}} & \rdTo^{\ec_{m_+}} & & \ldTo^{\e_{m_-}} &
        \dTo^{\e_{m_3}} & \rdTo^{\ec_{m_4}} & \\
X^{m_1} & & X^{m_2} & & X^{m_-}& & X^{m_3} & & X^{m_4},
\end{diagram}
where $\phi$ and $\psi$ are the obvious inclusions.  Note that the
diamond in the middle is Cartesian and that the usual inclusions
$\epsilon_i:X^{\bm} \to X^{m_i}$ factor as
\begin{align}\label{eq.factor}
  \epsilon_1 &= \e_{m_1} \circ \phi \quad & \epsilon_2 = \e_{m_2}
  \circ \phi\\ \epsilon_3 &= \e_{m_3} \circ \psi \quad & \epsilon_4 =
  \e_{m_4} \circ \psi. \label{eq.factortwo}
\end{align}
Finally, we define $$\check{\epsilon}_4 = \sigma\circ \epsilon_4 =
\ec_{m_4} \circ \psi.$$ 

\medskip

For any $\cf_1 \in \ck_{m_1}$, $\cf_2 \in \ck_{m_2}$, $\cf_3 \in
\ck_{m_3}$, we have
\newlength{\LHSwidth}\setlength{\LHSwidth}{7em}
\newlength{\RHSoffset}\setlength{\RHSoffset}{5em}
\begin{equation*} 
\begin{split}
\makebox[\LHSwidth]{ $(\cf_1 \multipl \cf_2) \multipl \cf_3$ } & =
(\ec_{m_4})_*\Bigg(\e_{m_-}^*(\ec_{m_+})_* \left[ \e_{m_1}^* \cf_1
\tensor \e_{m_2}^* \cf_2 \tensor \lambda_{-1}( \crr(m_1, m_2,
m_+)^*)\right] \\ & \hspace*{\RHSoffset} {\tensor\, \e_{m_3}^* \cf_3
\tensor \lambda_{-1} (\crr(m_-, m_3,m_4)^*)\Bigg)}
\end{split}
\end{equation*}
\begin{equation*}
\begin{split}
  \rule{\LHSwidth}{0ex} 
&= (\ec_{m_4})_*\Bigg(\psi_* \Big( \phi^* \left[ \e_{m_1}^* \cf_1 
\tensor   \e_{m_2}^* \cf_2  
\tensor   \lambda_{-1}( \crr(m_1, m_2, m_+)^*) \right]     \\
  & \hspace*{\RHSoffset} \tensor \lambda_{-1}(E^*_{m_1,m_2}) \Big)
  \tensor \e_{m_3}^* \cf_3 \tensor \lambda_{-1} (\crr(m_-,
  m_3,m_4)^*)\Bigg)
\end{split}
\end{equation*}
\begin{equation*}
\begin{split}
\rule{\LHSwidth}{0ex}& =  (\ec_{m_4})_*\Bigg(\psi_* 
\Big( \phi^* \e_{m_1}^* \cf_1   
\tensor  \phi^* \e_{m_2}^* \cf_2  
\tensor  \phi^* \left( \lambda_{-1}( \crr(m_1, m_2, m_+)^*) \right)  \\ 
  & \hspace*{\RHSoffset} \tensor \lambda_{-1}(E^*_{m_1,m_2}) \tensor
  \psi^* \e_{m_3}^* \cf_3 \tensor \psi^*(\lambda_{-1} (\crr(m_-,
  m_3,m_4)^*))\Big)\Bigg) \\
\end{split}
\end{equation*}
\begin{equation*}
\begin{split}
\rule{\LHSwidth}{0ex} & = (\check{\epsilon}_4)_* \Bigg( \epsilon_1^*
\cf_1 \tensor \epsilon_2^* \cf_2 \tensor \phi^*\lambda_{-1}( \crr(m_1,
m_2, m_+)^*) \tensor \lambda_{-1}(E^*_{m_1,m_2}) \\
\end{split}
\end{equation*}
\begin{equation*}
\begin{split}
 & \hspace*{\RHSoffset} \tensor \epsilon_3^* \cf_3 \tensor
 \psi^*(\lambda_{-1} (\crr(m_-, m_3,m_4)^*))\Bigg) \\
\end{split}
\end{equation*}
\begin{equation}\label{eq:assoclhs}
\begin{split}
\rule{\LHSwidth}{0ex} &= (\check{\epsilon}_4)_* \Bigg( \epsilon_1^*
\cf_1 \tensor \epsilon_2^* \cf_2 \tensor \epsilon_3^* \cf_3 \tensor
\res{\lambda_{-1}( \crr(m_1, m_2, m_+)^*)}{X^{\bm}} \\
 & \hspace*{\RHSoffset} \tensor \res{\lambda_{-1} (\crr(m_-,
 m_3,m_4)^*)}{X^{\bm}} \tensor \lambda_{-1}(E^*_{m_1,m_2}) \Bigg),
\end{split}
\end{equation}
where the first equality is the definition, the second equality
follows from Theorem~\ref{thm:excessintersection}, the third equality
follows from the projection formula, and the fourth and fifth
equalities follow from Equations~(\ref{eq.factor}) and
(\ref{eq.factortwo}) and the definitions of $\psi$ and $\phi$.

A similar argument shows that the product $\cf_1 \multipl (\cf_2 \multipl
\cf_3) $ is given by
\begin{equation}\label{eq:assocrhs}
\begin{split}
\cf_1 \multipl (\cf_2 \multipl \cf_3) & = (\check{\epsilon}_4)_* \Big(
\epsilon_1^* \cf_1 \tensor \epsilon_2^* \cf_2 \tensor \epsilon_3^*
\cf_3 \tensor \res{\lambda_{-1}(\crr(m_1,m_2 m_3,m_4)^*)}{X^\bm} \\ &
\hspace*{5em} \tensor \res{\lambda_{-1}(\crr(m_2,m_3,(m_2
m_3)^{-1})^*)}{X^\bm}
 \tensor \lambda_{-1}(E_{m_2,m_3}^*)\Big).
\end{split}
\end{equation}
By equation (\ref{eq:KAssoct}) these two
expressions~(\ref{eq:assoclhs}) and (\ref{eq:assocrhs}) are equal, hence
associativity holds.
\end{proof}

\subsection{The trace axiom}\label{subsec:TraceAxiom}

We now prove the trace axiom in a similar way. Throughout this
section, we fix elements $a$ and $b$ in $G$ and let $m_1 := [a,b]$.
Let $\bmp := (\mmp_1, \mmp_2,\mmp_3) := ([a,b], b a b^{-1}, a^{-1})$.
Let $H := \langle a, b \rangle 
$ and let $H'
:= \langle \bmp \rangle \le H$ be the subgroup generated by the elements of $\bmp$.
Let $\crr(\bmp)$ denote the
element in $K(X^{H'})$ from Equation~(\ref{eq:R}).

Consider the commutative diagram
\begin{equation}
\label{eq:excessp}
\begin{diagram}
X^H  & \rTo^{\jp_2} & X^{H'} \\
\dTo^{\jp_1} & & \dTo^{\Deltap_2}
\\
X^{a}  & \rTo^{\Deltap_1} & X^{b a b^{-1}} \times X^{a^{-1}}.\\
\end{diagram}
\end{equation}
Here $\jp_1$ and $\jp_2$ are the obvious inclusion morphisms, $\Deltap_2$ is
the diagonal map, and $\Deltap_1$ is the composition of the morphisms
\[
X^a\rTo^{\Delta} X^a\times X^a\rTo^{\rho(b)\times\sigma} X^{b a b^{-1}}
\times X^{a^{-1}},
\]
where $\Delta$ is the diagonal map and 
$\rho(b)$ 
is the action of $b$. 
Let $\cep$ be the excess intersection bundle $E(X^{b a b^{-1}}\times X^{a^{-1}}, X^{H'}, X^a)$.

\begin{thm}\label{thm:GenusOne}
The following equality holds in $K(X^H)$:
\begin{equation}\label{eq:crrtGenusOne}
\jp_2^*\crr(\bmp)\oplus \cep =TX^{H} \oplus \res{\cs_{m_1}}{X^{H}}.
\end{equation}
\end{thm}

The previous theorem together with Corollary~\ref{cor:RVB} and the fact that
$\cep$ is an excess intersection (vector) bundle, yields the following.
\begin{crl}\label{crl:TraceVB}
$TX^H\oplus \res{\cs_{m_1}}{X^H}$ can be represented in $K(X^H)$ by a vector
bundle. 
\end{crl}

\begin{proof}(of Theorem~\ref{thm:GenusOne})
All equalities in this proof are understood to be in
$K(X^{H})$. Observe that
\begin{eqnarray*}
\jp_2^*\Deltap_2^*T(X^{b a b^{-1}}\times X^{a^{-1}}) &= &
\res{TX^{b a b^{-1}}}{X^{H}} \oplus \res{TX^{a^{-1}}}{X^{H}}\\
&=& 
 \res{\rho(b) (TX^{a})}{X^{H}} \oplus \res{\sigma^*TX^{a}}{X^{H}}
 \\
&=& \res{TX^{a}}{X^{H}} \oplus \res{TX^{a}}{X^{H}},
\end{eqnarray*}
where the third equality follows from the fact that 
$\rho(b)\times \sigma$ 
is an isomorphism. 
Plugging this into the definition of the excess
intersection bundle yields 
\[
\cep =
TX^{H}\oplus\res{TX^a}{X^{H}}\oplus\res{TX^{a}}{X^{H}} 
\ominus \res{TX^a}{X^{H}} \ominus\res{TX^{H'}}{X^{H}}, 
\]
which simplifies to 
\begin{equation}\label{eq:excessbundlep}
\cep = TX^{H}\oplus\res{TX^a}{X^{H}}\ominus\res{TX^{H'}}{X^{H}}.
\end{equation}

On the other hand, Equation~(\ref{eq:Magic}) yields the equality
\begin{equation*}
\crr(\bmp) = \res{TX^{{H'}}}{X^{H}}\ominus \res{TX}{X^{H}}
\oplus \res{\cs_{m_1}}{X^{H}} \oplus \res{\cs_{b a b^{-1}}}{X^{H}} \oplus
\res{\cs_{a^{-1}}}{X^{H}}.
\end{equation*}
Together with the equality
\[
\res{\cs_{b a b^{-1}}}{X^{H}} = \res{\rho(b)(\cs_a)}{X^{H}}
\]
and Equation~(\ref{eq:SSigma}), we obtain
\begin{equation}\label{eq:obsp}
\crr(\bmp) = \res{TX^{H'}}{X^{H}}\ominus\res{TX^a}{X^{H}}\oplus
\res{\cs_{m_1}}{X^{H}}.
\end{equation}
Combining Equations~(\ref{eq:excessbundlep}) and~(\ref{eq:obsp}) yields the
identity 
\begin{equation}
\jp_2^*\crr(\bmp)\oplus \cep = TX^{H} \oplus \res{\cs_{m_1}}{X^{H}}.
\end{equation}

\end{proof}

\begin{rem}
If we make the replacement $(a,b)\mapsto (\ta, \tb) := (a b
a^{-1},a^{-1})$ everywhere in the Theorem~\ref{thm:GenusOne}, then since
$\langle a,b \rangle = \langle \ta,\tb \rangle$ and $m_1 = [a,b] =
[\ta,\tb]$, the
right hand side of  Equation~(\ref{eq:crrtGenusOne}) stays the same, while the
left hand side changes. Hence, one obtains an interesting equality between
the corresponding
left hand sides of Equation~(\ref{eq:crrtGenusOne}) before and after
making the substitution $(a,b)\mapsto (\ta, \tb)$. The resulting equality is
the analogue of Equation~(\ref{eq:MagicRelation}).
\end{rem}

\begin{rem}\label{rem:GenusOne}
For the reader familiar with $G$-stable maps, we note that the element
$TX^{{H}} \oplus \res{\cs_{m_1}}{X^{{H}}}$ from
Equation~(\ref{eq:crrtGenusOne}) is the restriction of the obstruction
bundle over $\xi_{1,1}(m_1,a,b)\times X^{{H}}$ to $\{ q \}\times
X^{H}$, where $q$ is any point in $\xi_{1,1}(m_1,a,b)$.  The details
of this are given in \cite[Prop. 6.21]{JKK}.
\end{rem}

\begin{lm}\label{thm:traceaxiom}
If $X$ is a smooth, projective variety with an action of a finite
group $G$, then stringy K-theory $((\ck(X,G),\rho),\multipl,\vac,\eta,\tauk)$
and the stringy Chow ring $((\ca(X,G),\rho),\multipl,\vac,\eta,\taua)$ 
satisfy the trace axiom for pre-$G$-Frobenius algebras.
\end{lm}
\begin{proof}
Consider the case of $\ck(X,G)$. Corollary~\ref{crl:TraceVB} implies that
$\lambda_{-1}(T\Xab\oplus\res{\cs_{[a,b]}}{\Xab})^*$ is
well-defined, hence the trace element $\tauk$ given in
Equation~(\ref{eq:tauk}) is well-defined. The trace axiom for a
pre-$G$-Frobenius algebra follows immediately from the observation that $m =
[a b a^{-1},a^{-1}] = [a,b]$ and $\langle a,b \rangle = \langle a b a^{-1},
a^{-1} \rangle$ for all $a,b$ in $G$.

The case of $\ca(X,G)$ is analogous.
\end{proof}

This completes the proof of Theorems~\ref{thm:ChowGFA} and \ref{thm:KGFA} that stringy Chow and stringy K-theory are pre-$G$-Frobenius algebras.

\section{The stringy Chern character}
\label{sec:chern}

As mentioned in the introduction, for general $G$ the ordinary Chern character fails
to be a ring homomorphism; however, this drawback can be overcome
through the introduction of the appropriate correction terms to give
what we call the \emph{stringy Chern character}
$\CCh:\ck(X,G)\to\ca(X,G)$.  The main purpose of this section is to 
prove the key
properties of $\CCh$, and especially to prove that the map $\CCh$ is an
allometric isomorphism for any smooth,  projective variety $X$ with
an action of a finite group $G$. When $G$ is the trivial group, $\CCh$
reduces to the usual Chern character mapping from ordinary K-theory to the
ordinary Chow ring of $X$.

Recall that for any smooth, projective variety $X$ with an action of
  $G$, we have defined (see Definition~\ref{eq:DefCCh})
  the \emph{stringy Chern character} $\CCh:\ck(X,G)\to \ca(X,G)$
  to be 
\begin{equation}
\CCh(\cf_m) := \Ch(\cf_m)\cup \td^{-1}(\cs_m)
\end{equation}
for all $m$ in $G$ and $\cf_m$ in $\ck_m(X)$, where $\cs_m$ is defined in
Equation~(\ref{eq:S}), $\td$ is the Todd class, and $\Ch$ is the ordinary
Chern character.

The main result of this section is the following 
theorem.
\begin{thm}\label{thm:ChernHomo}
  The stringy Chern character $\CCh:\ck(X,G)\to \ca(X,G)$ is an
allometric isomorphism. In particular, $\CCh$ is a
$G$-equivariant algebra isomorphism.
\end{thm}
\begin{proof}
To see that $\CCh$ is an isomorphism of $G$-graded
$G$-modules, note first that $\td$ is invertible (it is a series starting with $\vac$), so $\CCh$ is an isomorphism of $G$-graded vector spaces.  The
equivariance under the $G$-action follows from naturality properties
of $\td$, $\Ch$, the cup product, and $\cs_m$.

We now prove that $\CCh$ respects multiplication. We suppress the cup
and tensor product symbols to avoid notational clutter.  Let
$\cf_{m_i}$ belong to $\ck_{m_i}(X)$ for $i=1,2,3$, where $m_1 m_2
m_3 = 1$. Let $\e_{m_i}$ denote the inclusion $X^\bm\to X^{m_i}$ and
$\ec_{m_i}:=\sigma\circ\e_{m_i}:X^\bm\to X^{m_i^{-1}}$ for all
$i=1,2,3$.  We have
\begin{eqnarray*}
\CCh(\cf_{m_1}\multipl\cf_{m_2})  
&=& \Ch(\cf_{m_1}\multipl\cf_{m_2})  \td(\ominus\cs_{m_3^{-1}}) \\
&=& \Ch(\ec_{m_3 *}(\e_{m_1}^*\cf_{m_1} 
\e_{m_2}^*\cf_{m_2}\lambda_{-1}(\crr^*)))
\td(\ominus\cs_{m_3^{-1}})  \\
&=& \ec_{m_3 *}(\Ch(\e_{m_1}^*\cf_{m_1} 
\e_{m_2}^*\cf_{m_2}\lambda_{-1}(\crr^*)) \td(T\ec_{m_3}))
\td(\ominus\cs_{m_3^{-1}}))  \\ 
&=& \ec_{m_3 *}(\e_{m_1}^*\Ch(\cf_{m_1}) 
\e_{m_2}^*\Ch(\cf_{m_2}) \Ch(\lambda_{-1}(\crr^*))
\td(T\ec_{m_3})) \td(\ominus\cs_{m_3^{-1}}) \\
&=& \ec_{m_3 *}(\e_{m_1}^*\Ch(\cf_{m_1}) 
\e_{m_2}^*\Ch(\cf_{m_2}) \ctop(\crr) \td^{-1}(\crr) \td(T\ec_{m_3}))
\td(\ominus\cs_{m_3^{-1}}) \\
&=& \ec_{m_3 *}(\e_{m_1}^*\Ch(\cf_{m_1}) 
\e_{m_2}^*\Ch(\cf_{m_2}) \ctop(\crr) \td(\ominus\crr\oplus T\ec_{m_3}))
\td(\ominus\cs_{m_3^{-1}}) \\
&=& \e_{m_3 *}(\e_{m_1}^*\Ch(\cf_{m_1}) \e_{m_2}^*\Ch(\cf_{m_2})
\ctop(\crr) \td(\ominus\crr\oplus T\ec_{m_3}) 
\ec_{m_3}^*\td(\ominus\cs_{m_3^{-1}})) \\
&=& \ec_{m_3 *}(\e_{m_1}^*\Ch(\cf_{m_1}) \e_{m_2}^*\Ch(\cf_{m_2})
\ctop(\crr) \td(\ominus\crr\oplus T\ec_{m_3}\ominus
\ec_{m_3}^*\cs_{m_3^{-1}})), 
\end{eqnarray*}
where the first two equalities follow from the definition of the
multiplication and $\CCh$, the third from the Grothendieck-Riemann-Roch
theorem, the fourth from the fact that the usual Chern character $\Ch$
commutes with pull back and is a homomorphism with respect to the
usual products in the Chow ring, and the fifth from 
Equation~(\ref{eq:UsefulMix}). The sixth and eighth equalities follow
from multiplicativity of $\td$, and the seventh follows from the
projection formula. 

If we let $\ct \in K(X^{\bm})$ be
\begin{equation*}
\ct := \ominus\crr\oplus 
TX^\bm\ominus \res{TX^{m_3^{-1}}}{X^\bm} \ominus
\ec_{m_3}^* \cs_{m_3^{-1}},
\end{equation*}
then by plugging in Equation~(\ref{eq:SSigma}), we obtain
\begin{equation}
\ct = \ominus\crr\oplus 
TX^\bm\ominus\res{TX}{X^\bm}\oplus \res{\cs_{m_3}}{X^\bm}.
\end{equation}
Therefore, we obtain the equality
\begin{equation}\label{eq:HomoLHS}
\CCh(\cf_{m_1}\multipl\cf_{m_2})  = \ec_{m_3
  *}(\e_{m_1}^*\Ch(\cf_{m_1}) \e_{m_2}^*\Ch(\cf_{m_2})
\ctop(\crr) \td(\ct)).
\end{equation}

Similarly, we see that
\begin{eqnarray*}
\CCh(\cf_{m_1})\multipl\CCh(\cf_{m_2}) &=& 
\ec_{m_3 *}( \e_{m_1}^*\CCh(\cf_{m_1})
\e_{m_2}^*\CCh(\cf_{m_2})\ctop(\crr)) \\
&=& \ec_{m_3 *}( \e_{m_1}^*(\Ch(\cf_{m_1}) e_{m_1}^*\td(\ominus\cs_{m_1}))
\e_{m_2}^*(\Ch(\cf_{m_2}) \td(\ominus\cs_{m_2})) \ctop(\crr)) \\
&=& \ec_{m_3 *}( \e_{m_1}^*\Ch(\cf_{m_1})
\e_{m_2}^*\Ch(\cf_{m_2})\ctop(\crr) \td(\ominus e_{m_1}^*\cs_{m_1}\ominus
e_{m_2}^*\cs_{m_2})), 
\end{eqnarray*}
where the first two equalities are by definition and the third is by
multiplicativity of $\td$. Thus, if 
\begin{equation}
\ct' := \ominus\res{\cs_{m_1}}{X^\bm} \ominus \res{\cs_{m_2}}{X^\bm},
\end{equation}
then
\begin{equation}\label{eq:HomoRHS}
\CCh(\cf_{m_1})\multipl\CCh(\cf_{m_2}) = \ec_{m_3 *}( \e_{m_1}^*\Ch(\cf_{m_1})
\e_{m_2}^*\Ch(\cf_{m_2})\ctop(\crr) \td(\ct')).
\end{equation}

$\CCh$ is therefore an algebra homomorphism if and only if the right
hand sides of Equations~(\ref{eq:HomoLHS}) and~(\ref{eq:HomoRHS}) are
equal. A sufficient condition for this equality to hold is that $\ct =
\ct'$, 
but this follows immediately from the definition of $\crr(\bm)$ (See Equation~\ref{eq:R}).

We will now prove that $\CCh$ preserves the trace element. For all $a,b$ in
$G$, for $m = [a,b]$, and for all $\cf_m$ in $\ck_m(X)$, we have
\begin{eqnarray*}
\tauk_{a,b}(\cf_m) 
&=& \chi(\Xab,\res{\cf_m}{\Xab}\otimes\lambda_{-1}(T\Xab\oplus
\res{S_m}{\Xab})^*) \\ 
&=& \int_{\Xab} \Ch(\res{\cf_m}{\Xab}\otimes\lambda_{-1}(T\Xab\oplus
\res{S_m}{\Xab})^*)\cup\td(T\Xab) \\ 
&=& \int_{\Xab} \Ch(\res{\cf_m}{\Xab})\cup\Ch(\lambda_{-1}(T\Xab\oplus\res{S_m}{\Xab})^*)\cup\td(T\Xab) \\ 
&=& \int_{\Xab} \Ch(\res{\cf_m}{\Xab})\cup\ctop(T\Xab\oplus
\res{S_m}{\Xab})\\
 & & \qquad \qquad \qquad \cup\,\td^{-1}(T\Xab\oplus\res{S_m}{\Xab})\cup\td(T\Xab) \\
&=& \int_{\Xab} \Ch(\res{\cf_m}{\Xab})\cup\ctop(T\Xab\oplus
\res{S_m}{\Xab})\cup\td^{-1}(\res{S_m}{\Xab}) \\
&=& \int_{\Xab} \res{\CCh(\cf_m)}{\Xab}\cup\ctop(T\Xab\oplus
\res{S_m}{\Xab})\\
&=& \taua_{a,b}(\CCh(\cf_m)),
\end{eqnarray*}
where we have used the Hirzebruch-Riemann-Roch Theorem in the second
equality, the fact that $\Ch$ preserves the ordinary multiplications in the
third, Equation~(\ref{eq:UsefulMix}) in the fourth, the multiplicativity of
$\td$ in the fifth, and the definition of $\CCh$ (Equation~(\ref{eq:DefCCh})) in
the sixth.
\end{proof}

\begin{rem}
It is instructive to consider the homomorphism property of $\CCh$ when
the obstruction bundle $\crr$ on $X^\bm$ is trivial. When $m_1 = 1$
and $m_2 m_3 = 1$, it is trivial to verify from
Equation~(\ref{eq:m3Equals1}) that
\begin{equation}\label{eq:CChHomo}
\CCh(\cf_{m_1}\multipl\cf_{m_2}) = \CCh(\cf_{m_1})\multipl \CCh(\cf_{m_2}).
\end{equation}
Indeed, Equation~(\ref{eq:CChHomo}) continues to hold even if $\CCh$ were
replaced by the ordinary Chern character $\Ch$. A similar result holds if
$m_2 = 1$ and $m_1 m_3 = 1$.  However, when $m_1 m_2 = 1$
and $m_3 = 1$, then Equation~(\ref{eq:CChHomo}) would fail to hold if $\CCh$ were
replaced by the ordinary Chern character $\Ch$ because of the presence 
of the nontrivial pushforward map $\ec_{m_3 *}$ in Equation
(\ref{eq:m3Equals1}). This shows that the stringy corrections to the
Chern character are necessary even when the obstruction bundle is trivial.
\end{rem}

Finally, $\CCh$ satisfies the usual functorial properties with respect
to equivariant \'etale morphisms.

\begin{thm}\label{thm:RRfunct}
Let $f:X\to Y$ be a $G$-equivariant, \'etale morphism between
smooth, projective varieties $X$ and $Y$ with $G$-action.  The following
properties hold.
\begin{enumerate}
\item \label{thm:Functoriality} (Pullback) The pullback maps
  $$f^*:((\ca(Y,G),\rho),\multipl,\vac)\to ((\ca(X,G),\rho),\multipl,\vac)$$ 
and 
  $$f^*:((\ck(Y,G),\rho),\multipl,\vac)\to ((\ck(X,G),\rho),\multipl,\vac)$$ 
are equivariant morphisms of $G$-graded associative algebras.

\item \label{thm:FunctorialityCCh} (Naturality)
The following diagram commutes.
\begin{equation}
\begin{diagram}
\ck(Y,G)  & \rTo^{f^*} & \ck(X,G) \\
\dTo^{\CCh} & & \dTo^{\CCh}
\\
\ca(Y,G)  & \rTo^{f^*} & \ca(X,G)
\end{diagram}
\end{equation}
\item \label{thm:FunctorialityGRR} (Grothendieck-Riemann-Roch) For
  all $m$ in $G$ and $\cf_m$ in $\ck_m(X)$,
\begin{equation}
f_*(\CCh(\cf_m)\cup\td(TX^m)) = \CCh(f_* \cf_m)\cup\td(TY^m).
\end{equation}
\end{enumerate}
\end{thm}

\begin{proof}
The proof of part~(\ref{thm:Functoriality}) 
follows immediately from the fact that since $f$ is $G$-equivariant and
\'etale, the bundle $f^*TY^m$ is isomorphic to $TX^m$. 

Part~(\ref{thm:FunctorialityCCh}) follows from the naturality of the
ordinary Chern character and the fact that if $f$ is \'etale, then
$f^*\cs_m^Y = \cs_m^X$, where $\cs_m^X$ and $\cs_m^Y$ are as defined
in Equation (\ref{eq:S}) for $X$ and $Y$, respectively.

Part~(\ref{thm:FunctorialityGRR}) follows from these same
considerations, since
\begin{eqnarray*}
f_*(\CCh(\cf_m)\cup\td(TX^m)) 
&=& f_*(\Ch(\cf_m)\cup\td(\ominus\cs^X_m)\cup\td(TX^m)) \\
&=& f_*(\Ch(\cf_m)\cup\td(TX^m)\cup\td(\ominus\cs^X_m)) \\
&=& f_*(\Ch(\cf_m)\cup\td(TX^m)\cup\td(\ominus f^*\cs^Y_m)) \\
&=& f_*(\Ch(\cf_m)\cup\td(TX^m)\cup f^*\td(\ominus\cs^Y_m)) \\
&=& f_*(\Ch(\cf_m)\cup\td(TX^m))\cup \td(\ominus\cs^Y_m) \\
&=& \Ch(f_*\cf_m)\cup\td(TY^m)\cup \td(\ominus\cs^Y_m) \\
&=& \CCh(f_*\cf_m)\cup\td(TY^m),
\end{eqnarray*}
where the projection formula was used in the fifth equality and the ordinary Grothendieck-Riemann-Roch Theorem was used in the sixth. 

\end{proof}

\section{Discrete torsion}\label{sec:DiscTors}

At this point, we wish to make a short comment about discrete
torsion. As discussed in \cite{Ka3}, any $G$-Frobenius algebra 
can be twisted by a discrete torsion, which is 
a $2$-cocycle $\alpha \in Z^2(G,\field^*)$, to obtain a $G$-Frobenius algebra with
twisted sectors of the same dimension. Of course, the same is true
for any pre-$G$-Frobenius algebra
with trace $\tau$, provided the trace $\tau$ is appropriately
twisted, as we explain below.

This procedure allows us to ``twist'' the stringy Chow ring
$\ca(X,G)$ and the stringy $K$-theory $\ck(X,G)$. If one twists both
rings by the same element $\alpha$, then the stringy Chern character
$\CCh$ again provides an allometric isomorphism.

We briefly recall the main points of the construction of twisting by
discrete torsion, omitting the proofs which all follow from rather
straightforward computations.  A reference for the proofs is
\cite{Ka3}.

For $\alpha \in Z^2(G,\field^*)$, let $\field^{\alpha}[G]$ be the
twisted group ring, i.e., $\field^{\alpha}[G]=\bigoplus_{m\in G}
\field e_m$ with the multiplication $e_{m_1}\star
e_{m_2}=\alpha(m_1, m_2)e_{m_1m_2}$.

Set $\epsilon(\gamma,m):=\alpha(\gamma,m)/\alpha( \gamma m
\gamma^{-1} ,\gamma)$ and define $\rho(\gamma)(e_m)=
\epsilon(\gamma,m)e_{\gamma m\gamma^{-1}}.$ Define a bi-linear form
$\eta$ by $\eta(e_{m_+},e_{m_-})=0$ unless $m_+m_-=1$, and
$\eta(e_{m},e_{m^{-1}})=\alpha(m,m^{-1}).$ Lastly, let $\vac=e_1$.

\begin{lm}
 $((\field^{\alpha}[G],\rho),\star,\vac,\eta)$ is a $G$-Frobenius
 algebra.
\end{lm}

\begin{df}
We define the tensor product $\hat \otimes $ of two
pre-$G$-Frobenius algebras $((\ch,\varphi), \star,\vac, \eta,\tau)$
and $((\ch',\varphi'),\star', \vac', \eta',\tau')$ to be the
pre-$G$-Frobenius algebra $\ch\hat \otimes \ch'= \bigoplus_{m\in G}
(\ch\hat \otimes \ch')_m$ with $(\ch\hat \otimes \ch')_m:= \ch_m
\otimes_{\field} \ch'_m$, diagonal multiplication $\star \otimes
\star'$, diagonal $G$-action $\rho \otimes \rho'$, the tensor
product 
pairing
$\eta\otimes \eta'$, unity $\vac\otimes\vac'$,
and trace $\tau \otimes \tau'$.
\end{df}

\begin{prop}
The tensor product of two pre-$G$-Frobenius algebras is a
pre-$G$-Frobenius algebra. Similarly, the tensor product of two
$G$-Frobenius algebras is a $G$-Frobenius algebra.
\end{prop}

\begin{df}
For a pre-$G$-Frobenius algebra $\ch$ and an element $\alpha \in
Z^2(G,\field^*)$, we set
\begin{equation}
\ch^{\alpha}:= \ch\hat\otimes \field^{\alpha}[G].
\end{equation}
\end{df}

Notice that as vector spaces

$\ch^{\alpha}_{m}= \ch_{m} \otimes_{\field} \field \simeq \ch_{m}$.

\begin{lm}
Using the identification $\ch^{\alpha}_m \cong \ch_m$, the
$G$-Frobenius structures for $((\ch^{\alpha},\rho^{\alpha}),
\star^{\alpha}, \vac^{\alpha}, \eta^{\alpha}))$ are
\begin{align*}
v_{m_1} \star^{\alpha} v_{m_2} &:= \alpha(m_1,m_2) v_{m_1} \star v_{m_2},\nonumber \\
\rho^{\alpha}(\gamma) v_m & := \epsilon(\gamma,m) \rho(\gamma) v_m,
\nonumber\\
\eta^{\alpha}(v_{m},v_{m^{-1}}) &:=
\alpha(m,m^{-1})\eta(v_m,v_{m^{-1}})
\nonumber \\
\intertext{and} \tau_{a,b}^{\alpha}( v_{ [a,b] } ) &:=
\frac{\alpha([a,b], b a b^{-1}) \alpha(b,a)}{\alpha(b a b^{-1},b)}
\tau_{a,b}( v_{ [a,b] } ) \nonumber
\end{align*}
for all $v_{m_i}$ in $\ch^\alpha_{m_i}$, $v_m$ in $\ch^\alpha_m$,
and $v_{m^{-1}}$ in $\ch^\alpha_{m^{-1}}$.
\end{lm}

\begin{prop}
The pre-$G$-Frobenius algebras $\ch$ and $\ch^{\alpha}$ are
isomorphic if and only if $\alpha$ is a coboundary; that is,
$[\alpha]=0\in H^2(G,\field^*).$
\end{prop}

\begin{prop}
If $\Phi: \ch \rightarrow \ch'$ is an isomorphism (or allometric
isomorphism) of pre-$G$-Frobenius algebras, then $\Phi\otimes id$ is
an isomorphism (respectively allometric isomorphism) between
$\ch^{\alpha}$ and $\ch^{\prime\alpha}.$
\end{prop}

\begin{crl}\label{crl:twistCh}
Let $\CCh:\ck(X,G)\rightarrow \ca(X,G)$ denote the stringy Chern
character.  For all $\alpha\in Z^2(G,\field^*)$, the map
$\CCh^{\alpha}= \CCh\otimes id:\ck^{\alpha} (X)\rightarrow
\ca^{\alpha}(X)$ is an allometric  isomorphism.
\end{crl}

\section{Relation to Fantechi-G\"ottsche, Chen-Ruan, and Abramovich-Graber-Vistoli}
\label{sec:obs}\label{sec:FG=us}

In \cite{FG} Fantechi and G\"ottsche  describe a ring $\ch^{\bullet}(X,G)$,
which we call the \emph{stringy cohomology}, associated to every manifold $X$
with an action by a finite group $G$.  The stringy cohomology is also a
$G$-Frobenius algebra \cite{JKK}, and the ring of $G$-invariants of
$\ch^{\bullet}(X,G)$ is known to be isomorphic to the Chen-Ruan orbifold
cohomology $H^{\bullet}_{orb}([X/G],\nq)$ of the quotient stack $[X/G]$.
Abramovich, Graber, and Vistoli \cite{AGV} have an algebraic construction,
similar to that of Chen and Ruan, of what we would call an orbifold Chow
ring.   

\begin{rem}
Note that Abramovich, Graber, and Vistoli call their ring the ``stringy 
Chow ring,'' but we prefer to reserve the word \emph{stringy} for $G$-Frobenius
structures associated to a manifold with a specific group action, and use the
word \emph{orbifold} for Frobenius algebras that are associated to orbifolds,
and are thus presentation independent.  
The general philosophy is that a {stringy} construction associated to a finite group $G$ acting on $X$ should have, as its ring of invariants, an {orbifold} construction for the quotient stack $[X/G]$.  The orbifold construction should also generalize to stacks which are not global quotients by finite groups.
\end{rem}

Just as our constructions rely on the special bundle $\crr(\bm)$, the
constructions of Fantechi-G\"ottsche, Chen-Ruan, and
Abramovich-Graber-Vistoli all use an obstruction bundle
arising in the theory of stable maps---either stable maps into an orbifold or
$G$-stable maps (see \cite{JKK}) into a manifold with $G$-action.  The
description of these obstruction bundles is rather technical and is generally
difficult to use for computation. 

In this section, we prove that the obstruction bundle of Fantechi-G\"ottsche is equivalent to our bundle $\crr(\bm)$.  It is known that their construction 
agrees (after taking $G$-invariants) with that of Chen-Ruan and
Abramovich-Graber-Vistoli \cite[\S2]{FG}.  In Section~\ref{sec:orb}, and
especially in Theorem~\ref{thm:orbObs}, we will generalize this to general
orbifolds---not just those which are global quotients by finite groups.   

For our purposes the most important consequence of the equivalence of our construction with that of Fantechi-G\"ottsche is that the element $\crr(\bm)\in K(X^\bm)$
 is actually represented by a vector bundle.  But another important consequence is that their obstruction bundle may now be described solely in terms of the $G$-action on the tangent  bundle of $X$, restricted to various fixed-point loci.  This greatly simplifies the computation of stringy cohomology, orbifold cohomology, and orbifold Chow. In particular, it allows us to circumvent all of the technical details of those constructions, including stable curves, stable maps, admissible covers, and moduli spaces.

\subsection{The obstruction bundle of Fantechi and
G\"ottsche}\label{sec:FGObs}

We briefly review the construction of the obstruction bundle of Fantechi and
G\"ottsche \cite{FG}.  For each triple $\bm=(m_1,m_2,m_3) \in G^3$ such that
$m_1 m_2 m_3=1$,  
let $\langle \bm \rangle$ be the subgroup generated by the elements
$m_1,m_2$, and $m_3$.  There is a presentation of the fundamental
group $\pi_1(\mathbb{P}^1-\{0,1,\infty\})$ as $\langle c_1,c_2,c_3|
c_1 c_2 c_3 =1 \rangle$, where $c_1,c_2$ and $c_3$ are 
based loops around $p_1=0, p_2=1,$ and $ p_3=\infty$, respectively.  We define a natural homomorphism $\pi_1(\mathbb{P}^1-\{0,1,\infty\}) \to \langle
\bm\rangle$, taking $c_i$ to $m_i$.  This defines a principal $\langle
\bm \rangle$-bundle over $\mathbb{P}^1 -\{0,1,\infty\}$ which extends
to a smooth connected curve $E$. The curve $E$ has an action of
$\langle \bm \rangle$ such that the quotient $E/\langle\bm\rangle$ has
genus zero, and the natural map $E \to E/\langle\bm\rangle$ is
branched at the three points $p_1,p_2, p_3$ with monodromy
$m_1,m_2,m_3,$ respectively.

Let $\pi:E \times X^{\bm} \to X^{\bm}$ be the second projection.  The obstruction bundle 
of Fantechi and G\"ottsche, which we denote by $\crrfg(\bm)$, on $X^{\bm}$ is
\begin{equation}\label{eq:obstr-bundle}
  \crrfg(\bm):=R^1\pi^{\langle\bm\rangle}_*(\co_E \boxtimes
  TX|_{X^{\bm}}).
\end{equation}
One can check that the restriction of the bundle
$\crrfg(\bm) \in K(X^\bm)$ to a connected component $U$ of $X^\bm$ has
rank \begin{equation}\label{eq:RankR}
a(m_1,U)+a(m_2,U)+a(m_3,U)-\codim(U\subseteq X).
\end{equation}

\begin{rem}
For those familiar with quantum cohomology, this obstruction bundle is
the analogue of the obstruction bundle for stable maps, but with
additional accounting for the structure of the group action on $X$.
That is, $c_{top}(\crrfg)$ is the virtual fundamental class on
(distinguished components of) the moduli space of genus-zero,
three-pointed $G$-stable maps into $X$.  The base space $X^{\bm}$ in
the definition of the obstruction bundle is the distinguished
component $\xi_{0,3}(X,0,\bm)\cong{pt}\times X^{\bm}$ of
$\MM_{0,3}(X,0,\bm)$.  The interested reader may refer to
\cite[\S6]{JKK} for more details.
\end{rem}

\begin{thm}\label{thm:magic}
  Let $X$ be a smooth variety (not necessarily projective, or even
  proper) with an action of a finite group $G$.
  If $\bm = (m_1, m_2, m_3) \in G^3$ is such that
$m_1m_2m_3 = 1$, 
then on the fixed
point locus $X^\bm := X^{m_1}\cap X^{m_2}$, we have 
\begin{align}\label{eq:Magic}
\crrfg(\bm) &= \crr(\bm)\\
\notag				&= TX^\bm \ominus \res{TX}{X^\bm}\oplus\bigoplus_{i=1}^3
  \res{\cs_{m_i}}{X^{\bm}},
\end{align} in the K-theory $K(X^\bm)$ of $X^\bm$.
\end{thm}

\begin{crl}\label{cor:RVB}
For each triple $\bm = (m_1,m_2,m_3)$ with $m_1m_2m_3 = 1$, the element $\crr(\bm) \in K(X^\bm)$ is represented by a vector bundle on $X^\bm$.
\end{crl}

As a first step to proving Theorem~\ref{thm:magic}, we prove   Lemma~\ref{thm:rep-magic}.  
The basic setup for Lemma~\ref{thm:rep-magic} is as follows.  Let $E$ be a
smooth algebraic curve of genus $\tilde{g}$, not necessarily
connected, with a finite group $G$ acting effectively on $E$.  Assume
that the quotient $E/G$ has genus $g$.  Denote the orbits where the
action is not free by $p_1, \dots, p_n \in E/G$.  A choice of base
point $\pt \in E$ induces a homomorphism of groups
$$\varphi_{\pt}:\pi_1(E/G - \{p_1, \dots, p_n\}, p) \to G,$$ where $p$ is
the image of $\pt$ in $E/G$ (we assume $p \notin \{p_1, \dots, p_n\})$.
Denote by $H$ the image of $\varphi_{\pt}$ in $G$.  Note that the number
$\alpha$ of connected components of $E$ is the index $[G:H]$.  There is a
presentation of $\pi_1(E/G-\{p_1, \dots, p_n\},p)$ of the form $\langle a_1,
\dots, a_g, b_1, \dots, b_g, c_1, \dots, c_n |\prod^n_{i=1}
c_i=\prod^g_{j=1}[a_j,b_j]\rangle$, where the $c_i$ are loops around the
points $p_i$.  For each $i \in \{1, \dots, n\}$ we call the image
$m_i:=\varphi_{\pt}(c_i) \in G$ of $c_i$ the \emph{monodromy} around $p_i$,
and we denote the order of $m_i$ in $G$ by $r_i$.  Of course, a different choice of
$\pt$ will change all of the $m_i$ by simultaneous conjugation with an
element of $G$.    

The following lemma describes the $G$-module structure of the cohomology
$H^1(E;\co_E)$.  It has recently come to our attention that a 
related result can be found in \cite{Kani}. 
\begin{lm}\label{thm:rep-magic}
Given the setup described above, and letting $\nc[G]$ denote the group
ring regarded as a $G$-module under 
multiplication, we have
the following equality in the representation ring of $G$,
\begin{equation}\label{eq:magic2}
H^1(E;\co_E) = 
\nc[ G/H] 
\oplus (g-1) \nc[G]\oplus 
\bigoplus_{i=1}^n\bigoplus_{k_i=0}^{r_i-1}
\frac{k_i}{r_i} \ind^G_{\langle m_i\rangle }\bv_{m_i,k_i},
\end{equation}
where 
$\bv_{m_i,k_i}$ is the irreducible representation of $\langle m_i \rangle$
such that $m_i$ acts by multiplication by 
$\exp(-2\pi i k_j/r_j)$,
and 
$\ind^G_{\langle m_i\rangle }\bv_{m_i,k_i}$ 
is the induced representation
$\nc[G]\otimes_{\nc[\langle m_i \rangle ]} \bv_{m_i,k_i}.$
\end{lm}

\begin{proof}
  It suffices to check that these two virtual representations have the same
  virtual character.  The trace of the action of an element $\gamma \in G$ on
  the right hand side is $$
  \chi_{\nc[G/H]} (\gamma) + 
(g-1) |G| \delta_{\gamma,e} + \sum^n_{i=1} \sum^{r_i-1}_{k_i=0}
\frac{k_i}{r_i} \chi_{\ind^G_{\langle m_i\rangle } \bv_{m_i,k_i}} (\gamma).
$$
  
  It is well known (e.g., \cite[ex 3.19]{FH}) that, for a
  representation $V$ of a subgroup $H \leq G$, we have $$\chi_{\ind_H^G{V}}
  (\gamma) =\frac{|G|}{|H|} \sum_{{\sigma \in H \cap [\![\gamma]\!] }}
  \frac{\chi_V(\sigma)}{|[\![\gamma]\!]|},$$
  where $[\![\gamma]\!]$ is the
  conjugacy class of $\gamma$ in $G$.  In our case, with $H=\langle
  m_i\rangle$ of order $r_i$, and 
$V=\bv_{m_i,k_i}$,
we have
\begin{align*}
  \chi_{\ind^G_HV}(\gamma)
  & = \frac{|G|}{r_i|[\![\gamma]\!]|} \sum_{m^l_i \in [\![\gamma]\!] } \chi_V(m^l_i)\\
  &= \frac{|G|}{r_i |[\![\gamma]\!]|} \sum_{m^l_i \in [\![\gamma]\!] } \zeta^{lk_i}_i
\end{align*}
where 
$\zeta_j=\exp(-2 \pi i/r_j)$ 
for each $j \in \{1, \dots, n\}$.
Thus the trace of the right hand side of equation~(\ref{eq:magic2})
becomes 
\begin{equation}\tr_{RHS}(\gamma) = 
\chi_{\nc[G/H]}(\gamma) 
+ (g-1) |G|
\delta_{\gamma,e} + \sum^n_{i=1} \sum^{r_i-1}_{k_i=0}
\frac{k_i|G|}{r^2_i{|[\![\gamma]\!]|}} \sum_{m^l_i \in [\![\gamma]\!] }
\zeta^{lk_i}_i.
\end{equation}
 
If $\gamma=e$ is the identity element of $G$, we have
 \begin{align}
   \tr_{RHS} (e)
   &=\alpha + |G|(g-1) + \sum^n_{i=1} \sum^{r_i-1}_{k_i=0} \frac{k_i|G|}{r^2_i}\nonumber\\
   &= \alpha + |G|(g-1) +|G| \sum^n_{i=1} \frac{r_i-1}{2r_i}\nonumber\\
   &= \dim_{\nc}H^1(E,\co_E),
\end{align}
where the last equality follows from the Riemann-Hurwitz formula and the fact that the genus of $E/G$ is $g$.

If $\gamma \neq e$, then
\begin{align}
  \tr_{RHS}(\gamma)
  &= 
\chi_{\nc[G/H]}(\gamma)
 + \sum_{i=1}^n \frac{|G|}{r^2_i|[\![\gamma]\!]|} \sum_{m^l_i \in [\![\gamma]\!] } \sum^{r_i-1}_{k_i=0} k_i \zeta^{lk_i}_i \nonumber\\
  &=
\chi_{\nc[G/H]}(\gamma) 
+ \sum^n_{i=1} \frac{|G|}{r^2_i|[\![\gamma]\!]|} \sum_{m^l_i \in [\![\gamma]\!]}
  r_i\frac{\zeta^{-l}_i}{1-\zeta^{-l}_i}\nonumber\\ 
  &= 
\sum_{\substack{\sigma \in G/H\\ \gamma \sigma=\sigma}}
  1 + \sum^n_{i=1} \frac{|G|}{r_i|[\![\gamma]\!]|} \sum_{m^l_i \in
    [\![\gamma]\!] } \frac{\zeta^{-l}_i}{1-\zeta^{-l}_i},
\end{align}
where the last equality follows from standard results on induced
representations \cite[3.18]{FH}.

This formula is related to fixed points of the action of $\gamma$ on
$E$ as follows.  The element $\gamma$ can only fix points that lie over
the $p_i$, for $i \in \{1, \dots, n\}.$ If $\pt_i$ is a point over
$p_i$ fixed by $\gamma$, then $\pt_i$ has monodromy
conjugate to $m_i$,
and thus $\gamma$ must be conjugate to $m^l_i$ for some $l$.
Conversely, if $\gamma$ is conjugate to $m^l_i$ for some $l$, then
$\gamma$ fixes all points $\pt_i$ that lie over $p_i$, and $\gamma$
acts on the tangent space $T_{\pt_i}E$ by 
$\zeta^{-l}_i$.

If both $\pt_i$ and $\pt_i'$ are fixed by $\gamma$ with action
$\zeta^{-l}_i$ 
on the tangent space, then $\pt_i'=\pt_i\sigma$ for some
$\sigma \in G$, such that $\sigma$ commutes with $\gamma$, but if
$\sigma \in \langle m_i\rangle $, then $\pt_i=\pt_i'$.  So the number
of such points lying over $p_i$ is exactly
$\frac{|Z_G(\gamma)|}{|\langle m_i\rangle|}
=\frac{|G|}{r_i|[\![\gamma]\!]|}$, where $Z_G(\gamma)$ is the
centralizer of $\gamma$ in $G$.

The term 
$\sum_{\substack{\sigma \in G/H\\ \gamma\sigma=\sigma}} 1$ 
counts connected components of $E$ which are
mapped to themselves by $\gamma$; that is, it is the trace of $\gamma$
for the natural representation of $G$ on $H^{0}(E,\co_E)$.  If we now
denote by 
$d \gamma_{\pt}=\zeta^{-l}_i$ 
the action of $\gamma$ on the
tangent space $T_{\pt} E$ at a fixed point $\pt \in E$, the above
argument shows that
$$\tr_{RHS} (\gamma)
= \chi_{H^{0}(E,\co_E)}(\gamma) + \sum_{\pt \text{ fixed by } \gamma} 
\frac{d \gamma_{\pt}}{1-d\gamma_{\pt}}.
$$

But the Eichler trace formula says that this is precisely the trace of
the action of $\gamma$ on $H^1(E,\co_E);$ that is, the traces of
equation ~(\ref{eq:magic2}) agree (see \cite[\S V.2.0]{FK} for $E$
connected with $\tilde{g}>1$, and \cite[\S17]{Sh} in general).
\end{proof} 

\begin{proof}[Proof of Theorem~\ref{thm:magic}] 

For any $\bm=(m_1,m_2,m_3)$ with $m_1 m_2 m_3=1$, the curve $E$ in the
definition of $\crr(\bm)$ is connected and has an effective action of
$G':=\langle m_1,m_2,m_3\rangle $ with quotient $\mathbb{P}^1=E/G'$
and three branch points $p_1,p_2,p_3$.

We have
$$\crrfg (\bm) =R^1 \pi^{G'}_*(\co_E \boxtimes TX|_{X^{\bm}}) \cong (H^1(E, \co_E)\otimes 
TX|_{X^{\bm}})^{G'},$$
and by Lemma~\ref{thm:rep-magic} this is 
\begin{align}
  (H^1(E, \co_E)\otimes TX|_{X^{\bm}})^{G'} &= \left(\left(\nc \ominus
      \nc[G'] \oplus \bigoplus_{i=1}^3\bigoplus_{k_i = 0}^{r_i-1}
      \frac{k_i}{r_i} \ind_{\langle m_i\rangle }^{G'}\bv_{m_i,k_i}
      \right)\tensor TX|_{X^{\bm}}\right)^{G'}\nonumber\\ 
  &=TX^{\bm} \ominus TX|_{X^{\bm}} \oplus
  \bigoplus_{i=1}^3\bigoplus_{k_i = 0}^{r_i-1}
  \frac{k_i}{r_i} \left( \ind_{\langle m_i\rangle }^{G'}\bv_{m_i,k_i}
  \tensor TX|_{X^{\bm}}\right)^{G'}\nonumber\\ 
  &=TX^{\bm} \ominus TX|_{X^{\bm}} \oplus
  \bigoplus_{i=1}^3\bigoplus_{k_i = 0}^{r_i-1}
  \frac{k_i}{r_i}  \left(TX|_{X^{\bm}}\right)_{m_i,k_i}\nonumber\\
  &= TX^{\bm} \ominus TX|_{X^{\bm}} \oplus
  \bigoplus_{i=1}^3\res{\cs_{m_i}}{X^{\bm}}\nonumber \\
  &=\crr(\bm).
\end{align}
\end{proof}

\subsection{The Abelian case}
It is instructive to consider the special case where $G$ is an Abelian
group. In this case, our analysis of the obstruction bundle $\crrfg$
yields, as a simple corollary, a result originally due to Chen and Hu
\cite{ChHu}.

Consider the obstruction bundle $\crrfg$ over $X^\bm$, where $\bm =
(m_1,m_2,m_3)$ in $G^3$ satisfies $m_1 m_2 m_3 = 1$. Let us assume
without loss of generality that $G = \langle \bm \rangle$. Since $G$
is Abelian, one can simultaneously diagonalize the actions of $m_i$,
for $i=1,2,3$ on $\crrfg$. If $W_\bm$ denotes the normal bundle of
$X^\bm$ in $X$, then we have the simultaneous eigenbundle
decomposition
\begin{equation}\label{eq:W}
W_\bm = \bigoplus_\bk W_{\bm,\bk},
\end{equation}
where the sum is over all $\bk = (k_1,k_2,k_3)$ such that $k_i =
0,\ldots,r_i-1$, and $r_i$ is the order of $m_i$ for all $i \in \{1,2,3\}$.
The eigenbundle $W_{\bm,\bk}$ of $W_\bm$ is the bundle where for all
$j\in \{1,2,3\}$ each $m_j$ has an eigenvalue 
$\exp(2\pi i k_j/r_j)$.  
The
following proposition is an easy corollary of Theorem~\ref{thm:magic}.
\begin{prop} \cite{ChHu}
When $G$ is Abelian, under the above assumptions, then
\begin{equation}\label{eq:AbelianR}
       \crrfg = \bigoplus_\bk W_{\bm,\bk}
\end{equation}
in  $K(X^\bm)$, where the sum is over triples $\bk = (k_1,k_2,k_3)$, for
$k_i=0,\ldots,r_i-1$ and $i=1,2,3$, such that 
\begin{equation}\label{eq:AbelianDegree}
\frac{k_1}{r_1} + \frac{k_2}{r_2} + \frac{k_3}{r_3} = 2.
\end{equation}
\end{prop}
\begin{proof}
It is a straightforward exercise to verify that the right hand side of
Equation~(\ref{eq:Magic}) agrees with Equation~(\ref{eq:AbelianR})
when $G$ is Abelian.
\end{proof}

\begin{rem}
Chen and Hu use this characterization of the obstruction bundle to
give a de~Rham model for Chen-Ruan orbifold cohomology when the
orbifold arises as the quotient of a variety by an Abelian group. It
would be interesting to see how their constructions can be generalized
to non-Abelian groups in light of Equation~(\ref{eq:Magic}).
\end{rem}

\section{The orbifold K-theory of a stack}\label{sec:stack}\label{sec:orb}

In this section, we introduce two variants of orbifold K-theory. The first
is given by extending the main constructions of stringy K-theory to an
orbifold $\cx$.  We use a special vector bundle $\crrt$ on the double inertia
stack $\cccx$, an orbifold analogue to $\crrt$ from stringy K-theory, to
define a new product $\multipl$ on the K-theory $\okf(\cx) := K(\ccx)$ of the
inertia stack $\ccx$.  We call the resulting algebra the 
\emph{full orbifold K-theory} of $\cx$. The second construction is associated
to a global quotient $\cx = [X/G]$ of a smooth, projective variety $X$  by a 
finite group $G$. The ring $\ok(\cx)$, which we call \emph{small orbifold
K-theory}, is the algebra of $G$-invariants of its stringy K-theory. We show
that, after tensoring with $\nc$, the orbifold K-theory $\ok(\cx)$ is a
sub-algebra of the pre-Frobenius algebra $\okf(\cx)$.

We also make the analogous constructions for Chow rings, but 
unlike in the case of K-theory, the (full) orbifold Chow ring
is isomorphic, as a pre-Frobenius algebra, to the invariants of the stringy
Chow ring.

We also define a ring homomorphism,
the \emph{full orbifold Chern character,}
$\CChf:\okf(\cx) \to \oa(\cx)$. The construction of $\okf(\cx)$ is to Givental
and Lee's quantum K-theory \cite{YPL}, as Chen-Ruan \cite{CR2} and
Abramovich-Graber-Vistoli's \cite{AGV} orbifold  quantum cohomology is to
quantum cohomology.  Furthermore, as a vector space, our construction  of
$\ok(\cx)$ agrees with the construction of Adem and Ruan 
\cite{AR1}, but unlike theirs, our orbifold  product has the virtue that the associated
orbifold Chern character $\och$ is a \emph{ring} isomorphism---not just an
additive isomorphism.

The main results of this section are Theorems~\ref{thm:ccx-preFA} and
\ref{thm:pistuff} which say, 
among other things,
that $\okf(\cx)$ with the new product is a 
pre-Frobenius algebra, that $\ok(\cx)$ is a
sub-pre-Frobenius algebra of $\okf(\cx)$, and that the 
orbifold Chern character $\CChf:\okf(\cx)\to \oa(\cx)$ is a
homomorphism of pre-Frobenius algebras  
which induces the isomorphism $\och:\ok(\cx)\rTo^{\sim} \oa(\cx)$.

\subsection{The Full Orbifold K-theory and the Orbifold Chow Ring}

First, we need to establish some notation.  We denote by $Z_G(g)$ the centralizer of an
element $g$ in a group $G$, and we denote by $Z_G(g,g')$ the intersection
$Z_G(g,g'):=Z_G(g) \cap Z_G(g')$.  For any group $G$, we denote the set of all conjugacy classes in $G$ by ${\Gb}$.

Recall that the inertia stack $\ccx$ of $\cx$ is defined to be 
the stack whose $T$-valued points are pairs $(x,g)$, where $x$ is a $T$-valued point of $\cx$ and $g$ is an automorphism of $x$ in $\cx(T)$.  An equivalent definition is 
$$\ccx:=\cx \times_{(\cx\times\cx)}\cx$$ corresponding to the diagram 
$$
\begin{diagram}
\ccx 	&\rTo 	& \cx\\
\dTo	& 		& \dTo^{\Delta}\\
\cx		&\rTo^{\Delta}	& \cx \times \cx.
\end{diagram}
$$
We can write 
$$\ccx = \coprod_{[\![g]\!]} \cx_{[\![g]\!]},$$
where the indices run over conjugacy classes ${[\![g]\!]}$ of local automorphisms and $\cx_{[\![g]\!]}$ is the locus of pairs $(x,g')$ with $g' \in {[\![g]\!]}$. 
There is an obvious inclusion $$j:\cx \rInto \ccx,$$ taking $x$ to $(x,\vac_x)$.  There is also a canonical involution $\sigma: \ccx \to\ccx$ given by $\sigma:
(x,g) \mapsto (x,g^{-1})$. 

The double inertia stack $\cccx$ is defined to be the stack whose $T$-valued points are triples $(x,g,h)$, where $x$ is a $T$-valued point of $\cx$ and $g$ and $h$ are both automorphisms of $x$ in $\cx(T)$.  Again we can decompose $\cccx$ as
 $$\cccx:=\coprod_{[\![g_1,g_2]\!]} \cx_{[\![g_1,g_2]\!]},$$ where the indices run over (diagonal) conjugacy classes of pairs of local automorphisms $[\![g_1, g_2]\!]$, and
$\cx_{[\![g_1,g_2]\!]}$ is the locus of triples $(x,g'_1,g'_2)$ with $(g'_1, g'_2) \in [\![g_1,g_2]\!]$.
 
There are three ``evaluation'' morphisms $$\cccx \rTo^{ev_i} \ccx$$ defined
by  
\begin{align*}
ev_1&:  (x,g_1,g_2) \mapsto (x,g_1)\\
ev_2&:  (x,g_1,g_2) \mapsto (x,g_2)\\
ev_3&:  (x,g_1,g_2) \mapsto (x,(g_1g_2)^{-1}).
\end{align*}
More generally, given any elements $a,b,m\in G$, such that $m$ is in the subgroup $\langle a,b\rangle$ generated by $a$ and $b$, there is a morphism $$\varepsilon_{a,b,m}:\cx_{[\![a,b]\!]} \rTo \cx_{[\![m]\!]}.$$
We define $$\check{ev}_i = \sigma \circ ev_i,$$ where $\sigma: \ccx \to\ccx$ is the canonical involution.

Also, forgetting automorphisms entirely gives  morphisms $$i:\ccx \rTo \cx$$ and $$J:\cccx \rTo \cx.$$  We have $$J = i\circ ev_1 = i \circ ev_2 = i\circ ev_3,$$ and $$i\circ j = \vac_{\cx}.$$

For most of our constructions we will need to impose an additional condition on the
stacks.
\begin{df}
We say that a stack $\cx$ \emph{satisfies the KG-condition} if the Grothendieck group $K^{\mathrm{naive}}(\cx,\nz)$ of (orbi-)vector bundles is isomorphic to the Grothendieck group $K^0(\cx,\nz)$ of perfect complexes and to the Grothendieck group $G_0(\cx,\nz)$ of coherent sheaves on $\cx$. 

If $\cx$ satisfies the KG-condition, we will simply write $K(\cx)$ to denote this group with rational coefficients: $$K(\cx):=K^0(\cx,\nz)\otimes \nq \cong K^{\mathrm{naive}}(\cx,\nz)\otimes \nq \cong  G_0(\cx,\nz)\otimes \nq.$$

If the stack $\cx$, its inertia stack $\ccx$, and its double inertia stack $\cccx$ all satisfy the KG-condition, then we say that \emph{$\cx$ satisfies condition \Kcond}.
\end{df}

If a stack $\cx$ is smooth, then it is always true that $K(\cx) \cong G(\cx)$ \cite[\S2]{JoshHI}.  Moreover, if $\cx$ has the \emph{resolution property} that every coherent sheaf is a quotient by a vector bundle, then $K^{naive}(\cx) \cong K(\cx)$ \cite[Prop 2.3]{JoshK}.  Smooth Deligne-Mumford stacks with a finite stabilizer group which is generically trivial and with a coarse moduli space which is a separated scheme satisfy the resolution property \cite[Thm 1.2] {Tot}, thus they satisfy the KG-property.   In particular, condition \Kcond\ holds for the quotient $\cx = [X/G]$ of a smooth projective variety $X$ by a finite group $G$.

As in the stringy case, for each conjugacy
class $[\![g]\!]$, with $g$ of order $r$, the element $g$ acts by $r$-th
roots of unity on $W_{[\![g]\!]}:=\res{T\cx}{\cx_{[\![g]\!]}}$.  We define
$W_{[\![g]\!],k}$ to be the eigenbundle of $W_{[\![g]\!]}$, where $g$ acts by
multiplication by 
$\zeta^k = \exp(2\pi i k/r)$.
Note that this eigenbundle is determined only by the conjugacy class
$[\![g]\!]$ rather than by the particular representative $g$. Finally,
we define \begin{equation}\cso_{[\![g]\!]}: = \bigoplus_{k=0}^{r-1}\frac{k}{r} W_{[\![g]\!],k}
\in K(\cx_{[\![g]\!]}).
\end{equation} This allows us to define 
$\cso \in K(\ccx)$ as 
\begin{equation}\label{eq:orbS}\cso =
\bigoplus_{[\![g]\!]}\cso_{[\![g]\!]}.\end{equation} 

\begin{df}
The \emph{orbifold obstruction bundle} $\ob\in K(\cccx)$ 
is
\begin{equation}\ob:= T\cccx \ominus J^*T\cx \oplus \bigoplus^3_{i=1}
ev^*_i \cso.
\end{equation}
\end{df}

We  can now define $\okf(\cx)$ as follows.
\begin{df}
As a vector space, the \emph{full orbifold K-theory} $\okf(\cx)$ is $K(\ccx)$. 
The orbifold product of
  two bundles $\cf$ and $\cf'$ in $\okf(\cx)$
is defined to be 
  \begin{equation}\label{eq:multi}
  \cf \multipl
  \cf' :=(\check{ev}_3)_* \left( ev_1^*(\cf) \otimes ev_2^* (\cf') \otimes
  \lambda_{-1}(\ob^*)\right).
  \end{equation}
  
  The trace element $\taukof\in 
\okf(\cx)^*
$ is defined by
  \begin{equation}
	\taukof (\cf_{[\![m]\!]}):=
	\sum_{\substack{[\![a,b]\!]\\ [a,b]
	\in [\![m]\!]
}}
	\chi\left(\cx_{[\![a,b]\!]}, \varepsilon_{a,b,m}^*(\cf_{[\![m]\!]})\otimes \lambda_{-1}\left( \left(T\cx_{[\![a,b]\!]} \oplus \varepsilon_{a,b,m}^*(\cs_{[\![m]\!]})\right)^*\right)\right),
  \end{equation}
  where the sum runs over all (diagonal) conjugacy classes of pairs $[\![a,b]\!]$ 
  of local automorphisms 
  such that $[a,b] \in [\![m]\!]$. 
  
Finally, the symmetric bilinear form
$\etakof$ on 
$\okf(\cx)$ 
is 
$$\etakof(\cf, \cg) := \chi\left( \ccx, \cf \otimes \sigma^*(\cg)\right).$$
\end{df}

We make a similar definition for the Chow ring.
\begin{df}
As a vector space, we let $\oa(\cx) := A^{\bullet}(\ccx)$. 
The orbifold product of  $\cf$ and $\cf'$ in 
$\oa(\cx)$ is defined to be 
  \begin{equation}
  \cf \multipl
  \cf' :=(\check{ev}_3)_* \left( ev_1^*(\cf) \otimes ev_2^* (\cf') \otimes
c_{top}(\ob)
\right).
  \end{equation}
The trace element, denoted by $\tauao$, where $\tauao \in \oa(\cx)^*$ is
given by 
\begin{equation}
	\tauao (\cf_{[\![m]\!]}):=
	\sum_{\substack{[\![a,b]\!]\\ [a,b]
	\in [\![m]\!]
}}
	\int_{\cx_{[\![a,b]\!]}}
	\varepsilon_{a,b,m}^*(\cf_{[\![m]\!]})\cup \ctop\left(
	T\cx_{[\![a,b]\!]} \oplus
	\varepsilon_{a,b,m}^*(\cs_{[\![m]\!]})\right),
  \end{equation}
and the pairing, denoted by $\etaao$, on $\oa(\cx)$ 
as 
$$\etaao(\cf, \cg) := \int_\ccx \cf \cup \sigma^*(\cg).$$
\end{df}

A related construction has now appeared in \cite{AR2} in the topological
category although they not describe a Frobenius structure or a Chern
character.  It would be interesting to explore the connections between their
construction and ours.

\begin{thm}\label{thm:ccx-preFA}

Let $\cx$ be a smooth Deligne-Mumford stack.
\begin{enumerate}
\item\label{fullorbChow} $(\oa(\cx),\multipl,\etaao,\tauao)$ is a pre-Frobenius
algebra called \emph{the orbifold Chow ring of $\cx$}. 
Moreover, $(\oa(\cx),\multipl)$ is isomorphic to the 
``orbifold Chow ring'' of
\cite{AGV}.
\item\label{fullorbK} If $\cx$ satisfies condition \Kcond, then
      $(\okf(\cx),\multipl,\etakof,\taukof)$ is a pre-Frobenius algebra called
      \emph{the full orbifold K-theory of $\cx$}. 
\item\label{strchern} 
If $\cx$ satisfies condition \Kcond, then 
the \emph{full orbifold Chern character} 
$\CChf: \okf(\cx) 
\rTo 
\oa(\cx)
$ is an allometric homomorphism of pre-Frobenius algebras, 
where $\CChf := \cch \cup \td(\ominus\cso)$, and $\cso$ is given in Equation~(\ref{eq:orbS}). 
\end{enumerate}
\end{thm}

\begin{proof}
Parts~(\ref{fullorbChow}) and  (\ref{fullorbK}) will be proved in subsection~\ref{subsec:proofs}.

Part~ (\ref{strchern}) follows from the definition of $\ob$, of $\cs$, and of
$\CChf$, in a manner similar to the proof of Theorem~\ref{thm:ChernHomo}. 
\end{proof}

We conclude this section with a full orbifold version of
Theorem~\ref{thm:RRfunct}.

\begin{thm}\label{thm:orbRRfunct}
Let $f:\cx \to \cy$ be an \'etale morphism 
of smooth, Deligne-Mumford
stacks such that $\cx$ and $\cy$ both satisfy condition \Kcond.
The following properties hold.
\begin{enumerate}
\item \label{thm:orbFunctoriality} (Pullback) The pullback maps
  $f^*:(\oa(\cy),\multipl,\vac_\cy) \to (\oa(\cx),\multipl,\vac_\cx) $ 
  and   $f^*:(\okf(\cy),\multipl,\vac_\cy) \to (\okf(\cx),\multipl,\vac_\cx) $ are 
ring homomorphisms.
\item \label{thm:orbFunctorialityCCh} (Naturality)
The following diagram commutes.
\begin{equation}
\begin{diagram}
\okf(\cy)  & \rTo^{f^*} & \okf(\cx) \\
\dTo^{\CChf} & & \dTo^{\CChf}
\\
\oa(\cy)  & \rTo^{f^*} & \oa(\cx)
\end{diagram}
\end{equation}
\item \label{thm:orbFunctorialityGRR} (Grothendieck-Riemann-Roch) For
  all $\cf$ in $\okf(\cx)$,
\begin{equation}
f_*(\CChf(\cf)\cup\td(T\cx)) = \CChf(f_* \cf)\cup\td(T\cy). 
\end{equation}
\end{enumerate}
\end{thm}

The proof of this theorem is a straightforward adaptation of the proof of its
stringy counterpart, Theorem~\ref{thm:RRfunct}.

\subsection{Small Orbifold K-theory}
We now introduce the algebra called the \emph{small orbifold K-theory}
$\ok(\cx)$ when $\cx  = [X/G]$ is the global quotient of a smooth, projective
variety $X$ by a finite group $G$. We will then explain its relationship to
to $\okf(\cx)$.

For the rest of this section, assume that $\cx$ is a global quotient $[X/G]$
of a smooth, projective variety $X$ 
by the action of a finite group $G$.

In this case, the inertia stack $\ccx$ and
double inertia stack $\cccx$ can also be written as global quotients $$\ccx =
\left[I_G(X)/G\right] = \left[\left(\coprod_{g\in G} X^g\right)/G\right],$$ and 
$$\cccx = \left[\II_G(X)/G\right] := \left[\left(\coprod_{g, h \in G} X^{g,h}\right)/G\right],$$
where each $X^g$ and $X^{g,h}$ is also smooth. As usual, there are morphisms $ev_i: \II_G(X) \rTo I_G(X)$ and $\check{ev}_i:
\II_G(X) \rTo I_G(X)$ for $i\in\{1,2,3\}$. 

Flat descent \cite[VIII\S1]{SGA6}
shows that $$K(\cx) \cong K_G(X).$$ The projection $$\pi: I_G(X) \rTo \ccx$$
induces a ring homomorphism (with respect to the usual product
$\otimes$) $$\pi^* : K(\ccx) = K_G(I_G(X)) \rTo K(I_G(X))^G.$$
Similarly, we also have a ring \emph{isomorphism} (with respect to the
usual product $\cup$) $$\pi^* : A^\bullet(\ccx) = A^\bullet_G(I_G(X)) \rTo
A^\bullet(I_G(X))^G$$ because $G$ is a finite group.

It is straightforward to see that
$\pi^*$ commutes with both pullback and pushforward along the morphisms $ev_i$ and $\check{ev}_i$ for all $i \in \{1,2,3\}$.

\begin{df}
For any stack $\cx$ which is a global quotient  $[X/G]$ of a smooth, projective 
variety $X$ by a finite group $G$, the \emph{small orbifold
K-theory} $\ok(\cx)$ of $\cx$ is the pre-Frobenius algebra $\ckb(X,G)
=\ck(X,G)^G$ of coinvariants of stringy K-theory:  
$$\ok(\cx):= \ckb(X,G).$$
This is linearly isomorphic to $K(\cx)$ and to $K(I_G(X))^G$, but with the stringy product instead of the tensor product.  
\end{df}

\begin{thm}\label{thm:pistuff}\label{thm:orbFACh}
Let $X$ be a smooth, projective variety with the action of a 
finite group $G$.
\begin{enumerate}
\item\label{piringhom} $\pi^*$ is a ring homomorphism from
$(\okf(\cx),\multipl)$ to $(\ok(\cx),\multipl)$. 
\item\label{piringhomchow} $\pi^*$ is an isomorphism of pre-Frobenius
algebras from
$(\oa(\cx),\multipl)$ to $( \ca^\bullet(X,G)^G,\multipl) 
= \cab(X,G)$. 
\item\label{chorb} The stringy Chern character $\CCh:
\ck(X,G) \rTo \ca^{\bullet}(X,G)$ induces an allometric isomorphism on the
invariants $$\och:\ok(\cx) =
\ckb(X,G)
\rTo 
\cab(X,G)
\cong
\oa(\cx)$$ which we call the \emph{small orbifold Chern character}. 

\item\label{pichern} The following diagram of ring homomorphisms (with respect to the orbifold, or stringy product $\multipl$) commutes: 
$$\begin{diagram}
\okf(\cx) &= &K_G(I_G X)	&\rTo^{\pi^*}      & K(I_G X)^G & = & \ok(\cx)\\
\dTo^{\CChf} &			     &&&& & \dTo^{\och}\\
 \oa(\cx) &=& A^{\bullet}_G(I_G X) &  \rEq^{\pi^*}     & \ca^\bullet(X,G)^G 
 & = &\cab(X,G).
\end{diagram}$$
\item\label{PresInd} The ring
$\ok(\cx)$ is independent of the choice of presentation of the stack $\cx$ as a
global quotient $[X/G]$.
\item\label{Embed} There is an embedding of pre-Frobenius algebras
\[
\ee:\ok([X/G])\otimes\nc\rTo\okf([X/G])\otimes\nc,
\]
such that $\pi^* \circ \ee = \vac_{\ok([X/G])}$.
\end{enumerate}
\end{thm} 

\begin{lm}\label{lm:stringypi} Let $\cso \in K(\ccx)$ be the (virtual) sheaf
given in Equation~(\ref{eq:orbS}), and let $\cs\in 
\ckb(X,G)$
be the
K-theoretic age sheaf given in Definition~\ref{df:kage}. 
Similarly, let $\ob$ be the obstruction bundle in $K(\cccx) = K_G(\II_G(X))$, and let $\crr$ be the obstruction bundle in $K(\II_G(X))$ arising in the stringy K-theory $\ck(X,G)$.  We have
\begin{equation}\label{eq:piCs}\cs = \pi^* \cso\end{equation} and 
\begin{equation}\label{eq:piOb}\crr = \pi^* \ob.\end{equation}
\end{lm}
\begin{proof}
Equation~(\ref{eq:piCs}) follows immediately from the definitions of $\pi^*,$ $\cs$, and $\cso$.  Equation~(\ref{eq:piOb}) follows from the fact that the $\cs_{m}$, $\cso_{m}$, and the normal bundles in the definition of $\crr$ and $\ob$ match term by term.
\end{proof}
\begin{proof}[Proof (of Theorem~\ref{thm:pistuff})]
Parts~(\ref{piringhom}) and ~(\ref{piringhomchow}) follow immediately from Lemma~\ref{lm:stringypi} and the fact that $\pi^*$ commutes with pullback and pushforward along the maps $ev_i$ and $\check{ev}_i$ for all $i\in \{1,2,3\}$. 

Part~(\ref{chorb}) follows from (taking invariants of)
Theorem~\ref{thm:ChernHomo}.

Part~(\ref{pichern}) follows from naturality of the classical Chern character
and Lemma~\ref{lm:stringypi}. 

Part~(\ref{PresInd}) follows from the fact that $\oa(\cx)$ is presentation
independent and $\och$ is an allometric isomorphism. 

Finally, Part~(\ref{Embed}) follows from the fact that if $Y$ is a smooth,
projective variety with the action of a finite group $G$, then there is a
canonical  isomorphism \cite{AS} (see also \cite[\S1]{Vi}, \cite[\S7]{VV} and 
\cite[\S3.4]{EG}) of algebras (with the ordinary multiplication $\otimes$):
\begin{equation}\label{eq:ASdecomp} 
K_{G}(Y) \otimes \nc \cong (K(I_G Y)\otimes\nc)^G 
\cong \bigoplus_{[\![m]\!]\in {\Gb}} (K(Y^m)\otimes \nc)^{Z_G(m)}.
\end{equation} 
The isomorphism takes $\cf \in K_G(Y)$ to an eigenbundle decomposition of $\cf|_{Y^m}$ in the $[\![m]\!]$-sector  $(K(Y^m)\otimes \nc)^{Z_G(m)}$

Setting $Y = I_G X$ then we obtain
\begin{align*}
\okf(\cx)\otimes \nc 		& := K([I_G(X)/G])\otimes \nc\\
& \cong (K(I_G(I_G(X))) \otimes \nc)^G\\ 
& = \bigoplus_{[\![g,m]\!]}(K(X^{\langle g,m\rangle }) \otimes\nc)^{Z_{G}(g,m)},  
\end{align*}
where the indices run over diagonal conjugacy classes of 
commuting
pairs $(g,m)$ in
$G^2$.  We denote the sector $(K(X^{\langle
g,m\rangle }) \otimes \nc)^{Z_{G}(g,m)}$ by $K^{\ccx}_{[\![g,m]\!]}$ 
Similarly, setting $Y=\II_G(X)$ we obtain
\begin{align*}
K(\cccx)\otimes \nc 	& = K([\II_G(X)/G])\otimes \nc\\
& \cong (K(I_G(\II_G(X))) \otimes \nc)^G\\
& = \bigoplus_{[\![g,h,m]\!]} (K(X^{\langle g,h,m\rangle }) \otimes
\nc)^{Z_{G}(g,h,m)}\\ & =: \bigoplus_{[\![g,h,m]\!]} K^{\cccx}_{[\![g,h,m]\!]},
\end{align*}
where the indices run over diagonal conjugacy classes of triples $(g,h,m)$ in
$G^3$
such that $m$ commutes with both $g$ and $h$.
It is easy to see that pullback along the morphism $ev_1:\cccx \rTo
\ccx$ takes the sector  
$K^{\ccx}_{[\![a,m]\!]}$ to the sum of sectors 
$\bigoplus_{[\![a,h,m]\!]} K^{\cccx}_{[\![a,h,m]\!]}$ where the indices run over conjugacy classes $[\![a,h,m]\!]$ with fixed pair $[\![a,m]\!]$.
Similarly, pullback along the morphism $ev_2:\cccx \rTo \ccx$ takes the sector 
$K^{\ccx}_{[\![a,m]\!]}$ to the sum of sectors 
$\bigoplus_{[\![g,a,m]\!]} K^{\cccx}_{[\![g,a,m]\!]}$ where the indices run over conjugacy classes $[\![g,a,m]\!]$ with fixed pair $[\![a,m]\!]$.
Pushforward along $\check{ev}_3$ maps sectors of the form $K^{\cccx}_{[\![a,b,m]\!]}$ to the sector $K^{\ccx}_{[\![ab,m]\!]}$.

Define the map $\ee:\ok(\cx)\otimes \nc \rTo \okf(\cx)\otimes \nc$ to send
the sector  $(K(X^g)\otimes \nc)^{Z_G(g)} \subseteq \ok(\cx)$ identically to
the ``untwisted'' sector $K^{\ccx}_{[\![g,1]\!]} \subseteq K(\ccx)\otimes \nc
= \okf(\cx)\otimes \nc$. Similarly, define a map $\eee: K(\II_G(X))^G =
\bigoplus_{[\![g,h]\!]}(K(X^{\langle g,h\rangle }) \otimes \nc)^{Z_G(g,h)} \otimes
\nc \rTo K(\cccx)\otimes \nc$ by taking the ${[\![g,h]\!]}$-sector
identically to the ${[\![g,h,1]\!]}$-sector.  

The map $\pi^*:\okf(\cx) \rTo \ok(\cx)$ sends a sector of the form $K^{\ccx}_{[\![g,m]\!]}$ to $(K(X^g)\otimes \nc)^{Z_G(g)}$, and on sectors of the form $K^{\ccx}_{[\![g,1]\!]}$ the map $\pi^*$ sends $K^{\ccx}_{[\![g,1]\!]} = $ identically to the sector $(K(X^g)\otimes \nc)^{Z_G(g)}$.
From this it is clear that $\pi^*\circ \ee = \vac$ and that $\ee$ is injective.

Clearly, $\ee$ and $\eee$ commute with $ev_i^*$ and $(\check{ev}_i)_*$, and $\lambda_{-1}$ preserves the decomposition.  Moreover, the tensor product of two homogeneous elements will vanish unless both elements lie in the same sector.  Thus, although $\eee(\crr)$ is \emph{not} equal to $\ob$ in $K(\cccx)$, it is true that for any $\cf$ and $\cf'$ in $\ok(\cx)$ we have
\begin{align*}
\ee(\cf)\multipl \ee(\cf') 	&= (\check{ev}_3)_* \left( ev_1^*(\ee(\cf)) \otimes ev_2^* (\ee(\cf')) \otimes \lambda_{-1}(\ob^*)\right)\\
  									&= (\check{ev}_3)_* \left( ev_1^*(\ee(\cf)) \otimes ev_2^* (\ee(\cf')) \otimes \eee(\lambda_{-1}(\crr*))\right)\\
									&= (\check{ev}_3)_* \left( \eee(ev_1^*(\cf)) \otimes \eee(ev_2^* (\cf')) \otimes \eee(\lambda_{-1}(\crr*))\right)\\
  									&= \ee\left((\check{ev}_3)_* \left( ev_1^*(\cf) \otimes ev_2^* (\cf') \otimes (\lambda_{-1}(\crr*))\right)\right)\\
									&= \ee(\cf\multipl\cf')
 \end{align*}
This shows that $\ee$ is a ring homomorphism.  Similar arguments show that $\ee$ preserves the trace and pairing.
 
\end{proof}
 
\subsection{Proof of Theorem~\ref{thm:ccx-preFA}}\label{subsec:proofs}
In this subsection we will prove Theorem~\ref{thm:ccx-preFA}.

The key step is showing that, as in the stringy case, the element $\ob\in K(\cccx)$
is represented by a vector bundle.  Associativity of the multiplication
follows from this fact.  To prove that $\ob$ is a vector bundle, we will show
that it is equal to the obstruction bundle for genus-zero, 
three-pointed orbifold stable maps into $\cx$. 

It is not hard to see that the stack $\M_{0,3}(\cx,0)$ of degree-zero, genus-zero, $3$-pointed orbifold stable maps into $\cx$ is naturally isomorphic to the double inertia stack $\cccx$.  We will equate the two from now on.   We denote the universal curve over it by $\varpi:\cc \to \M_{0,3}(\cx,0)$, and the universal stable map by $\bar{f}:\cc \to \cx$.  The evaluation maps from
$\M_{0,3}(\cx,0)$ to $\ccx$ are given by $ev_i([\bar{f}:\cc \to\cx]) = (f(p_i), [\![g_{p_i}]\!]) \in
\cx_{[\![g_{p_i}]\!]}$, where $p_i$ is the $i$-th marked (gerbe) point of
$\cc$, and $g_{p_i}$ is the image of the canonical generator of
$\stab(p_i)$ in $\stab(f(p_i))$.  Of course, this image is only defined up to conjugacy, since if $\cx$
is locally presented as $[X/G]$ near a point $p_i\in \cx$, then a
representative $\pt_i\in X$ of $p_i$ can be replaced by another
representative 
$\gamma\pt_i$ 
for any $\gamma\in G$, which replaces
$g_{p_i}$ by 
$\gamma g_{p_i} \gamma^{-1}$. Because of this, the 
$i$-th evaluation 
map $\M_{0,3}(\cx,0) \rTo \ccx$ agrees with the $ev_i$ described above for $\cccx \rTo \ccx$
for all $i \in \{1,2,3\}$. 

\begin{thm}\label{thm:orbObs}
  In the K-theory of $\cccx=\M_{0,3}(\cx,0)$, the following relation holds
  for the bundle $\ob$:
$$\ob\cong R^1\varpi_*(\bar{f}^*T\cx).$$ 
\end{thm}

\begin{proof}
  The idea of the proof is to use distinguished components of the
  stack of pointed admissible covers $\xi_{0,3}$ \cite[\S2.5.1]{JKK}
  to produce an \'etale cover of the moduli stack $\M_{0,3}(\cx,0)$.
  On this cover, we can easily produce an isomorphism of equivariant
  bundles, but it is not unique---it is only determined up to
  conjugacy.  However, the bundles we really want are the
invariant sub-bundles of these equivariant bundles, and on these subbundles the
 induced isomorphism is independent of conjugation.  Thus, \'etale
 descent applies, and we obtain the desired isomorphism.

  We first recall the definitions from \cite[\S2.5.1 and \S6]{JKK} of
  $\xi^G_{0,3}(\bm)$ and $\xi^G_{0,3}(X,\bm)$.  
As briefly described in Section~\ref{sec:FG=us}, to each triple $\bm = (m_1,
m_2, m_3) $ with $m_1m_2m_3=1$, there is a canonical choice of pointed admissible $\langle m \rangle$-cover $(E,\pt_1,\pt_2,\pt_3)\rTo(\mathbb{P}^1,0,1,\infty)$,
ramified only over the points $0$, $1$, and $\infty$, and with
monodromy $m_1$, $m_2$, and $m_3$, respectively, at those points.  
There is also a canonical choice of pointed admissible $G$-cover $(\tilde{E},\pt_1,\pt_2,\pt_3)\rTo(\mathbb{P}^1,0,1,\infty)$ with $\tilde{E} = E \times_{\langle m \rangle} G$.
For a complete discussion of these constructions, see \cite[\S2.5.1]{JKK}.  

\begin{df}
  We define $\xi^{\langle \bm \rangle }_{0,3}$ to be the connected
  component of the stack of $3$-pointed admissible
$\langle\bm\rangle$-covers of genus zero that contains the canonical
  admissible cover $(E,\pt_1,\pt_2,\pt_3)$.  Similarly, we define
  $\xi^G_{0,3}$ to be the connected component of the stack of
  three-pointed admissible $\langle \bm\rangle$-covers of genus zero that
  contains the admissible cover $(\tilde{E},\pt_1,\pt_2,\pt_3)$.
  
  If $X$ is any variety with a $G$-action, a degree-zero, $3$-pointed
  $G$-stable map of genus zero is a $G$-equivariant morphism
  $\tilde{E} \to X$ from a $3$-pointed admissible $G$-cover to $X$,
  such that the induced morphism $\tilde{E}/G \to X/G$ is a
  $3$-pointed stable map of genus zero.  We define
  $\xi^G_{0,3}(X,0,\bm)$ to be the component of the stack of pointed
  $G$-stable maps whose underlying $3$-pointed admissible $G$-covers
  $\tilde{E}$ correspond to points of $\xi^G_{0,3}(\bm)$.
\end{df}

It is easy to see that there is a canonical isomorphism
$\xi^{\langle \bm \rangle }_{0,3}(\bm) \irightarrow \xi^G_{0,3}(\bm)$, and
that $\xi^G_{0,3}(\bm)$ is the stack quotient 
$\cb H_\bm = [pt/H_\bm]$ of
a point modulo the group $H_\bm:=\langle m_1\rangle \cap \langle
m_2\rangle \cap \langle m_3\rangle$ (see \cite[Prop 2.20]{JKK}).
Moreover, in \cite[Lemma 6.7]{JKK} it is shown that
$\xi^G_{0,3}(X,0,\bm)$ is canonically isomorphic to $\xi^G_{0,3}
\times X^\bm$.  Finally, we have a morphism $q:
\xi_{0,3}(X,0,\bm) \to \M_{0,3}([X/G],0)$ given by sending a
pointed $G$-stable map $[f:E \to X]$ to the induced map of quotient stacks
$[\bar{f}:[E/G] \to [X/G]\!]$.  This morphism is easily seen to be \'etale.

Now we may begin the proof.  First consider any \'etale
cover $U$ of $\cx$, consisting of a disjoint union of smooth varieties
$X_{\alpha}$ with finite groups $G_{\alpha}$ acting to make
$q_{\alpha}:X_{\alpha}\to \cx$ induce an isomorphism
$[X_{\alpha}/G_{\alpha}]$ to a neighborhood in $\cx$ (that is,
$\{(X_{\alpha}, G_{\alpha}, q_{\alpha})\}$ form a uniformizing system).
We may construct an \'etale cover
$$p:\coprod_{\alpha, \bm} X_{\alpha}^{\bm} \to \coprod_{\alpha,
  \bm}\xi^{G_\alpha}_{0,3} \times X_{\alpha}^{\bm} = \coprod_{\alpha,
  \bm}\xi^{G_\alpha}_{0,3}(X_{\alpha},0,\bm) \to \M_{0,3}(\cx,0),$$
where for each $\alpha$, the $\bm$ runs through all triples in
$G_\alpha$ whose product is $1$, and the first morphism is induced by
the obvious (\'etale) map $pt \times X_{\alpha}^{\bm} \to [pt/H_\bm]
\times X_{\alpha}^{\bm} = \xi^{G_\alpha}_{0,3} \times
X_{\alpha}^{\bm}$.

For each $\alpha$ and $\bm$, the pullback $p^*\ob$ is easily seen to
be the usual obstruction bundle $\crr(\bm)$ on $X_{\alpha}^{\bm}$, and
the pullback $p^* \left(T\M_{0,3}(\cx,0) \ominus J^*T\cx \oplus
  \bigoplus^3_{i=1} ev^*_i \cs\right)$ is clearly equal to $
TX_{\alpha}^\bm\ominus
\res{TX_{\alpha}}{X_{\alpha}^\bm}\oplus\bigoplus_{i=1}^3
\res{\cs_{m_i}}{X_{\alpha}^{\bm}}$.  But to prove the theorem, we will
need to provide a canonical isomorphism between these bundles.

The fibered product $X_{\alpha}^{\bm}\times_{\M_{0,3}(\cx,0)}
X_{\beta}^{\bm'}$ is non-empty if and only if
the stack $\xi^{G_\alpha}_{0,3}(X_{\alpha},0,\bm)
\times_{\M_{0,3}(\cx,0)} \xi^{G_\beta}_{0,3}(X_{\beta},0,\bm')$ is
non-empty; and it is straightforward to see that this occurs only if
there is a group $G$ with injective homomorphisms $G\hookrightarrow
G_\alpha$ and $G\hookrightarrow G_\beta$, such that the triple $\bm'$
is diagonally (i.e., all three terms simultaneously) conjugate to
$\bm$ in $G^3$.  Moreover, for each connected component of
$\M_{0,3}(\cx,0)$, there is a well-defined diagonal conjugacy class of
such triples.

For each such conjugacy class, choose a representative $\bm$ and let
$K = \langle m_1, m_2,m_3\rangle$ be the group generated by the
triple.  As described above, this triple determines a well-defined
distinguished component $\xi^K_{0,3}(\bm)$ of the stack of
three-pointed, admissible $K$-covers of genus zero.

Choose, once and for all, an isomorphism $\Phi_{\bm}$ of
$K$-representations giving the (virtual) equality of
Equation~(\ref{eq:magic2}) in Lemma~\ref{thm:rep-magic},
but where $G$ is replaced by $K$.  For any
other triple $\bm'$ in the same conjugacy class, there is a canonical
isomorphism of groups $K' =\langle m'_1, m'_2,m'_3 \rangle
\irightarrow K$ taking $\bm'$ to $\bm$, and a canonical (equivariant)
isomorphism of representations $H^1(E'; \co_{E'}) \cong H^1(E;
\co_{E})$, where $E \to C \to \xi^K_{0,3}(\bm)$ is the three-pointed
admissible $K$-cover with holonomy $\bm$, and $E'\to C' \to
\xi^{K'}_{0,3}(\bm')$ is the three-pointed admissible $K'$-cover with
holonomy $\bm'$.  Similarly, we have canonical (equivariant)
isomorphisms of the representations
\begin{equation}
\nc \ominus \nc[K]\oplus\bigoplus_{i=1}^n\bigoplus_{k_i=0}^{r_i-1}
\frac{k_i}{r_i} \ind^K_{\langle m_i\rangle }\bv_{m_i,k_i}
\cong
\nc \ominus \nc[K']\oplus\bigoplus_{i=1}^n\bigoplus_{k_i=0}^{r_i-1}
\frac{k_i}{r_i} \ind^{K'}_{\langle m_i\rangle }\bv_{m_i,k_i}.
\end{equation}  Thus $\Phi_{\bm}$ induces an isomorphism $\Phi_{\bm'}$ 
for each triple $\bm'$ which is conjugate to $\bm$.

If $G$ is any group containing both $K$ and $K'$, with $K'$ a
conjugate (say by $\gamma \in G$) of $K$, then letting $\tilde{E}\to C
\to \xi^G_{0,3}(\bm)$ denote the distinguished three-pointed $G$-cover
with holonomy $\bm$, and $\tilde{E}'\to C \to \xi^G_{0,3}(\bm')$
denote the distinguished universal three-pointed $G$-cover with
holonomy $\bm'$, the group action $\rho(\gamma):\xi^G_{0,3}(\bm)
\irightarrow \xi^G_{0,3}(\bm')$ identifies the base ($\gamma$ acts on
$E$ and $E'$).  Furthermore, we have canonical isomorphisms of
$G$-representations
\begin{equation}\label{eq:indH}
H^1(\tilde{E};\co_{\tilde{E}}) \cong \ind^G_K (H^1 (E; \co_E))
\end{equation}
and 
\begin{equation}\label{eq:indHprime}
H^1(\tilde{E}';\co_{\tilde{E}'}) \cong \ind^G_{K'} (H^1 (E'; \co_{E'})).
\end{equation}
As $G$-representations, $H^1(\tilde{E};\co_{\tilde{E}})$ and
$\rho(\gamma)^*H^1(\tilde{E}';\co_{\tilde{E}'})$ are not identical,
but rather are conjugate; that is, $H^1(\tilde{E}';\co_{\tilde{E}'})$
is the representation of $G$ arising from conjugating the action of
$G$ on $H^1(\tilde{E};\co_{\tilde{E}})$ by $\gamma$.  The same holds
for the induced representations
\begin{equation}\label{eq:indMagic}
\nc[G/K]
\ominus \nc[G]\bigoplus^3_{i=1}\bigoplus^{r_i-1}_{k_i=0}
\frac{k_i}{r_i} \ind^G_{\langle m_i \rangle}\bv_{m_i,k_i}  
\end{equation}
$$\cong \ind^G_K(\nc \ominus \nc[K] \bigoplus^3_{i=1}
\bigoplus^{r_i-1}_{k_i=0} \frac{k_i}{r_i} \ind^K_{\langle m_i
\rangle}\bv_{m_i,k_i}),$$ and 
\begin{equation}
\nc[G/K'] \ominus \nc[G]\bigoplus^3_{i=1}\bigoplus^{r_i-1}_{k_i=0}
\frac{k_i}{r_i} \ind^G_{\langle m'_i \rangle}\bv_{m'_i,k_i}  
\end{equation}
$$\cong \ind^G_{K'}(\nc \ominus \nc[K'] \bigoplus^3_{i=1}
\bigoplus^{r_i-1}_{k_i=0} \frac{k_i}{r_i} \ind^{K'}_{\langle m'_i \rangle}
\bv_{m'_i,k_i}).$$ 

Finally, for an open subset $V$ of any $X_{\alpha}$ with $G$ acting on
$V$, pulling back by the action $$\rho(\gamma) : V^{\bm} \irightarrow
V^{\bm'}$$ makes the $G$-bundle $\rho(\gamma)^*f^* TV =
\co_{\tilde{E}'} \boxtimes \ind^G_{K'}TV|_{V^{\bm'}}$ on $V^{\bm}$
isomorphic to the conjugate by $\gamma$ of the $G$-bundle $f^* TV =
\co_{\tilde{E}} \boxtimes \ind^G_{K}TV|_{V^{\bm}}$.  Thus the
isomorphisms $\Phi_{\bm'}$ and the induced isomorphisms
\begin{align}
  \tilde{\Phi}_{\bm}: R^1\pi_*(f^*TV) & =
  H^1(\tilde{E};\co_{\tilde{E}}) \otimes \ind^G_{K}TV|_{V^{\bm}}
  \notag \\
  & \irightarrow 
\left(\nc[ G/K] \ominus \nc[G] \oplus
    \bigoplus_{i=1}^3\bigoplus_{k_i = 0}^{r_i-1} \frac{k_i}{r_i}
    \ind_{\langle m_i\rangle }^{G}\bv_{m_i,k_i} \right)\tensor
  \ind^G_{K} TV|_{V^{\bm}}
\end{align}
on $\xi^G_{0,3}(\bm) \times V^{\bm}$ are determined up to conjugacy by
an element in $G$.

However, a representation and any conjugate of that representation
have canonically identified 
invariants, so the isomorphisms
$\tilde{\Phi}_{\bm}$ induce isomorphisms of the 
invariant bundles
\begin{align}
  \bar{\Phi}_{\bm}: R^1\pi^G_*(f^*TV) & \irightarrow
  \left(\left(\nc[G/K] \ominus \nc[G] \oplus
      \bigoplus_{i=1}^3\bigoplus_{k_i = 0}^{r_i-1} \frac{k_i}{r_i}
      \ind_{\langle m_i\rangle }^{G}\bv_{m_i,k_i}
    \right) \tensor TV|_{V^{\bm}}\right)^G \notag\\
  &= TV^{\bm} \ominus TV|_{V^{\bm}} \oplus \bigoplus_{i=1}^3
  \cs_{m_i}|_{V^{\bm}},
\end{align}
which are independent of conjugation.

In summary, we have chosen an explicit isomorphism
$$\overline{\Phi}: p^* \ob \irightarrow p^*\left(T \M_{0,3}(\cx,0)
  \ominus J^* T\cx \oplus \bigoplus^3_{i=1} ev^*\cs\right)$$
on the \'etale cover $\coprod_{\alpha, \bm} X_{\alpha}^{\bm} \rTo^{p}
\M_{0,3}(\cx,0),$ with the particular property that on the product
$$\coprod X_{\alpha}^{\bm} \times_{\M_{0,3}(\cx,0)} \coprod
X_{\alpha}^{\bm} \pile{\rTo^s \\ \rTo_{t}} \coprod X_{\alpha}^{\bm}$$
we have $s^*\overline{\Phi} = t^*\overline{\Phi}$.  Thus by \'etale
descent the isomorphism $\overline{\Phi}$ descends from the cover
$\coprod X_{\alpha}^\bm$ 
to the stack $\M_{0,3}(\cx,0)$.
\end{proof}

The proof of associativity given in
Lemma~\ref{thm:SecondAssoc} is now easily adapted to give a proof of
associativity for the orbifold product in 
$\oa(\cx)$ and 
$\okf(\cx)$.  The rest of the properties of a pre-Frobenius algebra are straightforward to check.

The fact that $\oa(\cx)$ is isomorphic to the construction of \cite{AGV} also follows from Theorem~\ref{thm:orbObs} and from the equality $\M_{0,3}(\cx,0) = \cccx$, since the definition of the product in \cite{AGV} is the usual quantum product with the obstruction bundle $R^1\varpi_*(\bar{f}^*T\cx)$.

\section{Stringy topological K-theory and stringy cohomology}\label{sec:top}

All of the results in the previous sections have their counterparts in
the topological category.

\subsection{Ordinary topological K-theory and cohomology}
Throughout this section, unless otherwise stated, $G$ is a finite
group acting on a compact, almost complex manifold $X$, preserving
the almost complex structure.

Furthermore, let $H^\bullet(X)$ be the rational cohomology of $X$. It
is a Frobenius superalgebra: a Frobenius algebra with a multiplication
that is graded commutative.

Topological K-theory $\kt(X):= \kt(X;\nz)\tensor_\nz \nq$ is also a
Frobenius superalgebra with the $\nz/2\nz$-grading: $$\kt(X) =
\ktzero(X;\nz)\tensor_\nz \nq \oplus \ktone(X;\nz)\tensor_\nz \nq.$$
Here $\ktzero(X;\nz)$ is
defined exactly as $K(X;\nz)$ but in the topological category. That
is, $\ktzero(X;\nz)$ is additively generated by isomorphism classes of
complex topological vector bundles over $X$ modulo the relation of
Equation~(\ref{eq:KRelation}) whenever Equation (\ref{eq:ModExact})
holds.  The odd part $\ktone(X;\nz)$ is defined to be $\ktzero(X\times
\nr;\nz)$.  Equivalently, we may take $\ktone(X;\nz)$ to be the kernel
of the restriction map $i^*:\ktzero(X\times S^1) \to \ktzero(X \times
pt)$ induced in K-theory from the inclusion of a point $i:X\times pt
\to X\times S^1$.  

Associated to a differentiable proper map of almost complex manifolds
$f:X\to Y$, there is the induced pushforward morphism $f_*:\kt(X)\to
\kt(Y)$ (see \cite[IV 5.24]{Kar} and \cite[Sec. 4]{AH}).  In
particular, if $Y$ is a point and $f:X\to Y$ is the obvious map, we
again define the Euler characteristic as 
$$\chi(X, \cf) :=f_* \cf.$$  Associated to any continuous $f:X\to Y$, there is a pullback
homomorphism $f^*:\kt(Y)\to \kt(X)$ \cite[II.1.12]{Kar}.

For any compact, almost complex manifolds $X$ and $Y$, there are
natural morphisms $$ \nu:\ktn(X)\otimes \ktm(Y) \to \ktzero(X\times Y
\times \nr^{n+m}).
$$
Bott periodicity says that if $Y$ is a point, there is an
isomorphism $$\beta:\ktzero(X) \irightarrow \ktzero(X\times \nr^2)$$
\cite[III.1.3]{Kar}, which is natural with respect to both pullback
and pushforward.  Therefore, for any compact, almost complex manifold
$X$, composition of $\nu$ with pullback along the diagonal map
$\Delta:X \to X\times X$ gives a multiplication $$\mu: \ktn(X)\otimes
\ktm(X) \to \ktzero(X \times \nr^{n+m}) \subseteq \kt(X)$$
if $n+m
\leq 1$, and $$\mu:\ktone(X) \otimes \ktone (X) \rTo
\ktzero(X\times\nr^2) \rTo^{\beta^{-1}} \kt(X),$$
if $n=m=1$. Here
$\beta^{-1}$ is the inverse of the Bott isomorphism.  We will write
$\cf_1 \otimes \cf_2$ to denote $\mu(\cf_1,\cf_2)$.  This product
makes $\kt(X)$ into a commutative, associative superalgebra
\cite[II.5.1 and II.5.27]{Kar}.

One can define a metric on $\kt(X)$ by $$\eta_\kt(\cf_1,\cf_2) :=
\chi(X,\cf_1\otimes\cf_2),$$ and we define $\vac :=\co_X$.  It is
straightforward to check that $(\kt(X),\otimes,\vac,\eta_\kt)$ is a
Frobenius superalgebra. Moreover, the projection formula holds for
proper, differentiable maps with a compact target \cite[IV.5.24]{Kar}.

The Frobenius superalgebra of topological K-theory satisfies the usual
naturality properties with respect to pullback, is also a
$\lambda$-ring \cite[\S 7.2]{Kar}, and satisfies the splitting
principle \cite[Thm IV.2.15]{Kar}.

For all $i$, the $i$-th Chern class $c_i(\cf)$ associated to any $\cf$
in $\ktzero(X)$ belongs to $H^{2i}(X)$, and so $H^{2 p}(X)$ may be
regarded as the analogue of the Chow group $A^p(X)$. The associated
Chern polynomial $c_t$ satisfies the usual multiplicativity and
naturality properties, and the Chern character $\Ch:\kt(X)\to
H^\bullet (X)$, defined by Equation~(\ref{eq:ChernCharDef}), is an
isomorphism of commutative, associative superalgebras \cite[Thm.
V.3.25]{Kar}. The Todd classes are defined from the ordinary Chern
classes as before. In addition, Proposition (\ref{prop:ToddCh}) holds
in topological K-theory since it follows from the splitting principle,
the Chern character isomorphism, and the $\lambda$-ring properties
\cite[Prop. I.5.3]{FH}.

Finally, the Grothendieck-Riemann-Roch formula (see \cite[Cor
V.4.18]{Kar} or \cite[Thm. 4.1]{AH}) and the excess intersection
formula (Theorem~\ref{thm:excessintersection}) hold (see \cite[Prop
3.3]{Qu}, which is written for cobordism, but the proof works as well
for topological $K$-theory).

\begin{rem}
  Let $X$ be a compact $G$-manifold with a smoothly varying 
  one-parameter family of $G$-equivariant, almost complex structures
  $J_t:TX\to TX$ for all $t$, say, in the interval $[0,1]$. Because of the
  homotopy invariance of characteristic classes, the resulting
 $G$-Frobenius algebras $\ch(X;G)$ and $\ck(X;G)$, and the stringy Chern
  character are all independent of $t$. Therefore, these stringy algebraic
  structures depend only upon the homotopy class of the
  $G$-equivariant almost complex structure on the $G$-manifold $X$.
  
  In particular, when $X$ is a compact symplectic manifold with an
  action of $G$ preserving the symplectic structure, 
  since up to homotopy there exists a unique, $G$-equivariant, almost complex
  structure compatible with the symplectic form \cite[Ex. D.12]{GGK},
  these stringy algebraic structures are invariants of the symplectic
  manifold with $G$-action.
\end{rem}

\begin{rem}\label{rem:StableACM}
  While we are primarily interested in $G$-equivariant almost complex
  manifolds in this section, our constructions generalize in a
  straightforward way to the case where $X$ is a compact manifold with
  an oriented, $G$-equivariant stable complex structure (see
  \cite[App. D]{GGK}).  The key point \cite{GHK} is that 
for any subgroup $H\le G$, a  
  $G$-equivariant stable complex structure induces an almost complex
  structure on the normal bundle to the submanifold $X^H$ (the locus
  of points fixed by $H$).  Furthermore,
  both $\cs_m$ (see Remark \ref{rem:TangentNormal}) and the right hand
  side of Equation (\ref{eq:Magic}) only depend upon such normal
  bundles.
\end{rem}

\subsection{Stringy topological K-theory and stringy cohomology} 

Let $X$ be a compact, almost complex manifold with an action of a
finite group $\rho:G\to \Aut(X)$ preserving the almost complex
structure.

Fantechi and G\"ottsche's \cite{FG} stringy cohomology $\ch(X,G)$ of
$X$ is given by
$$\ch(X,G) := \bigoplus_{m\in G} \ch_m(X),$$
where $\ch_m(X) :=
H^\bullet(X^m)$, and the definition of the multiplication is still
given by Equation~(\ref{eq:DefChowMult}), the trace element $\tau$ by
Equation~(\ref{eq:taua}), and similarly for the metric
and unity.
However, the $\nq$-grading here is not quite
that defined by Equation (\ref{eq:DefGrading}), but is defined instead
by the equation
\begin{equation}\label{eq:TopGrading}
\grq{v_m} := 2 a(v_m) + \grz{v_m},
\end{equation}
where $\grz{v_m} := p$ when $v_m$ belongs to $H^p(X^m)$ and $a(v_m) := a(m,U)$.

Furthermore, Theorem~(\ref{thm:ChowGFA}) holds, provided that
$\ca(X,G)$ is everywhere replaced by $\ch(X,G)$, $\dim X$ is
understood to be the dimension of $X$ as a real manifold, and
           ``pre-$G$-Frobenius algebra'' is replaced by ``$G$-Frobenius
superalgebra.'' 
However, we need to establish the following Proposition to complete the
proof. 
\begin{prop}\label{prop:TracesAgree}
Let $X$ be a compact, almost complex manifold with the action of a finite
group $G$ with stringy cohomology $\ch(X,G)$. If the trace element $\tau$ 
is given by Equation~(\ref{eq:taua}), then Equation~(\ref{eq:FAtrace}) holds, where $\ch
:= \ch(X,G)$. Consequently, the trace axiom,
Equation~(\ref{eq:OldTraceAxiom}), is satisfied, and $\ch(X,G)$ is a
$G$-Frobenius superalgebra. The characteristic $\ttau(\vac)$ of $\ch(X,G)$ 
also satisfies Equation~(\ref{eq:Characteristic}) and is an integer.
\end{prop}
\begin{proof}
To avoid distracting signs, we assume that $\ch(X,G)$ has only even
dimensional cohomology classes.

We first
prove Equation~(\ref{eq:FAtrace}). We henceforth adopt the
notation from Section~(\ref{subsec:TraceAxiom}). 
Let $\{ \nu_{\alpha[a]} \}$
be a homogeneous basis for $\ch_a(X)$ with 
$\alpha[a] \in \{1,\ldots,d_a\}$, 
where
$d_a$ is the dimension of $\ch_a(X)$.
Similarly, let $\{\mu_{\beta[a^{-1}]}\}$ with 
$\beta[a] \in \{1,\ldots,d_a\}$ be a homogeneous basis for 
$\ch_{a^{-1}}(X)$. 
Let $\eta_{\alpha[a],\beta[a^{-1}]}$ be the matrix of the metric pairing
$\ch_a(X)\times\ch_{a^{-1}}(X) \rTo \nq$, with respect to these bases, 
and let $\eta^{\alpha[a] \beta[a^{-1}]}$ 
be the inverse of 
$\eta_{\alpha[a],\beta[a^{-1}]}$.
We observe, from the K\"unneth theorem, that 
\begin{equation}\label{eq:twisteddiagonalp}
\Deltap_2^*\Deltap_{1 *}\vac_{X^a} = \eta^{\alpha[a] \beta[a^{-1}]}
\res{(\rho(b) \nu_{\alpha[a]})}{X^{H'}} \tensor
\res{\nu_{\beta[a^{-1}]}}{X^{H'}}.
\end{equation}

Thus, we have 
for $m_1 = [a,b]$,
\begin{eqnarray*}
&& \tr_{\ch_a(X)}(L_{v_{m_1}}\circ\rho(b))  
\\
&=&
\int_{X^{H'}}\ctop(\crr(\bmp)) \cup \e_{m_1}^*
v_{m_1}\cup
\res{(\rho(b) \nu_{\alpha[a]})}{X^{H'}} \cup
\res{\nu_{\alpha[a^{-1}]}}{X^{H'}} \eta^{\alpha[a]\alpha[a^{-1}]} \\ 
&=&
\int_{X^{H'}}\ctop(\crr(\bmp)) \cup \e_{m_1}^*
v_{m_1}\cup \Deltap_2^*\Deltap_{1 *}\vac_{X^a} \\ 
&=&
\int_{X^{H'}}\ctop(\crr(\bmp)) \cup \e_{m_1}^*
v_{m_1}\cup \jp_{2   *} \ctop(\ce')
 \\ 
&=&
\int_{X^{H'}} \jp_{2 *}(\jp_2^*\e_{m_1}^*v_{m_1})\cup
\jp_2^*\ctop(\crr(\bmp)) \cup
\ctop(\ce')\\
&=&
\int_{X^{H'}} \jp_{2 *}(\jp_2^*\e_{m_1}^*v_{m_1})\cup
\ctop(\jp_2^*\crr(\bmp)\oplus \ce')\\
&=&
\int_{X^{H}} \res{v_{m_1}}{X^{H}}
\cup \ctop(\jp_2^*\crr(\bmp)\oplus\cep) \\
&=&
\int_{X^{H}} \res{v_{m_1}}{X^{H}}
\cup \ctop(TX^H\oplus \res{\cs_{m_1}}{X^H}) \\
&=&
\tau_{a,b}(v_{m_1}),
\end{eqnarray*}
where the first equality holds by definition of the trace, the second by
Equation~(\ref{eq:twisteddiagonalp}), the third by
Theorem~(\ref{thm:excessintersection}), the fourth by the projection
formula, the sixth by properties of the top Chern class, and the seventh
by Theorem~(\ref{thm:GenusOne}).

The trace axiom (Equation~(\ref{eq:OldTraceAxiom})) for the $G$-Frobenius algebra $\ch(X,G)$ follows from the trace axiom for a
pre-$G$-Frobenius algebra together with Equation~(\ref{eq:FAtrace}). 

The integrality of the characteristic follows from
Equation~(\ref{eq:Canonicaltaub}) by plugging in $\overline{v} = \vac$
for the Frobenius superalgebra of $G$-coinvariants $\chb = \ch(X,G)^G$.
\end{proof}

Stringy topological K-theory $\ckt(X,G) := \bigoplus_{m\in G}
\ckt_m(X)$ is defined additively by $\ckt_m(X) := \kt(X^m)$ for all
$m$ in $G$. The stringy multiplication, metric, 
the trace element $\tau$ and unity
are defined just as in the case of $\ck(X,G)$. 
This is compatible with
the $\nz/2\nz$-grading because 
$\crr$ is an
element of $\ktzero(X^{\bm})$.

Since the Eichler trace formula holds for all compact Riemann
surfaces, our formula (\ref{eq:Magic}) for the obstruction bundle, and
indeed the entire analysis in 
Sections~\ref{sec:AssocTrace}  and~\ref{sec:obs},
holds in topological K-theory. The K-theoretic
version of Proposition~\ref{prop:TracesAgree} holds as the arguments are
purely functorial. 
Now an argument essentially identical to that for stringy cohomology shows that  
$((\ckt(X,G),\rho),\multipl,\vac,\eta)$ is a
$G$-Frobenius superalgebra. 
We state this formally in the following proposition.
\begin{prop}
Let $X$ be a compact, almost complex manifold with the action of a finite
group $G$ with stringy topological K-theory $\ckt(X,G)$. If the trace
element $\tauk$ is given by Equation~(\ref{eq:tauk}), then Equation~(\ref{eq:FAtrace})
holds, where $\ch := \ckt(X,G)$. Consequently, the trace
axiom, Equation~(\ref{eq:OldTraceAxiom}), is satisfied, and $\ckt(X,G)$ is a
$G$-Frobenius superalgebra. 
\end{prop}

Furthermore, the stringy Chern character 
$\CCh:\ckt(X,G)\to\ch(X,G)$ is still defined by Equation~(\ref{eq:DefCCh}). The
rest of the analysis in Section~\ref{sec:chern} holds, provided that
$\ca(X,G)$ is everywhere replaced by $\ch(X,G)$ and $K$-theory is
everywhere replaced by topological K-theory. Therefore,
$\CCh:\ckt(X,G)\to 
\ch(X,G)
$ is an allometric isomorphism.

Finally, the analysis in Section~\ref{sec:stack} holds after replacing
Chow groups by cohomology everywhere. 
In particular, since $\chb(X,G)$ is isomorphic \cite{FG} to the
stringy (or Chen-Ruan orbifold) cohomology $H^\bullet_{orb}([X/G])$, the stringy Chern character
$\CChb:\cktb(X,G)\to\chb(X,G)$ gives a ring isomorphism 
$\och:\okt([X/G])\to H^\bullet_\mathrm{orb}([X/G])$, where $\okt([X/G])$ is the topological small orbifold K-theory of $[X/G]$.

\subsection{The symmetric product and crepant
  resolutions}\label{eq:SubsectionSP}

One of the most interesting examples of stringy K-theory and
cohomology is the symmetric product. Let $X := Y^n$, where $Y$ is 
an almost
complex manifold of complex dimension $d$ with the symmetric group
$S_n$ acting on $Y^n$ by permuting its factors. In this case, for any
$m \in S_n $ it is easy to see that the age $a(m)$ is related to the
length of the permutation $l(m)$:
$$a(m)= l(m) d /2.$$
Consequently, by Equation~(\ref{eq:TopGrading}),
the $\nq$-grading on $\ch(X,G)$ is, in fact, a grading by (possibly
odd) integers.

Consider stringy topological K-theory $\ckt(Y^n,S_n)$ of the
$S_n$-variety $Y^n$. Choose the $2$-cocycle (discrete torsion)
$\alpha$ in $Z^2(S_n,\field^*)$
\[
\alpha(m_1,m_2) := (-1)^{\varepsilon(m_1,m_2)},
\]
where $\varepsilon$ is defined by
\[
\varepsilon(m_1,m_2) := \frac{1}{2} (l(m_1)+l(m_2) - l(m_1 m_2)).
\]
It is straightforward to verify that $\varepsilon(m_1,m_2)$ is an
integer.  Now, twist the $S_n$-Frobenius algebra $\ckt(Y^n,S_n)$ by
$\alpha$, as in Section~\ref{sec:DiscTors}, to yield a new
$S_n$-Frobenius algebra
$((\ckt(Y^n,S_n),\rho),\star,\vac,\eta^\alpha)$, which we will denote
by $\ctkt(Y^n,S_n)$.  Notice that the $G$-action is unchanged by the
twist, but the twisted multiplication $\ckt(Y^n,S_n)$ is given by the
formula
\begin{equation}\label{eq:SPDiscreteTorsion}
v_{m_1} \star v_{m_2} :=  
v_{m_1} \multipl^{\alpha} v_{m_2}
= \alpha(m_1,m_2) v_{m_1} \multipl v_{m_2},
\end{equation}
where $\multipl$ denotes the stringy multiplication in $\ckt(Y^n,S_n)$.

Twisting the multiplication on the stringy cohomology of $Y^n$ in the
same fashion, we obtain the $S_n$-Frobenius algebra
$((\ch(Y^n,S_n),\rho),\star,\vac,\eta^\alpha)$, which we will denote
by $\cth(Y^n,S_n)$.  By the obvious topological analogue of
Corollary~\ref{crl:twistCh}, the stringy Chern character
$\CCh:\ctkt(Y^n,S_n)\to \cth(Y^n,S_n)$ is an isomorphism of
$S_n$-commutative algebras.  After taking $S_n$-coinvariants, we
obtain a ring isomorphism
\[
\och :\mathbf{K}^\mathrm{top}_\mathrm{orb}([Y^n/S_n])\to
      \mathbf{H}^\bullet_\mathrm{orb}([Y^n/S_n]),
\]
where $\mathbf{K}^\mathrm{top}_\mathrm{orb}([Y^n/S_n])$ is
the topological small orbifold K-theory
$\okt([Y^n/G])$, but with the $\alpha$-twisted multiplication, and
similarly for $\mathbf{H}^\bullet_\mathrm{orb}([Y^n/S_n])$.

What makes these particular twisted rings interesting is the following
theorem.
\begin{thm}\label{thm:SymmetricProduct} Let $Y$ be a complex, projective
  surface such that $c_1(Y) = 0$. Consider $Y^n$ with $S_n$ acting by
  permutation of its factors. If $Y^{[n]}$ denotes the Hilbert scheme
  of $n$ points in $Y$, then
  $\mathbf{K}^\mathrm{top}_\mathrm{orb}([Y^n/S_n])$
is isomorphic as a Frobenius superalgebra to $\kt(Y^{[n]})$.
\end{thm}
\begin{proof}
We define $\psi$ so that the following diagram commutes
\begin{equation}
\begin{diagram}
\mathbf{K}^\mathrm{top}_\mathrm{orb}([Y^n/S_n])  & \rTo^{\och}
&
\mathbf{H}^\bullet_\mathrm{orb}([Y^n/S_n])  \\
\dTo^{\psi} & & \dTo^{\psi'}\\
\kt(Y^{[n]})  & \rTo^{\cch} & H^\bullet(Y^{[n]})\\
\end{diagram},
\end{equation}
where $\psi'$ is the ring isomorphism $\Psi^{-1}$ in \cite[Thm 3.10]{FG}. 
This uniquely defines $\psi$, since $\cch$ and $\och$
are ring isomorphisms. 

The homomorphism $\psi$ also preserves the metrics because of the
Hirzebruch-Riemann-Roch Theorem and the fact that $\psi'$ preserves
the metrics.
\end{proof}

\begin{rem}
The rings
$\mathbf{K}^\mathrm{top}_\mathrm{orb}([Y^n/S_n])\otimes_\nq\nc$ and
$\okt([Y^n/S_n])\otimes_\nq \nc $ are isomorphic (see \cite{Ru}).
Since $Y^{[n]} \to Y^n/S_n$ in the previous theorem is a crepant (and
hyper-K\"ahler) resolution, this is an example of our K-theoretic
version of Conjecture~\ref{conj:KHRC}. Our result is nontrivial
precisely because of the nontrivial definition of multiplication on
$\okt([Y^n/S_n])$ and the stringy Chern character.
\end{rem}

\bibliographystyle{amsplain}

\providecommand{\bysame}{\leavevmode\hboxto3em{\hrulefill}\thinspace}

\end{document}